\documentclass[11pt,leqno]{article}

\usepackage{amsmath,amsthm}
\usepackage {latexsym}
\usepackage{amssymb}

\def\th{\theta}
\newcommand{\les}{\lesssim}

\newcommand{\bse}{\begin{equation}}

\newcommand{\bea}{\begin{eqnarray}}
\newcommand{\eea}{\end{eqnarray}}
\newcommand{\be}{\begin{equation}}
\newcommand{\ee}{\end{equation}}

\newcommand{\supp}{\mbox{\rm supp}}
\newcommand{\spec}{{\rm spec}}

\newcommand{\half}{\frac{1}{2}}

\newcommand{\eps}{{\varepsilon}}
\newcommand{\R}{{\mathbb R}}
\newcommand{\C}{{\mathbb C}}

\newcommand{\Compl}{{\mathbb C}}

\newcommand{\calg}{{\mathfrak g}}
\newcommand{\calG}{{\mathcal G}}

\newcommand{\calS}{{\mathcal S}}
\newcommand{\calU}{{\mathcal U}}

\newcommand{\rootsp}{{\mathcal N}}

\newcommand{\sign}{\mbox{sign}}

\newcommand{\diag}{{\rm diag}}

\newcommand{\Laplace}{\triangle}

\newcommand{\trip}{|\!|\!|}
\newcommand{\la}{\langle}
\newcommand{\ra}{\rangle}
\newcommand{\si}{\sigma}

\newcommand{\spa}{{\rm span}}
\newcommand{\calN}{{\mathcal N}}

\def\pa{\partial}

\def\nn{\nonumber}
\def\Xl{{\cal{X}}}
\def\Yl{{\cal{Y}}}

\newtheorem{theorem}{Theorem}[section]
\newtheorem{lemma}[theorem]{Lemma}
\newtheorem{defi}[theorem]{Definition}
\newtheorem{cor}[theorem]{Corollary}
\newtheorem{prop}[theorem]{Proposition}
\newtheorem{proposition}[theorem]{Proposition}

\theoremstyle{remark}
\newtheorem{remark}[theorem]{Remark}

\def\wt{\widetilde}

\def\ga{\gamma}
\def\de{\delta}

\def\bm{\left( \begin{array}{cc}}
\def\endm{\end{array}\right)}
\def\ker{{\rm ker}}
\def\Ran{{\rm Ran}}
\def\Hil{{\mathcal H}}
\def\Dom{{\rm Dom}}

\def\vpsi{{\vec \psi}}

\def\ep{\epsilon}
\def\lam{\lambda}
\newcommand{\calM}{{\mathcal M}}
\newcommand{\Lapl}{\frac12{\triangle}}
\def\siginf{\sigma^{\infty}}

\numberwithin{equation}{section}

\begin{document}

\parskip=2pt
\baselineskip=14.5pt
\thispagestyle{empty}

\title {Asymptotic stability of $N$-soliton states of NLS}
\author{ I. Rodnianski, W. Schlag and A. Soffer
\thanks{The research of I.R. was partially conducted during the period
he served as a Clay Mathematics Institute Long-Term  
Prize Fellow.  He was also supported in part by the NSF grant DMS-0107791.
W.S.~was partially supported by the NSF grant DMS-0070538
and a Sloan Fellowship. A.S.~was partially supported 
by the NSF grant DMS-0100490.}}
\maketitle

\begin{abstract}
We prove the asymptotic stability and asymptotic completeness of  an 
 arbitrary number of weakly interacting solitons for NLS.
 \end{abstract}

\section{Introduction}

The nonlinear Schr\"{o}dinger equation 
\be
i \, \frac{\partial \psi}{\partial t} 
\, = \, - \triangle \psi \, - \,
\beta (|\psi|^2) \psi, \hspace{.5in} 
x \in \mathbb{R}^n \tag{NLS} 
\ee
has in general (exponentially) localized
solutions in space, provided the nonlinearity has a negative (attractive)
part.  This is due to a remarkable cancellation of the dispersive effect
of the linear part with the focusing caused by the attractive nonlinearity.
To find such solutions, we look for time periodic solutions
$\psi \, \equiv \, e^{i \omega t} \, \phi_\omega ( x ) $.
It follows that $\phi_\omega$, if it exists, is a nonzero solution  of the problem
\be
\tag{ENLS}
- \omega \phi_\omega \, = \, - \triangle \phi_\omega \, - \, \beta ( |
\phi_\omega|^2 ) \phi_\omega.
\ee
We shall refer to such solutions as {\em nonlinear eigenfunctions}. 
In general, for $\phi_\omega$ to be localized (at least as $L^2$
function) we need $\omega > 0$.  

The general existence theory for this elliptic problem has been studied
in great detail, see the work of  Coffman~\cite{Cof}, Strauss~\cite{Str},
and Berestycki, Lions~\cite{BL}. 

It is easy to see that if $\psi_\omega (t, x ) = e^{i\omega t} 
\phi_{\omega} (x)$ is a solution of
NLS, then for any vector $\vec{a} \in \mathbb{R}^n$ the 
function $\phi_\omega (t, x - \vec a)$ is also a solution.  More generally, 
NLS is invariant under Galilean transformations
\begin{equation}
\label{eq:Galtr}
{\calg}_{\vec v, D}(t) := e^{-i\vec v\cdot x -\frac 12 |v|^{2}t} e^{i(t\vec v 
+D)p},
\end{equation}
and therefore we can
construct solutions  from $\phi_\omega$ which are moving with arbitrary velocity
$\vec{v}$.  As a result we obtain a family of 
 exponentially localized solutions
$$
\psi_{\vec v, \ga, D, \omega} = e^{i \vec v\cdot x - i\frac 12 (|v|^{2} 
-\omega) t + i \ga} \phi_{\omega} ( x -\vec v t - D)
$$
parametrized by a constant $(2n+2)$--dimensional vector 
$(\vec v, \ga, D, \omega )$, which are known as solitons.

Solving the initial value problem for NLS
requires understanding of two fundamental questions. The first
is the existence of global in time solutions. 
Due to the focusing character of the nonlinear term the global 
existence theory ought to be based on the $L^{2}$ conservation
law. By the results of Kato~\cite{Kato1} and Tsutsumi~\cite{Tsu}
 we can construct unique global solutions for 
any $L^{2}$ initial data under the assumption that the nonlinearity
$\beta$ satisfies the condition $|\beta(s)|\les (1+|s|)^{q}$ with
$q<\frac 2n$.

The second problem is that of the asymptotic behavior of solutions as
$t\to +\infty$. Guided by the completely integrable models in one dimension
we expect to have solutions with the asymptotic profile of $N$
independently moving solitons:
\begin{equation}
\label{eq:sum_soli}
\psi(t,x)\approx
\sum_{k=1}^{N}\psi_{\vec v_{k}, \ga_{k}, D_{k}, \omega_{k}} (t,x).
\end{equation}
Moreover, given such a solution $\psi(t,x)$ we expect that initial data
$\psi(0,x)+R$ for a suitably small perturbation~$R$ leads to a solution
$\tilde{\psi}$ with $\|\psi(t,\cdot)-\tilde\psi(t,\cdot)\|\to0$ as $t\to\infty$
in an appropriate norm.  The latter property is
known as asymptotic stability. In the context of NLS one needs to modify~\eqref{eq:sum_soli}
since one needs to make $(\vec v_{k}, \ga_{k}, D_{k}, \omega_{k})$ time-dependent. 
This fact was already observed in the context of orbital stability, see Weinstein~\cite{W1}.

In this paper we give an affirmative answer to the question of 
 existence and asymptotic stability of solutions with 
$N$-soliton profiles under the assumption of weak mutual interaction
between the solitons. A superposition of $N$ spatially separated moving 
solitons gives only an approximate solution of NLS. Our goal
is to show that the initial data 
\begin{equation}
\label{eq:indata}
\psi_{0}(x) = \sum_{k=1}^{N} e^{i\vec v_{k}(0)\cdot x+i\gamma_{k}(0)}
\phi_{\omega_{k}(0)}(x-D_{k}(0)) + R_{0},
\end{equation}
give rise to an asymptotically stable solution with the 
profile of $N$ independent solitons, with perhaps different
parameters $(\vec v_{k}, \ga_{k}, D_{k}, \omega_{k})\ne 
(\vec v_{k}(0), \ga_{k}(0), D_{k}(0),\omega_{k}(0))$.

The function $R_{0}$ is a perturbation satisfying a smallness assumption 
on its $L^{1}\cap L^{2}$ norm together with its derivatives. An important
aspect of our main result is the assumption that the solitons are weakly 
interacting. This condition can be enforced in two ways. Firstly, one can
assume that the initial shifts
$D_{k}(0)$ and the initial velocities $v_{k}$ are  chosen to 
model the case of non-colliding solitons:
\begin{equation}
\label{eq:separo}
|D_{k} + \vec v_{k} t - D_{k'} - \vec v_{k'} t|\ge L +ct,
\qquad k\ne k'
\end{equation}
for some sufficiently large constant $L$. 

\noindent Alternatively, one can assume that the relative
initial velocities of the solitons are large, i.e.,
\begin{equation}
\label{eq:vel_diff}
\min_{j\ne k} |\vec v_j -\vec v_k| > L.
\end{equation}
For the most part, we give details only for the case
of~\eqref{eq:separo} and leave the simple modifications
required by~\eqref{eq:vel_diff} to the reader. 
Note that~\eqref{eq:vel_diff} does not rule out that the
solitons collide. However, in view of~\eqref{eq:vel_diff}
the time of interaction is of size~$L^{-1}$, and therefore
the overall interaction remains weak. 

We shall also 
require certain spectral assumptions on the nonlinear 
eigenstates $\phi_{\omega}$ which will be explained below.

We believe that the methods we use may be applied to other classes of
equations with solitary type solutions and other symmetry groups (e.g.
Lorentz instead of Galilean), if and whenever certain linear $L^p$ decay
estimates can be verified for the linearized operators around one  such
soliton.  A detailed analysis of such $L^p$ estimates for NLS was
recently given in~\cite{RSS}.

To explain our results we recall the precise notions of
stability. Suppose we take the initial data of NLS to be an exact nonlinear eigenstate 
$\phi_\omega$, plus a small perturbation $R_0$.  What is then the expected 
behavior of the solution?  If the solution $\psi(t)$ stays near
the soliton 
$\psi_{\omega}(t)=e^{i\omega t}\phi_\omega (x)$ up to a phase and 
translation for all times (in $H^1$ norm) we say that the soliton $\psi_\omega$ is 
{\em orbitally stable}.
If, as time goes to infinity, the solution in fact converges in $L^2$ 
to a nearby soliton plus radiation\footnote{a function with asymptotic
behavior $e^{it\triangle} f$. Alternatively, we can and will
replace the $L^{2}$ by the $L^{\infty}$ convergence. In the latter 
topology the contribution of the radiation can be ignored.} 
we say that the solution is {\em asymptotically stable}.

Orbital stability of one soliton solutions has been subject of
extensive work in the last 20 years. The first results date back to
the work of Cazenave and Cazenave-Lions on logarithmic and monomial 
nonlinearities. The general case has been treated in the defining
works of Shatah-Strauss~\cite{ShSt}, Weinstein~\cite{W1} and~\cite{W2}, 
and Grillakis-Shatah-Strauss~\cite{GSS1}.
The general phenomena that has emerged from their results is 
that the orbital stability is essentially controlled by the sign
of the quantity $\pa_{\omega} \|\phi_{\omega}\|_{L^{2}}$
(stable, if positive, and unstable, if negative). However, all 
these results addressed orbital stability of a special class of 
solitons generated by ground states: positive, radial solutions of 
the equations ENLS. 

In \cite{BL} 
Berestycki-Lions proved 
the existence of a ground state in three or more
dimensions for any $\omega\ne 0$ under the conditions that 
the nonlinearity $\beta$ verifies
$
\lim\limits_{s \to +\infty} \, \beta (s^{2}) s^{-\frac 4{n-2}}\geq 0
$ 
and such that there exists
$0 < s_{0} < \infty$, with $G ( s_{0} ) > 0$, for 
$G ( s) \equiv 2{\int^{s}_{0}} (\beta(\tau^{2}) \tau  -\omega \tau ) 
d\tau$.  In fact, in their work ground states are found as 
minimizers of the constrained variational problem:
\be
\label{eq:convar}
\inf J[u]=\inf \left\{\int_{\R^{n}}|\nabla u|^{2}\,:\,\, \int_{\R^{n}} G(u)=1 \right\}.
\ee
The question of uniqueness of a ground state has been studied in 
\cite{McS}, \cite{Kw}, \cite{ML}. 

The asymptotic 
stability of {\em one} ground state 
soliton solutions of NLS and other
equations
was first shown for NLS with an extra attractive potential term in \cite{SW1}, 
\cite{SW2}, and  \cite{PW},
for one NLS soliton in dimension one in \cite{BP1} and in dimensions $n\ge 3$ 
in \cite{Cuc}; for NLS-Hartree see~\cite{FTY}. 


While the arguments for orbital stability were essentially based on 
Lyapunov type analysis and relied only on some limited information about 
the spectrum of the associated linear problem, the proofs of asymptotic 
stability required much more detailed properties of the related 
linearized systems. In particular, it led
to the need to impose additional spectral assumptions on the linear 
operators associated with a soliton. 

To describe the problem of linear stability (and spectral theory)
we linearize NLS around a soliton 
$w=e^{i\theta}\phi_{\omega}$, using the ansatz
\be
\label{eq:onesol}
\psi = e^{i \theta} \bigg (\phi_\omega(x-\vec vt) + R(t,x-\vec v 
t)\bigg)
\ee 
The resulting "linear" operator
acting on $R$, also has a term containing  
$\overline{R}$.   After complexifying the
space to 
$(R , \overline{R})$ we are left
with a matrix non self-adjoint operator of the type 
\begin{equation}
\label{eq:Hind}
{\cal{H}} \, = \, \left(
\begin{matrix}
L_{+} & W(x) \\
- {W}(x) & - L_{+}
\end{matrix}
\right)
\end{equation}
acting on $L^2 \times L^2$. Here 
\begin{align*}
L_{+}= - \Laplace +  \omega - \beta(\phi_{\omega}^{2}) - 
\beta'(\phi_{\omega}^{2})\phi_{\omega}^{2},\qquad
W(x)= -\beta'(\phi_{\omega}^{2}) \phi_{\omega}^{2}
\end{align*}
The operator $L_{+}$ is a self-adjoint perturbation of 
$- \triangle$ by an exponentially 
localized function, and $W$ is an exponentially localized potential.
Moreover, centering the perturbation $R$ around the soliton as in 
\eqref{eq:onesol} ensures that ${\cal H}$ is {\it{time-independent}}.   
Note that by non self-adjointness of $\cal H$ one can
no longer guarantee that $\sup_t\|e^{it{\cal H}}\|_{2\to2}<\infty$ where
$U(t)=e^{it{\cal H}}$ is given by $i\pa_{t} U + {\cal {H}}U=0$, $U(0)={\rm Id}$. In fact,
the operator ${\cal {H}}$ has a zero root space $\rootsp:=\bigcup_{\ell\ge 1}\,
\ker \cal H^{\ell}$
of dimension  at least 
$2n+2$ containing the eigenfunction $\phi_{\omega}$ as well as 
the elements generated from $\phi_{\omega}$ by infinitesimal 
symmetries of the problem. 
We decompose 
$$
L^{2}\times L^{2} = \rootsp + {\rootsp^{*}}^{\perp}
$$
and let $P$ denote the projector on the second term in this decomposition
(here $\rootsp^*=\bigcup_{\ell\ge 1}\ker (\cal H^*)^{\ell}$). It is easy
to see that $\|U(t)f\|_2$  grows polynomially  for some $f\in\rootsp$. On
the other hand, it is known from work of Weinstein~\cite{W1} that
\[\sup_{t}
\parallel U ( t ) P\psi \parallel_{L^2} \,< \, \infty
\tag{LS}
\]  
under certain conditions on the nonlinearity and provided $\phi_\omega$
is the (positive) ground state of ENLS. Generally speaking, we refer
to the property~(LS) as {\em linear stability}. 

In the case when ${\cal {H}}$ has one negative eigenvalue and 
$\phi_{\omega}$ is the unique ground state, one can show
that the condition $\sigma ( {\cal{H}} ) \subset \mathbb{R}$ is 
equivalent to the orbital stability condition 
$\pa_{\omega}\|\phi_{\omega}\|_{L^{2}}>0$ (see \cite{Grill}, \cite{BP1}).
Although due to the lack of self-adjointness this is not sufficient 
for linear stability, additional arguments show that (LS) in fact 
holds just under the above conditions (see \cite{W1}, \cite{GSS1}).

Linear stability of ground states
has been considered for a large class of NLS in the work of 
Weinstein~\cite{W1},
\cite{W2}, and Shatah, Strauss, Grillakis , \cite{ShSt}, \cite{GSS1}, \cite{GSS2}, 
\cite{Grill}, see also~\cite{SuSu}, \cite{Stu}. However, unconditional results were 
established only in the case of monomial nonlinearities 
$\beta(s)=s^{p}$. Moreover, linear stability has been 
shown to be essentially equivalent to the orbital stability
(except in the case of an $L^{2}$ critical nonlinearity 
$\beta(s)=s^{2/n}$).

Linear stability plays an essential role in the results on 
asymptotic stability. Moreover, the proofs of asymptotic stability
of one-soliton solutions required even more stringent assumptions
on the structure of the spectrum of ${\cal{H}}$, see Buslaev, Perelman~\cite{BP1}, 
and Cuccagna~\cite{Cuc}.  
This can be linked to the fact that on the linearized level 
asymptotic stability requires dispersive estimates of the type
\be
\label{eq:displin}
\|U(t) P \psi_{0}\|_{L^{\infty}}\les t^{-\frac n2}\|\psi_{0}\|_{L^{1}}
\ee
\cite{Cuc} or similar 
$L^{2}$-weighted decay estimates. 
To prove such estimates one needs to impose additional {\it {spectral conditions}} such as: 
absence of the discrete spectrum
for ${\cal H}$ on the subspace ${\rootsp}^{\perp}$,
absence of embedded eigenvalues, and absence of resonances at the edges of the continuous
spectrum. The dispersive estimates for such matrix Hamiltonians in 
dimensions $n\ge 3$ were proved by Cuccagna by an extension of the method introduced by Yajima
in the scalar case~\cite{Ya1}. In that approach, the decay estimates follows as 
a consequence of the proof of the $L^{p}\to L^{p},\, \forall p\in [1,\infty]$ 
boundedness of the wave operators.
In our recent work~\cite{RSS} we suggested a perhaps more 
straightforward approach for proving such estimates which instead relies
on construction of the analytic extension of the resolvent of 
${\cal H}$ and goes back to the work of Rauch~\cite{Rauch} in the scalar
case. This method, however, requires that all potential terms in the 
Hamiltonian ${\cal H}$ are exponentially localized functions, which 
perfectly fits the problem at hand. We refer to Hamiltonians 
verifying the required spectral assumptions (as well as the linear 
stability condition) as {\it admissible}.
The proof of the dispersive estimates crucially 
relies on the time independence of the Hamiltonian ${\cal H}$,
which was ensured by the choice of the ansatz.

In our work we choose the initial data of the form \eqref{eq:indata},
with $\phi_{\omega_{k}}$ verifying the elliptic problem ENLS,
satisfying the separation condition \eqref{eq:separo} and study the 
time asymptotic behavior of the corresponding solution of the 
time-dependent~NLS.

As in the study of asymptotic stability of one-soliton solutions 
the first objective is 
the analysis of the linearized problem.
It is natural to use the following ansatz for the solution:
\begin{equation}
\label{eq:Rdec}
\psi(t,x) = w+ R(t,x) = \sum_{k=1}^{N} w_{k} +  R(t,x),\qquad
w_{k}(t,x):=e^{i\theta_{k}} \phi_{\omega_{k}}
(x-\vec v_{k} t - D_{k})
\end{equation}
where $\theta_{k}$ are the phases associated with $k$-th soliton 
and $\vec v_{k}, D_{k}, \omega_{k}$ are its parameters\footnote{which
in the true ansatz are to be made time-dependent}. 
We substitute this ansatz into the equation (NLS) and retain only 
those terms which are linear in $R, \bar R$.
The resulting linear problem for the unknown $(R,\bar R)$ is 
\begin{equation}
\label{eq:SHt}
i\pa_{t} U + {\cal H}(t) U=0.
\end{equation}
It 
contains a {\it time-dependent} complex Hamiltonian
\[
{\cal{H}}(t) \, = \, \left(
\begin{matrix}
L_{+} & W(t,x) \\
- \overline{W}(t,x) & - L_{+}
\end{matrix}
\right)
\]
where 
\begin{align*}
L_{+}= - \Laplace -  \beta(|w|^{2}) -
\beta'(|w|^{2})|w|^{2},\qquad W(t,x) = -\beta(|w|^{2}) w^{2}. 
\end{align*}
Because of the smallness assumption on the initial perturbation
$R_{0}$ we can assume the the parameters of the final asymptotic 
profile will lie in a small neighborhood of the initial parameters
$(\vec v_{k}(0), \ga_{k}(0), D_{k}(0), \omega_{k}(0))$.
The separation condition on the initial parameters and the exponential 
localization of functions $\phi_{\omega_{k}}$ guarantees that $w$ 
is decomposed into a sum of functions of essentially disjoint support.
We thus can replace the Hamiltonian ${\cal H}(t)$ with 
\begin{equation}
\label{eq:chH}
{\cal{H}}(t) \, = \, \left(
\begin{matrix}
-\Laplace & 0 \\
0 & \Laplace \end{matrix}
\right) + \sum_{k=1}^{N}\left(
\begin{matrix}
U_{k}(x-\vec v_{k}t -D_{k}) & W_{k}(x-\vec v_{k}t-D_{k}) \\
- \overline{W_{k}}(x-\vec v_{k}t-D_{k}) & - U_{k}(x-\vec v_{k} t-D_{k})
\end{matrix}
\right)
\end{equation}
where 
\begin{align*}
&U_{k}=-\beta(\phi^{2}_{\omega_{k}}(x)) -
\beta'(\phi^{2}_{\omega_{k}}(x))
\phi^{2}_{\omega_{k}}(x))^{2},\\
& W_{k} = -e^{2i\theta_{k}} \beta(\phi^{2}_{\omega_{k}}(x)) 
\phi^{2}_{\omega_{k}}(x)
\end{align*}
This Hamiltonian belongs to the class of the so called 
{\it matrix charge transfer Hamiltonians}. Each of the Hamiltonians
$$
{\cal H}_{k}(t) =  \left(
\begin{matrix}
-\Laplace & 0 \\
0 & \Laplace \end{matrix}
\right) + \left(
\begin{matrix}
U_{k}(x-\vec v_{k}t-D_{k}) & W_{k}(x-\vec v_{k}t-D_{k}) \\
- \overline{W_{k}}(x-\vec v_{k}t-D_{k}) & - U_{k}(x-\vec v_{k}t-D_{k})
\end{matrix}
\right) = H_{0} + V_{k}(x-\vec v_{k} t)
$$
represents a linearization around the $k$th soliton. Moreover,
a specially chosen Galilei transform 
of the type \eqref{eq:Galtr} maps a solution of the linear 
problem $i\pa_{t} U + {\cal H}_{k}(t) U =0$ into a solution 
of the problem $i\pa_{t} U + {\cal H}_{k} U =0$ 
with a time-independent Hamiltonian ${\cal H}_{k}$ of the form 
\eqref{eq:Hind}.  Here $H_{0} = \diag (-\Laplace, \Laplace)$ and
$V_{k}$ is an exponentially localized complex matrix potential. 

Once again, at the linear level, the heart of the problem 
of asymptotic stability of N-solitons are the dispersive estimates
for the soutions of the equation \eqref{eq:SHt} with a matrix charge
transfer Hamiltonian ${\cal H}(t)$. Observe that in contrast to the 
one soliton case, the linearized Hamiltonian is time-dependent.

Due to the separation condition the problem \eqref{eq:SHt} admits
"traveling" bound states generated by the discrete spectrum of each
of the Hamiltonians ${\cal H}_{k}$. These bound states are formed by 
eigenfunctions and elements of the root space of 
${\cal H}_{k}$ boosted by the Galilei transform corresponding to 
the parameters $\vec v_{k}, D_{k}$. 
The paper \cite{RSS} establishes dispersive estimates 
\begin{equation}
\label{eq:dispH}
\|U(t)\psi_{0}\|_{L^{2}+L^{\infty}}\les (1+t)^{-\frac 
n2}\|\psi_{0}\|_{L^{1}\cap L^{2}}
\end{equation}
for the solutions of the linear time-dependent Schr\"odinger 
equation~\eqref{eq:SHt} with a matrix charge transfer Hamiltonian 
${\cal H}(t)$ of the type~\eqref{eq:chH} in dimensions $d\ge3$.
These estimates hold under the assumption that each of the 
time-independent Hamiltonians ${\cal H}_{k}$ is admissible (requiring 
spectral assumptions and the linear stability condition) and 
that the solution $U(t)\psi_{0}$ is {\cal asymptotically orthogonal}
to all traveling bound states of ${\cal H}_{k}(t)$. The latter 
means that for all $k=1,\ldots,N$
$$
\|P_{b}({\cal H}_{k},t) U(t)\psi_{0}\|_{L^{2}}\to 0,\quad 
{\text as}\quad t\to +\infty
$$
with $P_{b}({\cal H}_{k},t)$ denoting the time-dependent projection on the
$k$-th subspace of traveling bound states, i.e., it is the conjugation of
the spectral projection of the stationary operator ${\cal H}_k$  onto its bound states
by suitable Galilei transforms. The estimate \eqref{eq:dispH} is our main linear estimate.
However, to ensure that the perturbation $R(t,x)$ in the decomposition
\eqref{eq:Rdec} is asymptotically 
orthogonal to the subspace of traveling bound states and thus decays 
in the linear approximation with the rate of $t^{-\frac n2}$, we need to
allow the soliton parameters $\sigma_{k} = (\vec v_{k}, D_{k}, \gamma_{k}, 
\omega_{k})$ to become time-dependent. This in turn makes it necessary to change
the form of $w_k(t,x)$, for example, $x-\vec v_k t - D_k$ becomes
\[ x-\int_0^t \vec v_k(s)\, ds - D_k(t).\]
The resulting {\it nonlinear} problem for $R$, after the 
complexification $Z=(R,\bar R)$, takes the form
\begin{equation}
\label{eq:Zquat}
i\pa_{t} Z + H(t,\sigma(t)) Z = \dot \sigma \cdot \pa_{\si} w + {\cal N}(Z,w),
\end{equation}
where $H(t,\sigma(t))$ is a time-dependent Hamiltonian. It is
of the matrix charge transfer type provided $\sigma(t)=$const. The term 
$\pa_{\si} w$ denotes the derivative of the solitary approximation 
$w(t,x;\sigma(t))=\sum_{k=1}^N w_{k}(t,x;\sigma(t))$ with respect to its parameters $\sigma_{k}$, the term 
$\dot \sigma$ denotes the time derivative of the soliton parameters 
$\sigma_{k}$, and 
${\cal N}(Z,w)$ is a nonlinear term in $Z$.
We introduce the notion of an admissible path $\sigma(t)$ in the 
space of parameters and a 
{\it reference} Hamiltonian $H(t,\sigma)$ at infinity 
corresponding to the matrix
charge transfer Hamiltonian with fixed constants 
$\sigma_{k}=\sigma_{k}(t=\infty)=(\vec v_{k}, D_{k}, \gamma_{k}, 
\omega_{k})$. 
At each time $t$ the solution $Z(t)$ is required to be orthogonal
to the traveling bound states of the charge transfer Hamiltonians 
$H(t,\sigma=\sigma(t))$ obtained by fixing the parameters 
$\sigma_{k}=\sigma_{k}(t)$
at a given time $t$. This leads to the so called {\em modulation equations }
for $\sigma_{k}$, which couple the PDE~\eqref{eq:Zquat} for $Z$ with an ODE 
for the modulation parameters $\sigma_{k}(t)$.
To impose the orthogonality condition we first need to verify that 
it is satisfied initially. Using standard arguments, see e.g. \cite{BP1},
one can ensure this property by modifying the soliton parameters $\sigma(0)$ slightly in the 
decomposition of the initial data $\psi_0$,
$$
\psi_0(x) = \sum_{j=1}^N w_j(0,x;\sigma(0))+ R(0,x).
$$ 
We later justify the orthogonality condition at any positive time by showing 
that it is propagated.
To handle the nonlinear equation \eqref{eq:Zquat} we introduce the   
Banach spaces $\Xl_{s}$ and $\Yl_{s}$ of functions of $(t,x)$ 
\bea
\|f\|_{\Xl_{s}} &=&  \sup_{t\ge 0}\Big (\|\psi(t,\cdot)\|_{H^{s}} +  
(1+t)^{\frac n2}\sum_{k=0}^{s}
\|\nabla^{k}f(t,\cdot)\|_{L^{2}+L^\infty}\Big )\nn \\
\|F\|_{\Yl_{s}} &=& \sup_{t\ge 0}\sum_{k=0}^{s}
\Big (\int_{0}^{t} \|\nabla^{k}F(\tau,\cdot)\|_{L^{1}}\,d\tau +
(1+t)^{\frac n2+1}\|\nabla^{k} F(t,\cdot)\|_{L^{2}}\Big ). \nn\\ 
\eea
The space $\Xl_{s}$ is designed to control the solution $Z(t)$ itself, while 
the space $\Yl_{s}$ takes care of the nonlinear terms appearing as 
inhomogeneous terms in the equation for $Z$. We rewrite the equation 
\eqref{eq:Zquat} replacing the Hamiltonian $H(t,\sigma(t))$ with the
reference charge transfer Hamiltonian $H(t,\sigma)$. The solution 
operator of the corresponding linear problem 
\begin{equation}
\label{eq:inhF}
i\pa_{t} U + H(t,\sigma) U = F
\end{equation}
maps $\Yl_{s}\to \Xl_{s}$ for any integer $s$ uniformly in $\sigma$, 
provided that the solution is asymptotically orthogonal to all traveling 
bound state of $H(t,\sigma)$.
We then show that the nonlinearity $F$ arising from the inhomogeneous 
terms in the equation $Z$, written relative to the reference Hamiltonian
$H(t,\sigma)$, maps $\Xl_{s}\to \Yl_{s}$ for any $s>[\frac n2]+1$.
Given the smallness assumption on the initial data for $Z=(R,\bar R)$
this allows one to conclude the desired properties of $R$.
The modulation equations for $\sigma(t)$ are then used in turn 
to control the path $\sigma(t)$. In particular, we show that 
there exist final values of the parameters $\sigma=\sigma(\infty)$,
thus justifying the introduction of the reference Hamiltonian at 
infinity.
The estimates for the inhomogeneous problem \eqref{eq:inhF} phrased 
in terms of the mapping between $\Yl_{s}$ and $\Xl_{s}$ are 
essentially the linear dispersive estimates for the time-dependent 
Schr\"odinger equation with a matrix charge transfer Hamiltonian ${\cal 
H}(t)= H(t,\sigma)$ proved in \cite{RSS}.
As it was mentioned before, these estimates only require the 
admissibility of the corresponding individual 
time-independent matrix Hamiltonians ${\cal H}_{k}$ (representing the 
linearization around each individual soliton). 
The admissibility conditions can be somewhat loosely divided
into two categories:

1) Conditions related to  linear stability: 
$\sup_{t} \|e^{it{\cal H}_{k}}f\|_{2\to2} <\infty $

2) Spectral assumptions on ${\cal{H}}_{k}$ (absence of the embedded 
eigenvalues, resonances, etc.). 
\begin{remark}
We also require
the absence of "spurious" eigenvalues. This means that we assume that 
all of 
the discrete spectrum of ${\cal H}_{k}$ is generated purely by the 
nonlinear eigenfunction $\phi_{\omega_{k}}$, i.e., it is described by the 
generalized 0 eigenspace space of ${\cal H}_{k}$ which has the precise 
dimension $(2n+2)$ (coinciding with the dimension of the parameter 
space of $\sigma_{k}$). This is motivated by the requirement that 
the solution $Z$ has to be orthogonal to all traveling bound states 
of ${\cal H}_{k}$ for the dispersive estimates to hold, which can only 
be achieved by the choice of the parameters $\sigma_{k}$. 
In general, it is believed that the states corresponding to the 
spurious eigenvalues decay in time but the mechanism of this decay is 
purely nonlinear. We do not pursue 
this issue here.
\end{remark}
The key ingredient in establishing the linear stability is known 
to be the monotonicity or convexity condition 
\begin{equation}
\label{eq:stabil}
\pa_{\omega} \|\phi_{\omega}\|_{L^{2}} > 0 \qquad 
{\text {or}}\qquad \la L_{+}^{-1} \phi_\omega , \phi_\omega \ra > 0 .
\end{equation}
The known examples of when the monotonicity condition can be 
verified are limited to the case of the ground states $\phi$ 
corresponding to the monomial subcritical nonlinearities 
\be
\label{eq:subcriti}
\beta (s^{2}) = s^{p-1},\quad 1< p< 1 +\frac 4n,\qquad \forall\omega\ne 0
\ee
and the nonlinearity of the mixed type (see \cite{Sh})
$$
\beta (s^{2}) = s^{2} - s^{4}
$$
for the values of $\omega$ close to $3/16$. 
In this paper we find a new class of nonlinearities satisfying 
condition \eqref{eq:stabil}. These nonlinearities lie "near" 
the subcritical monomials of \eqref{eq:subcriti} but vanish much faster 
near $s=0$. More precisely we consider functions
\be
\label{eq:beta}
\beta_{\theta}(s^{2}) = s^{p-1} \frac {f(s^2)}{\theta + f(s^2)}
\ee
with a constant $\theta >0$ and the function $f$ satisfying the conditions
\be\label{eq:condit-f}
C_1 s^{r+1}\le |f(s^2)|\le C_2 s^{r+1}, \qquad 
s^2 |f'(s^2)|\le C_2 |f(s^2)|, \qquad p\in (1, 1+\frac 4n),\,\,
r>-1
\end{equation} 
and prove that given a sufficiently small
neighborhood $U$ in the space of parameters $\omega$ there exists a 
sufficiently small value $\theta_{0}$ such that for all $\theta 
<\theta_{0}$ and all $\omega \in U$ the ground state of
$\beta_{\theta}$ corresponding to $\alpha$ satisfies the monotonicity
condition.

We note that the higher rate of vanishing of $\beta(s^{2})$ at $s=0$
is important for asymptotic stability. In particular, it should be 
mentioned that if the power $p$ in the monomial example 
is too low ($p< 1 + \frac 2n$) even the scattering 
theory (asymptotic stability of a trivial 0 solution) fails. 

We now describe the structure of the paper.

{\bf Section 2.}\,  contains the statement of the main result together 
with the definitions of some of the fundamental objects used in the 
proof. The latter include the definition of an ansatz 
$\psi= w_{\sigma} + R$, 
separation condition on the initial 
data ensuring that the solitons are only weakly interacting, the 
notion of an admissible parameter path $\sigma(t)$, and the 
 spectral assumptions. It also contains 
the first discussion of the linearized Hamiltonians appearing in 
the later sections. We note here that traditionally the results on 
the asymptotic stability require the smallness assumption on the 
initial data in weighted Sobolev space. Our result  uses instead the 
Sobolev space based on the intersection of $L^{1}\cap L^{2} (\R^{n})$,
which has a distinct advantage of being translation invariant.

{\bf Section 3.}\, gives a detailed description of the linearization 
of the equation \eqref{eq:SHt} around an N-soliton profile 
$w_{\sigma}$ and 
introduces the notion of the reference charge transfer Hamiltonian at infinity.

{\bf Section 4.}\, describes the structure of the nullspaces of the 
Hamiltonians $H_{j}(\sigma)$, associated with the linearization on 
each individual nonlinear eigenfunction $\phi_{j}$. 

{\bf Section 5.}\, recalls the dispersive estimates for solutions 
of the time-dependent Schr\"odinger equation with a charge transfer 
Hamiltonian.

{\bf Section 6.}\, derives a system of ODE's for the modulation 
parameters $\sigma (t) = (\sigma_{1},\ldots,\sigma_{N})$ with 
$\sigma_{k}= (\vec v_{k}, D_{k}, \gamma_{k}, \alpha_{k})$ by requiring 
the complexified perturbation $Z=(R,\bar R)$ to be orthogonal to the
the unstable manifold comprised of the elements of the nullspaces of $H_{j}(\sigma(t))$.

{\bf Section 7.}\, gives the bootstrap assumption on the size of the 
perturbation $Z$ and the admissible path $\sigma(t)$ and provides 
estimates on the difference between the linearized Hamiltonian 
$H(\sigma(t))$ and the reference Hamiltonian $H(\sigma, t)$ at infinity.

{\bf Section 8.}\, provides the solution of the modulation equations 
for the path $\sigma(t)$.

{\bf Section 9.}\,solves the nonlinear equation for the complexified
perturbation $Z$ in the space $\Xl_{s}$ for $s>[\frac n2]+1$.  
This includes algebra estimates designated to show 
that the  $\Yl_{s}$-norm of the nonlinear terms in the equation
for $Z$ can be controlled by the $\Xl_{s}$-norm of $Z$ itself.

{\bf Section 10.} proves the scattering result showing that as $t\to \infty$
the solution $\psi(t)$ decomposes into the sum of solitons and a solution
of the free linear Schr\"odinger equation.

{\bf Section 11.} discusses the existence for the coupled
PDE-ODE system for $Z$ and $\sigma(t)$.

{\bf Section 12.} returns to the detailed discussion of the associated 
linear problems and proves some of the assertions made in the first 
part of the paper. We give the precise definition of an admissible 
Hamiltonian and start the investigation of their spectral properties.
In particular, in Section~11.2 we describe the spectrum of an admissible 
Hamiltonian and prove exponential decay of the 
elements of its generalized eigenspaces. 
In Section~11.3 we specialize to the admissible Hamiltonians arising 
from linearization around a nonlinear eigenfunction $\phi$.
We introduce and discuss the associated self-adjoint operators 
$L_{+}$ and $L_{-}$ and show that the admissibility conditions on our 
Hamiltonian (excluding the assumptions on absence of the embedded 
spectrum and resonances) can be reduced to the monotonicity condition 
\eqref{eq:stabil}, and the statement that the null space of $L_{-}$ 
is spanned by $\phi$ while  $L_{+}$ has a unique negative eigenvalue 
with the corresponding null space spanned by $\pa_{x_{j}}\phi$.
Some of the arguments in this section follow those of \cite{W1}, \cite{BP1}.

{\bf Section 13.} shows that for a very large class of nonlinearities the 
operator $L_{+}$, obtained by linearizing on a ground state, 
has a unique negative eigenvalue.

{\bf Section 14.}\, establishes the desired properties of the 
operators $L_{+}$ and $L_{-}$ for a particular class of nonlinearities
(namely those used in our main result) in the case when the nonlinear 
eigenfunction $\phi$ is a ground state.

For the most part of the paper (until Section~12) we do not 
specify the nature of nonlinear eigenstates 
$\phi_{\omega_{k}}$
of the elliptic problem (ENLS). In particular, we do not require them to
be ground states. Instead, we choose to formulate a more general 
conditional result dependent upon verification of certain precise 
properties of the linearized operators ${\cal H}_{k}$ associated 
with each $\phi_{\omega_{k}}$. It is only in Section~12 that we 
verify some of these assumptions in the case when $\phi_{\omega_{k}}$ are ground states.
The reason for choosing this approach is to emphasize the method
which allows us to handle weak interactions of (non-colliding) solitons 
generated by $\phi_{\omega_{k}}$,
provided that certain properties of each {\it{individual}} eigenstate
$\phi_{\omega_{k}}$ hold. We believe that our method will have an even wider 
range of applications than described here.

\section{Statement of results}
\label{sec:results}

Consider the NLS
\be
\label{eq:NLS}
i\partial_t \psi+\half\Laplace \psi + \beta(|\psi|^2)\psi =0
\ee
in $\R^n$, $n\ge 3$, with initial data
\be
\label{eq:init}
\psi_0(x) = \sum_{j=1}^N w_j(0,x) + R_0(x).
\ee
Here $w_j(0,x)$ are nonlinear eigenfunctions 
 generating the solitons
\bea
\label{eq:sol}
w_j(t,x)=w(t,x;\sigma_j(0)) &=& e^{i\theta_j(t,x)}
\phi(x-x_j(t),\alpha_j(0)) \\
\theta_j(t,x) &=&  v_j(0)\cdot x - \frac12(|v_j(0)|^2-\alpha_j^2(0))t 
+ \gamma_j(0) \label{eq:theta} \\
x_j(t) &=& v_j(0)t + D_j(0) \label{eq:xt}. 
\eea
and $\phi=\phi(\cdot,\alpha)$ is a 
solution of 
\be
\label{eq:ell} 
\half\Laplace \phi -\frac{\alpha^2}{2}\phi + \beta(|\phi|^2)\phi =0.
\ee
The solitons $w_j$ as in \eqref{eq:sol} satisfy \eqref{eq:NLS} with 
arbitrary constant parameters 
$\sigma_j(0)=(v_j(0),D_j(0),\gamma_j(0),\alpha_j(0))$.
We assume that the nonlinearity $\beta$ satisfies for all integers $\ell\ge0$
\begin{align}
|\beta^{(\ell)}(s)| &\les s^{\big(\frac{p-1}{2}-\ell\big)_+} \text{\ \ for\ }0\le s\le1 \label{eq:beta_1} \\
|\beta^{(\ell)}(s)| &\les s^{\big(\frac{q-1}{2}-\ell\big)_+} \text{\ \ for\ } s\ge 1 \label{eq:beta_2}
\end{align}
where $p\ge 2+\frac{2}{n}$, $q<1+\frac{4}{n}$. 
The main results of our paper are the following theorems.

\begin{theorem}
\label{thm:main1} Let $\beta$ be as in~\eqref{eq:beta}, \eqref{eq:condit-f}
with a sufficiently small value of parameter $\theta >0$. 
In addition, assume that $\beta(s^2)$ is a smooth function of $s$ and that 
$p+r\ge 1+\frac 2n$. Let 
\[ 
\psi_0(x) = \sum_{j=1}^N w_j(0,x) + R_0(x),
\]
see~\eqref{eq:init} be initial data for NLS that satisfy the 
 separation condition~\eqref{eq:separat}. Moreover, assume that
$w_j$ are defined in terms of ground states (positive radial variational solutions) of ENLS.
Finally, suppose that for all $\sigma\in\R^{N(2n+2)}$ with
$|\sigma-\sigma(0)|<c$, the linearized operators $H_j(\sigma)$ from~\eqref{eq:statHam}  have zero as their only eigenvalue, 
no resonances, and no embedded eigenvalues in their continuous spectrum. 
 Then there exists a positive $\epsilon$ such that for $R_{0}$
satisfying the smallness assumption 
\be
\label{eq:smalldata}
\sum_{k=0}^{s}\|\nabla^{k} R_{0}\|_{L^{1}\cap L^{2}} <\epsilon
\ee
for some integer 
$s>\frac n2$, there exists an admissible path $\sigma(t)$ with the limiting 
value $\sigma^\infty$
such that 
\[ \Big\|\psi(t,x)-\sum_{j=1}^N w_j(t,x;\sigma_j(t))\Big\|_{L_x^\infty }
\les (1+t)^{-\frac{n}{2}} \]
as $t\to\infty$.   
Moreover, there exists $u_0\in L^2$ so that
\[ \Big\|\psi(t,\cdot)-\sum_{j=1}^N w_j(t,x;\sigma^\infty) - e^{i\frac{t}{2}\Laplace} u_0 \Big\|_{L^2} 
\to0\]
as $t\to\infty$. 
\end{theorem}
\begin{remark}
We assume that $\alpha \ne 0$.
The smallness of the parameter $\theta$ in the nonlinearity $\beta$ depends
on the values of $\sigma(0)$.
\end{remark}
Our methods also allow treating initial conditions that are 
defined in terms of nonlinear eigenfunctions which are not
ground states as well as more general nonlinearities~$\beta$. 
However, in contrast to the case of a ground state and the
nonlinearities~\eqref{eq:beta}, in this case we cannot verify the
convexity conditions as well as various conditions related to the
spectrum of the linearized operators $H_j(\sigma)$. Therefore, we
need to include them into the hypotheses.

\begin{theorem}
\label{thm:main2} Let $\beta$ be as in~\eqref{eq:beta_1} and~\eqref{eq:beta_2}.
Impose the separation and convexity condition
see~\eqref{eq:separat}
and~\eqref{eq:stabil_2}, the spectral assumption from Definition~\ref{def:spass}. 
Suppose $\psi$ is the solution of~\eqref{eq:NLS} 
with initial condition~\eqref{eq:init} where the $w_j$ are
generated by nonlinear eigenfunctions of~ENLS (without assuming that they
are ground states). Then there 
exists a positive $\epsilon$ such that for $R_{0}$
satisfying the smallness assumption \eqref{eq:smalldata}
the conclusion of the previous theorem holds.
\end{theorem}
In the remainder of this section we shall discuss the theorem and its assumptions in more detail.

{\bf The soliton profiles $w_j(t,x;\sigma(t))$:}\ \ \ 
In addition to the one-soliton solutions $w_j (t,x;\sigma_j)$ with constant $\sigma_j$ 
introduced above, we need functions 
\begin{equation}
\label{eq:wj_guess}
w_j(t,x;\sigma_j(t))= e^{i\theta_j(t,x;\sigma(t))} \phi(x-x_j(t;\sigma(t)),\alpha_j(t)).
\end{equation}
The phase $\theta_j(t,x;\sigma(t))$ and the path $x_j(t;\sigma(t))$ are defined in terms of the time-dependent
parameters $\sigma_j(t)=(v_j(t), D_j(t), \gamma_j(t), \alpha_j(t))$ as follows:
\bea
\theta_j(t,x;\sigma(t)) &=&  v_j(t)\cdot x - \int_0^t 
\frac12(|v_j|^2-\alpha_j^2)(s)\,ds + \gamma_j(t) 
\label{eq:thetanew} \\
x_j(t;\sigma(t)) &=& \int_0^t v_j(s)\,ds + D_j(t). \label{eq:xtnew}
\eea
Henceforth, the functions $w_j(t,x;\sigma(t))$ or 
simply $w_j(\sigma(t))$, correspond to the soliton moving 
along the time-dependent
curve $\sigma(t)$ in the parameter space according 
to~ \eqref{eq:thetanew} 
and~\eqref{eq:xtnew}, while $w_j(\sigma)$ is the true soliton 
moving along the straight line determined by an arbitrary
constant $\sigma$ as in \eqref{eq:theta} and \eqref{eq:xt}.

{\bf Admissible paths $\sigma(t)$:}\ \ \
We collect the individual parameter curves $\sigma_j(t)$ from above 
into a single curve $\sigma(t):=(\sigma_1(t),\ldots,\sigma_N(t))\subset \R^{(2n+2)N}$. 
Given the initial value $\sigma(0)$ we introduce the set of 
{\em admissible curves} $\sigma(t)$
as those $C^1$ curves that remain in a small neighborhood 
of $\sigma(0)$ for all times and converge
to their final value $\sigma(\infty)=\lim_{t\to +\infty}\sigma(t)$.  
We shall also impose 
the condition that for an admissible curve~$\sigma(t)$
\be
\int_0^\infty\int_s^\infty |\dot v_j(\tau)\cdot v_j(\tau)- 
\dot \alpha (\tau)_j\alpha_j(\tau)|\,d\tau\,ds <\infty,
\qquad \int_0^\infty \int_s^\infty |\dot v_j(\tau)|\,d\tau\,ds <\infty
\label{eq:integr}
\ee
for all $1\le j\le N$. Given an admissible curve $\sigma(t)$ we define 
the constant vector $\sigma^\infty$ in the following fashion:
\bea
&&v^\infty_j=v_j(\infty),\qquad
D^\infty_j=D_j(\infty) - \int_0^\infty \int_s^\infty 
\dot v_j(\tau)\,d\tau\,ds,\label{eq:vecnew}\\
&&\gamma^\infty_j = \gamma_j(\infty) +  \int_0^\infty 
\int_s^\infty (\dot v_j(\tau)\cdot v_j(\tau) - 
\dot \alpha_j(\tau) \alpha_j (\tau))\,d\tau\,ds,\qquad
\alpha^\infty_j=\alpha_j(\infty). \label{eq:scalnew}
\eea
{\bf Convexity condition:}\ \ \ 
We impose the {\em convexity condition} 
\be
\label{eq:stabil_2}
\la \partial_\alpha \phi(\cdot;\alpha), \phi(\cdot;\alpha)\ra > 
0,\qquad 
\forall \alpha : \, \min_{j=1,..,N} |\alpha - \alpha_j(0)| < c
\ee 
for some positive constant $c$. As noted in the Introduction, the
convexity condition is closely connected with the issue of orbital
stability of the individual solitons $w_j(t,x;\sigma_j(0))$.

{\bf Separation conditions:}\ \ \ 
Our theorem handles the case of so-called weakly interacting solitons. 
This means that the initial positions $D_j(0)$
and initial velocities $v_j(0)$ are such that for all~$t\ge 0$
one has the (physical) {\em separation condition}
\be
\label{eq:separat}
|D_j(0) + v_j(0) t - D_{\ell}(0) - v_{\ell}(0) t|\ge L + ct,\qquad 
\forall j\ne\ell=1,\ldots, N
\ee
with some sufficiently large constant $L$ and a positive constant $c$. 
Another assumption under which our theorems hold equally well is the 
condition of {\em large relative velocities} of the solitons. This means that
\begin{equation}
\label{eq:vel_sep}
 \min_{j\ne \ell}|v_j-v_\ell|> L
\end{equation}
for some large $L$. 

\noindent Let $\alpha_{\min}=\min_{1\le j\le k}\alpha_j(0)-c$. 
It will be understood henceforth that 
\be
\alpha_{\min} L \ge |\log \epsilon |.
\label{eq:quantsep}
\ee 
The small constant $\epsilon$ appears in our theorem as a measure of smallness
of the initial perturbation $R_0$. 

{\bf Spectral assumptions:}\ \ \
We will write $w_j(\sigma(t))=w_j(t,x;\sigma(t))= 
e^{i\theta_j(\sigma(t))}\,\phi_j(\sigma(t))$
where 
$\phi_j(\sigma(t))=\phi(t,x;\sigma_j(t))=\phi(x-x_j(t);\alpha_j(t))$. 
Linearizing the equation~\eqref{eq:NLS} around the 
state $w=\sum_{j=1}^N w_j,\,\,\psi=w+R$ one obtains the 
following system of equations for $Z=\binom{R}{\bar R}$:
\bea
\label{eq:linsystem}
i\partial_t Z + H(t,\sigma(t)) 	Z &=& F.
\eea
Here $H(t,\sigma(t))$ is the time-dependent matrix Hamiltonian
\bea
\label{eq:matrix}
&& \qquad\qquad H(t,\sigma(t))= H_0 +  \\ 
&&\sum_{j=1}^N \left (\begin{array}{cc}
\beta(|w_j(\sigma(t))|^2) 
+\beta'(|w_j(\sigma(t))|^2)|w_j(\sigma(t))|^2 &  
\beta'(|w_j(\sigma(t))|^2)w_j^2(\sigma(t)) \nn \\
-\beta'(|w_j(\sigma(t))|^2)\bar{w}^2_j(\sigma(t)) & 
-\beta(|w_j(\sigma(t))|^2) -
\beta'(|w_j(\sigma(t))|^2)|w_j(\sigma(t))|^2 
		    \end{array}\right ) \nn \\
&&\qquad \qquad \qquad = H_0 + \sum_{j=1}^N V_j(t,x;\sigma(t)),\label{eq:matrixV}\\
&&\qquad \qquad H_0=\left (\begin{matrix} \half\Laplace & 0 \\
                                              0 & -\half\Laplace 
                           \end{matrix} \right ) \nn 
\eea 
with complex matrix time-dependent potentials $V_j(t,x;\sigma(t))$ dependent on $w_j$ 
and $\sigma(t)$. 
The right-hand side $F$ in \eqref{eq:linsystem} depends on $\dot \sigma$, $w$, and 
nonlinearly on $Z$. For a given constant parameter vector~$\sigma$ 
we shall introduce the Hamiltonian
$H(t,\sigma)$
\bea
&& H(t,\sigma)=H_0 +  \label{eq:ref_hamil}\\ 
&&\sum_{j=1}^N \left (\begin{array}{cc}
\beta(|w_j(\sigma)|^2) +\beta'(|w_j(\sigma)|^2)|w_j(\sigma)|^2 &  
\beta'(|w_j(\sigma)|^2)w^2_j(\sigma)  \\
-\beta'(|w_j(\sigma)|^2)\bar{w}^2_j(\sigma) & -\beta(|w_j(\sigma)|^2) 
-
\beta'(|w_j(\sigma)|^2)|w_j(\sigma)|^2 
		    \end{array}\right ). \nn 
\eea
We refer to Hamiltonians of the  form~\eqref{eq:ref_hamil} 
as {\em matrix charge transfer} Hamiltonians. They are
discussed in more detail in Section~\ref{ap:lin}, as well
as in~\cite{RSS}. 
Recall that $w_j(\sigma)$ denotes the soliton moving along 
the straight line determined by the constant parameters~$\sigma_j$. 
The proof of our theorem relies on dispersive
estimates for matrix charge transfer Hamiltonians that
were obtained in~\cite{RSS}, see also Section~\ref{ap:lin}
below. For these estimates to hold, one needs to impose
certain {\em spectral conditions} on the 
stationary Hamiltonians 
\be
H_j(\sigma) := \left (\begin{array}{cc}
\half\Laplace -\frac{\alpha^2}2 + \beta(\phi_j(\sigma)^2) 
+\beta'(\phi_j(\sigma)^2)\phi_j(\sigma)^2 &  
\beta'(\phi_j(\sigma)^2)\phi^2_j(\sigma)  \\
-\beta'(\phi_j(\sigma)^2)\phi^2_j(\sigma) &-\half\Laplace 
+ \frac{\alpha^2}2 -\beta(\phi_j(\sigma)^2) -
\beta'(\phi_j(\sigma)^2)\phi_j(\sigma)^2 
	\end{array}\right ) \label{eq:statHam} 
\ee 
where $\phi_j(\sigma) = \phi(x,\alpha_j)$, see~\eqref{eq:ell}. 
These Hamiltonians arise from the matrix charge
transfer problem by applying a Galilei transform to 
the $j^{th}$ matrix potential in~\eqref{eq:ref_hamil} 
so that this potential becomes stationary (strictly
speaking, this also requires a modulation which leads to the 
spectral shift $\frac{\alpha^2}{2}$ in~\eqref{eq:statHam}). 
We impose the spectral assumption as described by
the following definition.  

\begin{defi}
\label{def:spass}
We say that the 
{\em spectral assumption} holds, provided for all $\sigma\in \R^{N(2n+2)}$
with $|\sigma-\sigma(0)|<c$ one has
\begin{itemize}
\item $0$ is the only point of the discrete spectrum of 
$H_j(\sigma)$ and the dimension of the corresponding root space
is $2n+2$, 
\item each of the $H_j(\sigma)$ is admissible in the sense
of Definition~\ref{def:spec_ass} below and the stability 
condition
\[ \sup_t \|e^{itH_j(\sigma)}P_s \|_{2\to 2} <\infty, \]
see~\eqref{eq:stab_cond}, holds (here $P_s$ is the projection
onto the scattering states associated with~$H_j$, see~\eqref{eq:split}).
\end{itemize}
\end{defi}

\noindent While the second condition is known to hold generically in 
an
appropriate sense, see Section~\ref{ap:lin},  the first
condition is more restrictive and not believed to hold generically.

\section{Reduction to the matrix charge transfer model}
 
For the sake of simplicity we consider the case of two solitons, 
i.e., $N=2$. Setting 
\begin{equation}
\label{eq:w_sum2}
 w_1(\sigma(t))+w_2(\sigma(t))=w,\qquad \psi=w+R,
\end{equation}
where $w_j$ are as in~\eqref{eq:wj_guess}, \eqref{eq:thetanew}, and~\eqref{eq:xtnew},
we derive from \eqref{eq:NLS} that
\bea
\label{eq:splitt}
&& i\partial_t R + \half\Laplace R + (\beta(|w|^2)
+\beta'(|w|^2)|w|^2)R + \beta'(|w|^2)w^2\,\bar{R} \\
&& = -(i\partial_t w + \half\Laplace w + \beta(|w|^2)w) 
+ O(|w|^{p-2}|R|^2) + O(|R|^p), \nn
\eea
using the assumptions \eqref{eq:beta_1}, \eqref{eq:beta_2}
on the nonlinearity~$\beta$. Observe that
\bea
i\partial_t w + \half\Laplace w + \beta(|w|^2)w &=& 
- \sum_{j=1}^2 \Big[\big( \dot{v}_j(t)\cdot x + 
\dot{\gamma}_j(t)\big)w_j(\sigma(t)) + i e^{i\theta_j(\sigma(t))}
\nabla \phi_j(\sigma(t))\cdot \dot{D}_j(t) \nn \\
&&  - ie^{i\theta_j(\sigma(t))}\partial_{\alpha}\phi_j(\sigma(t))
\dot{\alpha}_j(t) \Big] + O(w_1w_2). 
\label{eq:sigmadot}
\eea
In view of~\eqref{eq:w_sum2} one has
\bea
&& i\partial_t R + \half\Laplace R + (\beta(|w|^2)
+\beta'(|w|^2)|w|^2)R + \beta'(|w|^2)w^2\,\bar{R} \nn\\
&=&
i\partial_t R + \half\Laplace R + \sum_{j=1}^2 
\Bigl[\beta(|w_j(\sigma(t))|^2)+\beta'(|w_j(\sigma(t))|^2)
|w_j(\sigma(t))|^2\Bigr]R + \sum_{j=1}^2 
\beta'(|w_j(\sigma(t))|^2)w_j(\sigma(t))^2\,\bar{R} \nn \\
&&   + O(w_1(\sigma(t))w_2(\sigma(t)))R. \nn
\eea
Rewriting the equation \eqref{eq:splitt} as a system 
for $Z=(R,\bar{R})$ therefore leads to
\bea
\label{eq:system}
i\partial_t Z + H(\sigma(t))Z &=& \dot\Sigma W(\sigma(t)) 
+ O(w_1w_2)Z + O(w_1w_2) + O(|w|^{p-2}|Z|^2) + O(|Z|^p).
\eea
Here $H(\sigma(t))$ is the time-dependent matrix Hamiltonian 
from~\eqref{eq:matrix} and 
\be
\dot \Sigma W(\sigma(t))= \left(\begin{matrix} f \\
                                              -\bar f
                           \end{matrix} \right)
\label{eq:Sigmadot}
\ee
where
\be
f=\sum_{j=1}^2 \Big[( \dot{v}_j(t)\cdot x + 
\dot{\gamma}_j(t))w_j(\sigma(t)) + 
i e^{i\theta_j(\sigma(t))}\nabla \phi_j(\sigma(t))\cdot \dot{D}_j(t) 
  - ie^{i\theta_j(\sigma(t))}\partial_{\alpha}
\phi_j(\sigma(t))\dot{\alpha}_j(t)\Big] \label{eq:fdef}
\ee
We shall assume that $\sigma(t)$ is an admissible path with 
initial values $\sigma(0)$ in the sense of \eqref{eq:integr}.
Given an admissible curve $\sigma(t)$ we introduce the 
{\em reference Hamiltonian} $H(t,\siginf)$ ``at infinity''
\bea
&& H(t,\siginf)=H_0 +  \nn \\ 
&&\sum_{j=1}^2 \left (\begin{array}{cc}
\beta(|w_j(\siginf)|^2) +\beta'(|w_j(\siginf)|^2)|w_j(\siginf)|^2 &  
\beta'(|w_j(\siginf)|^2)w^2_j(\siginf)  \\
-\beta'(|w_j(\siginf)|^2)\bar{w}^2_j(\siginf) & -\beta(|w_j(\siginf)|^2) 
-
\beta'(|w_j(\siginf)|^2)|w_j(\siginf)|^2 
		    \end{array}\right ) \label{eq:Htsigma}
\eea
where $\siginf=(\siginf_1,\ldots,\siginf_N),\,\,
\siginf_j=(v^\infty_j, D^\infty_j, \gamma^\infty_j, \alpha^\infty_j)$ 
is the constant vector determined by the 
curve $\sigma(t)$ as in \eqref{eq:vecnew} and \eqref{eq:scalnew}:
\bea
&&v^\infty_j=v_j(\infty),\qquad
D^\infty_j=D_j(\infty) - \int_0^\infty \int_s^\infty 
\dot v_j(\tau)\,d\tau\,ds,\nn\\
&&\gamma^\infty_j = \gamma_j(\infty) +  \int_0^\infty 
\int_s^\infty (\dot v_j(\tau)\cdot v_j(\tau) - 
\dot \alpha_j(\tau) \alpha_j (\tau))\,d\tau\,ds,\qquad
\alpha^\infty_j=\alpha_j(\infty). \nn
\eea
Recall that $w_j(\siginf)$ is the soliton moving 
along the straight line determined by the constant
parameters $\siginf_j$.
For $j=1,\ldots,N$ we introduce the Hamiltonians 
\be
H_j(t,\siginf)=H_0 + \left (\begin{array}{cc}
\beta(|w_j(\siginf)|^2) +\beta'(|w_j(\siginf)|^2)|w_j(\siginf)|^2 &  
\beta'(|w_j(\siginf)|^2)w^2_j(\siginf)  \\
-\beta'(|w_j(\siginf)|^2)\bar{w}^2_j(\siginf) & -\beta(|w_j(\siginf)|^2) 
-
\beta'(|w_j(\siginf)|^2)|w_j(\siginf)|^2 
		    \end{array}\right ) \label{eq:timedep} 
\ee
together with their stationary counterparts 
$H_j(\siginf)$ as in \eqref{eq:statHam} with $\sigma=\siginf$.

The following lemma relates the evolutions corresponding to 
the Hamiltonians $H_j(t,\sigma)$ and 
$H_j(\sigma)$ for an arbitrary $\sigma$ by means of a modified Galilean transformation.

\begin{lemma}
\label{lem:conjug}
Let $U_j(t,\sigma)$ be the solution operator of the equation
\bea
&&i \partial_t U_j(t,\sigma) + H_j(t,\sigma)U_j(t,\sigma)=0,
\label{eq:evolU}\\
&& U_j(0,\sigma)=I\nn
\eea
and $e^{it H_j(\sigma)}$ be the corresponding 
propagator for the time-independent matrix Hamiltonian
$H_j(\sigma)$. Then
\be
\label{eq:conjug}
U_j(t,\sigma) = \calG^*_{v_j, D_j}(t) \calM_j^*(t,\sigma) 
e^{i t H_j(\sigma)} \calM_j(0,\sigma) \calG_{v_j, D_j}(0)
\ee
where $\calG_{v_j, D_j}(t)$ is the diagonal matrix 
Galilean transformation
\be
\label{eq:matrGal} 
\calG_{v_j, D_j}(t) \binom{f_1}{f_2} =
\binom{\calg_{v_j,D_j}(t)f_1}{\overline{\calg_{v_j,D_j}(t)\bar f_2}} 
\ee
and 
\be
\label{eq:matrM}
\calM_j(t,\sigma) = \left (
\begin{array}{cc}
e^{-i \frac {\alpha_j^2}2 t - i(v_j\cdot D_j + \gamma_j)} & 0\\
0 & e^{i \frac {\alpha_j^2}2 t + i(v_j\cdot D_j + \gamma_j)}
\end{array}\right ).
\ee
\end{lemma}
\begin{proof}
By definition,
\bea
i\dot U_j &=& i\dot \calG^*_{v_j, D_j}(t) \calM_j^*(t) 
e^{itH_j(\sigma)} \calM_j(0)\calG_{v_j, D_j}(0) + 
\calG^*_{v_j, D_j}(t) i \dot \calM_j^*(t) e^{itH_j(\sigma)} 
\calM_j(0)\calG_{v_j, D_j}(0) \nn \\
 && - \calG^*_{v_j, D_j}(t) \calM_j^*(t) H_j(\sigma) 
e^{itH_j(\sigma)} 
\calM_j(0)\calG_{v_j, D_j}(0). \label{eq:2ndline}
\eea
Clearly,
\[ 
i\dot \calM_j^*(t) = \left( \begin{matrix} 
 -\frac{\alpha_j^2}{2}\, e^{i \frac {\alpha_j^2}2 t 
+ i(v_j\cdot D_j + \gamma_j)} & 0 \\
0 &  \frac{\alpha_j^2}{2}\, e^{-i \frac {\alpha_j^2}2 t 
- i(v_j\cdot D_j + \gamma_j)}
                           \end{matrix}
                    \right),
\]
whereas one checks that
\[ 
i\dot \calG^*_{v_j, D_j}(t) \left( \begin{matrix} f_1 \\ {f}_2 
                                   \end{matrix}
                            \right)
= i\left( \begin{matrix} \dot\calg^*_{v_j,D_j} f_1 \\ 
\overline{\dot\calg^*_{v_j,D_j} \bar{f}_2} 
          \end{matrix}
   \right)
= \left( \begin{matrix} -(v_j^2/2 - v_j\cdot p)\;\calg^*_{v_j,D_j} 
f_1 \\ 
\;\;(v_j^2/2 + v_j\cdot p)\;\overline{\calg^*_{v_j,D_j} \bar{f}_2 } 
          \end{matrix}
   \right).
\]
Finally, we need to move $H_j(\sigma)$ to the left 
in~\eqref{eq:2ndline}. 
We consider the differential operator separately from the matrix potential, i.e., 
\bea
H_j(t,\sigma) &=& H_0+\left(\begin{matrix} U_j(x-x_j(t)) & e^{2i\theta_j(t,x)} W_j(x-x_j(t)) \\
 				 -e^{-2i\theta_j(t,x)} W_j(x-x_j(t)) & - U_j(x-x_j(t))
		       \end{matrix} \right)  \nn \\
H_j(\sigma) &=& H_0 + \left(\begin{matrix} -\frac{\alpha_j^2}{2} & 0 \\
  				      0 & \frac{\alpha_j^2}{2} 
                      \end{matrix}\right) + \left(\begin{matrix} U_j &  W_j \\
 			         	 - W_j & - U_j
		       \end{matrix}\right)   \label{eq:matr_rep}
\eea
where $x_j(t),\theta_j(t)$ are as in~\eqref{eq:xt}, \eqref{eq:theta},
and $U_j=\beta(\phi_j(\sigma)^2) +\beta'(\phi_j(\sigma)^2)\phi_j(\sigma)^2$, 
$W_j= \beta'(\phi_j(\sigma)^2)\phi^2_j(\sigma)$. Note that, on the one hand, $\calM_j$ commutes 
with all matrices in~\eqref{eq:matr_rep} that do not involve $U_j,W_j$. 
On the other hand, one has
\begin{align}
 & H_0 \calG^*_{v_j, D_j}(t) 
\left( \begin{matrix} f_1 \\ f_2 \end{matrix} \right)
 - \calG^*_{v_j, D_j}(t) \left[H_0 + \left(\begin{matrix} -\frac{\alpha_j^2}{2} & 0 \\
  				      0 & \frac{\alpha_j^2}{2} 
                      \end{matrix}\right) \right] 
 \left( \begin{matrix} f_1 \\ f_2 \end{matrix} \right)  = 
 \left( \begin{matrix} \frac12\alpha_j^2 \;\calg^*_{v_j, D_j}(t)f_1 \\
     -\frac12\alpha_j^2 \;\overline{\calg^*_{v_j, D_j}(t)\bar f_2} 
  \end{matrix} \right)      
\nn \\ 
 & + \half\left( \begin{matrix}
 e^{i\frac{v_j^2}{2}t}\,\Laplace \Bigl( e^{i{x\cdot v_j}}\,
e^{-i{v_j}\cdot(D_j+tv_j)} \,f_1(x-tv_j-D_j)\Bigr)  
- e^{i\frac{v_j^2}{2}t}\,e^{i{x\cdot v_j}}\,e^{-i{v_j}
\cdot(D_j+tv_j)} \,\Laplace f_1(x-tv_j-D_j)   \\
 -e^{-i\frac{v_j^2}{2}t}\,\Laplace \Bigl( e^{-i{x\cdot v_j}}
\,e^{i{v_j}\cdot(D_j+tv_j)} \, f_2(x-tv_j-D_j)\Bigr)  
+ e^{-i\frac{v_j^2}{2}t}\,e^{-i{x\cdot v_j}}\,e^{i{v_j}
\cdot(D_j+tv_j)} \,\Laplace  f_2(x-tv_j-D_j) 
          \end{matrix} 
    \right) \nn \\
 &= \left( \begin{matrix} 
               \frac12 (v_j^2 + \alpha_j^2)\;
\calg^*_{v_j, D_j}(t)f_1 - v_j\cdot p\, \calg^*_{v_j, D_j}(t)f_1 \\ 
               -\frac12 (v_j^2 + \alpha_j^2)\;
\overline{\calg^*_{v_j, D_j}(t)\bar f_2} - v_j\cdot p\, 
\overline{\calg^*_{v_j, D_j}(t)\bar{f_2}} 
            \end{matrix}
     \right). \nn
\end{align}
Finally, we need to deal with the matrix potentials. Write 
$\calM(t):=\calM_j(t,\sigma) = \left(\begin{matrix} e^{-i\omega(t)/2} & 0 \\ 0 & e^{i\omega(t)/2} 
\end{matrix}\right)$ and set $\rho=t|\vec v_j|^2 + 2x\cdot \vec v_j$. Then 
(omitting the index $j$ for simplicity)
\bea
&& 
\calM(t)\calG_{\vec v,D}(t) \bm U(\cdot-\vec vt-D) & e^{2i\theta}W(\cdot-\vec vt-D) \\
 -e^{-2i\theta}W(\cdot-\vec vt-D) & -U(\cdot-\vec 
vt-D) \endm \binom{f_1}{f_2} \nn \\
 &=&  
\bm e^{-i\omega(t)/2} & 0 \\ 0 & 
e^{i\omega(t)/2} \endm \binom{\calg_{\vec v,D}(t) U(\cdot-\vec v t-D)f_1 +
\calg_{\vec v,D}(t) e^{2i\theta} W(\cdot -\vec vt-D)f_2}
{\overline{-\calg_{\vec v,D}(t) 
e^{2i\theta} W(\cdot-\vec v t-D)\overline{f_1}} - 
\overline{\calg_{\vec v,D}(t)  U(\cdot -\vec vt-D)\overline{f_2}}} \nn \\
 &=&   
\binom{U\calg_{\vec v,D}(t) (e^{-i\omega(t)/2} f_1) 
+ W e^{-i(v^2 t+2 x\cdot \vec v)} e^{i(2\theta(t,\cdot+t\vec v+D)-\omega)} 
\overline{\calg_{\vec v,D}(t) \overline{e^{i\omega(t)/2} f_2}}}
{-We^{i(v^2t+2x\cdot\vec v)} e^{i(\omega-2\theta(t,\cdot+t\vec v+D))} \;\calg_{\vec v,D}(t) (e^{-i\omega(t)/2} f_1) 
- U \;\overline{\calg_{\vec v,D}(t) \overline{e^{i\omega(t)/2} f_2}}} \nn \\
&=& \bm U &  e^{i(2\theta(t,\cdot+t\vec v+D)-\omega-\rho)} W \\ -e^{-i(2\theta(t,\cdot+t\vec v+D)
-\omega-\rho)} W & -U \endm
\bm e^{-i\omega(t)/2} & 0 \\ 0 & e^{i\omega(t)/2} \endm \binom{\calg_{\vec v}(t) f_1}
{\overline{\calg_{\vec v,D}(t)\overline{f_2}}}\nn.
\eea
Now $2\theta(t,\cdot+t\vec v+D)-\rho-\omega= 2\vec v\cdot x+(|\vec v\,|^2+\alpha^2)t+2\gamma
+2\vec v\cdot D-t|\vec v\,|^2 - 2x\cdot \vec v - \omega =0$
by definition of $\omega$, i.e., $\omega=\alpha^2 t+2\gamma+2\vec v\cdot D$. 
Adding these expressions shows that
\[ i\dot U_j(t,\sigma) + H_j(t,\sigma)U_j(t,\sigma)=0,\]
as claimed.
\end{proof}

\section{The root spaces of $H_j(\sigma)$ and $H_j^*(\sigma)$}

In view of Section~\ref{ap:lin} below 
(see in particular Definition~\ref{def:spec_ass} as well 
as~\eqref{eq:split}) we will need to understand 
the generalized eigenspaces of the stationary
operators~$H_j(\sigma)$ from~\eqref{eq:stat}. 
By our spectral
assumption, see Definition~\ref{def:spass} above,
only generalized eigenspaces at~$0$ are allowed.
We denote these spaces by~$\rootsp_j(\sigma)$ and refer 
to them as root spaces. 
Thus, $\rootsp_j(\sigma)=\ker\big(H_j(\sigma)^2\big)$ and 
by~\eqref{eq:split} one has the direct 
(but not orthogonal) decomposition
\[ L^2(\R^3)\times L^2(\R^3) = {\rootsp_j^*(\sigma)}^\perp + \rootsp_j(\sigma),\]
where $\rootsp_j^*(\sigma)=\ker\big(H_j^*(\sigma)^2\big)$.
The (nonorthogonal) projection onto ${\rootsp_j^*(\sigma)}^\perp$ 
associated with  this decomposition is denoted 
by~$P_j(\sigma)$. While the evolution $e^{itH_j(\sigma)}$ 
is {\em unbounded} on~$L^2$ as $t\to\infty$, it is known 
in many cases that it remains bounded on~$\Ran(P_j(\sigma))$.
In Section~\ref{ap:lin} this is referred to as the 
{\em linear stability assumption}. 

\begin{prop}
\label{prop:null} Impose the hypotheses of Theorem~\ref{thm:main1} and
let $H_j(\sigma)$ be as in \eqref{eq:statHam}. Then 
\begin{itemize}
\item The  nullspace $\rootsp_j^*(\sigma)$
of $H^*_j(\sigma)$ is given by the following vector valued 
$2n+2$ functions $\xi^m_j(x;\sigma),\,\, m=1,\ldots,2n+2$:
\begin{align*}
&\xi^m_j(x;\sigma)=\binom{u^m_j(x;\sigma)}{\bar u^m_j(x;\sigma)},&\qquad & \\
& u^1_j(x;\sigma) = \phi_j(x;\sigma), &\qquad & H^*_j(\sigma) \xi^1_j(\cdot;\sigma) = 0,\\
& u^2_j (x;\sigma) = i \frac 2{\alpha_j} \partial_\alpha 
\phi_j(x;\sigma),&\qquad & 
H^*_j(\sigma)\xi^2_j(\cdot;\sigma) = -i \xi^1_j(\cdot;\sigma),\\
& u^m_j(x;\sigma) = i \partial_{x_{m-2}} \phi_j (x;\sigma),&\qquad & 
H^*_j(\sigma) \xi^m_j (\cdot;\sigma)= 0, \quad m=3,..,n+2,\\
& u^m_j (x;\sigma)=  x_{m-n-2} \phi_j(x;\sigma), &\qquad & 
H^*_j(\sigma) \xi^m_j(\cdot;\sigma) = -2i \xi^{m-n}_j(\cdot;\sigma),\quad m=n+3,..,2n+2
\end{align*}
\item Let
\[ J= \left(\begin{matrix} 0 & 1 \\
                    -1 & 0 
     \end{matrix} \right).
\]
Then $J$ is an isomorphism between the nullspaces of $H_j^*(\sigma)$ 
and $H_j(\sigma)$. In particular,
the nullspace of $H_j(\sigma)$ has a basis $\{J\xi_j^m(\cdot;\sigma)\:|\: 1\le m\le 
2n+2\}$. Moreover, $\rootsp_j^*(\sigma)$ is spanned by
\[ J\, \partial_{\sigma_r} W_j(t,x;\sigma) \text{\ \ for\ \ } 1\le r\le 2n+2,\]
where $W_j(t,x;\sigma) = \binom{w_j(t,x;\sigma)}{\bar{w}_j(t,x;\sigma)}$. 
\item One has the linear stability property 
\[ \sup_{t}\Big\|e^{itH_j(\sigma)} P_j(\sigma) 
\Big\|_{2\to2} < \infty\]
where $P_j(\sigma)$ is the projection onto 
${\rootsp_j^*(\sigma)}^\perp$ as introduced above.
\end{itemize}  
\end{prop}
\begin{proof} 
\end{proof}

For the case of monomial, subcritical nonlinearities these
results go back to Weinstein's work on 
modulational stability~\cite{W1}.

\section{Estimates for the linearized problem}

\noindent In~\eqref{eq:system}  we obtained the system 
\be
\label{eq:main}
i\partial_t Z + H(t,\siginf) Z = \big (H(\sigma(t))- H(t,\siginf)\big) Z+  
\dot\Sigma W(\sigma(t)) + O(w_1w_2)Z + O(w_1w_2) +
O(|w|^{p-2}|Z|^2) + O(|Z|^p),
\ee
The point of rewriting~\eqref{eq:system}
in this form is to be able to use the dispersive
estimates that were obtained in~\cite{RSS} for
(perturbed) matrix charge transfer Hamiltonians, 
see also Sections~\ref{ap:lin} and Section~8 in~\cite{RSS}. 

\begin{theorem}
\label{thm:charge0}
Let $Z(t,x)$ solve the equation
\bea
&& i \partial_t Z + H(t,\sigma) Z  = F,
\label{eq:chargeZ} \\
&& Z(0,\cdot) = Z_0(\cdot)\nn
\eea 
where the matrix charge transfer Hamiltonian $H(t,\sigma)$ 
satisfies the conditions of Definition~\ref{def:chargetransm}. 
Assume that $Z$ satisfies
\be 
\label{eq:asorthZ}
\|({\rm Id}-P_{j}(\sigma))\calM_j(\sigma,t) \calG_{v_j, D_j}(t) 
Z(t,\cdot)\|_{L^2} \le
B (1+t)^{-\frac n2},
\quad \forall j=1,\ldots,k, 
\ee
with some positive constant $B$, where 
$\calM_j(\sigma,t)$ and $ \calG_{v_j, D_j}(t)$ are as
in Lemma~\ref{lem:conjug}. 
Then $Z$ verifies the following decay estimate 
\be 
\label{eq:decayZ}
\|Z(t)\|_{L^{2}+ L^\infty} \les (1+t)^{-\frac n2} 
\Big(\|Z_0\|_{L^1\cap L^2} + \trip F\trip +B \Big)
\ee
for $t>0$ with 
\be
\nn
\trip F\trip := \sup_{t\ge 0} \Bigl[ \int_0^{t}\|F(s)\|_{L^{1}}\,ds 
+ (1+t)^{\frac n2+1}\|F(t)\|_{L^{2}}\Bigr].
\ee
In addition, we also have the $L^{2}$ estimate
\be
\label{eq:mass}
\|Z(t)\|_{L^{2}}\les \|Z_{0}\|_{L^{1}\cap L^{2}} + \trip F\trip + B
\ee
\end{theorem}

For the proof see \cite{RSS} and Section~\ref{ap:lin} below. 
In particular, note that~\eqref{eq:asorthZ} is related to the 
characterization of scattering states in Definition~\ref{def:asympm}.

In the applications the inhomogeneous term $F$ is a nonlinear
expression which depends on~$Z$. Therefore, in addition to the
estimates~\eqref{eq:decayZ} and~\eqref{eq:mass} we shall need 
corresponding estimates for the derivatives of~$Z$.

For an integer $s\ge 0$ we define  Banach spaces $\Xl_{s}$ and 
$\Yl_{s}$ of functions of 
$(t,x)$ 
\bea
\|\psi\|_{\Xl_{s}} &=&  \sup_{t\ge 0}\Big (\|\psi(t,\cdot)\|_{H^{s}} +  
(1+t)^{\frac n2}\sum_{k=0}^{s}
\|\nabla^{k}\psi(t,\cdot)\|_{L^{2}+L^\infty}\Big )\label{eq:Xl}\\
\|F\|_{\Yl_{s}} &=& \sup_{t\ge 0}\sum_{k=0}^{s}
\Big (\int_{0}^{t} \|\nabla^{k}F(\tau,\cdot)\|_{L^{1}}\,d\tau +
(1+t)^{\frac n2+1}\|\nabla^{k} F(t,\cdot)\|_{L^{2}}\Big )  \label{eq:Yl}
\eea
The generalization of the estimates of Theorem~\ref{thm:charge0} 
is given by the 
following theorem (see Section~9, in particular 
Proposition~9.3 in~\cite{RSS} for the proof).
\begin{theorem}
Under assumptions of Theorem \ref{thm:charge0} we have that
for any integer $s\ge 0$
\be
\label{eq:decZ}
\|Z\|_{\Xl_{s}}\les \sum_{k=0}^{s}\|\nabla^{k}
Z(0,\cdot)\|_{L^{1}\cap L^{2}} + \|F\|_{\Yl_{s}} + B
\ee
\label{thm:charge}
\end{theorem}

We apply Theorem~\ref{thm:charge} to the equation~\eqref{eq:main}.
This will, in particular, lead to our main result, i.e., that 
$\|Z(t)\|_\infty \les t^{-\frac n2}$ as $t\to\infty$.
We need to ensure that $Z$ is a scattering solution
relative to each of the channels of the charge transfer Hamiltonian~$H(t,\sigma)$, 
in the sense of the estimate~\eqref{eq:asorthZ}. 
Analogous to Buslaev, Perelman~\cite{BP1} this will
be accomplished by an appropriate choice of the
path~$\sigma(t)$, to be made in the following section.

In order to prove existence of solutions $Z$ and $\sigma$
we require another version of Theorem~\ref{thm:charge}, 
which follows easily from that result. 

\begin{remark} It is easy to see that a time-localized version
of the previous theorem also holds. Indeed, let 
\bea
\|\psi\|_{\Xl_{s}(T)} &=&  \sup_{0\le t\le T}\Big (\|\psi(t,\cdot)\|_{H^{s}} +  
(1+t)^{\frac n2}\sum_{k=0}^{s}
\|\nabla^{k}\psi(t,\cdot)\|_{L^{2}+L^\infty}\Big )\label{eq:Xl_loc}\\
\|F\|_{\Yl_{s}(T)} &=& \sup_{0\le t\le T}\sum_{k=0}^{s}
\Big (\int_{0}^{t} \|\nabla^{k}F(\tau,\cdot)\|_{L^{1}}\,d\tau +
(1+t)^{\frac n2+1}\|\nabla^{k} F(t,\cdot)\|_{L^{2}}\Big ).  \label{eq:Yl_loc}
\eea
Then assuming~\eqref{eq:asorthZ} for $0\le t\le T$ with a constant $B_T$, one has
\be
\label{eq:decZ_loc}
\|Z\|_{\Xl_{s}(T)}\les \sum_{k=0}^{s}\|\nabla^{k}
Z(0,\cdot)\|_{L^{1}\cap L^{2}} + \|F\|_{\Yl_{s}(T)} + B_T
\ee
\end{remark}

\begin{cor}
\label{cor:charge2}
Let $\tilde \xi_j^m(t,x)$, $1\le m\le 2n+2, 1\le j\le N$,
be a collection of smooth functions such that
\begin{equation}
\label{eq:diff_xi} 
\sup_{t\ge 0} \|\calM_j(\sigma,t) \calG_{v_j, D_j}(t) \tilde \xi_j^m(t,x) - 
\xi_j^m(\cdot;\sigma) \|_{L^1\cap L^2} \le \delta
\end{equation}
for some small $\delta>0$ and some given $\sigma$. 
Let $Z$ be a solution of 
\begin{align}
\label{eq:ctransV}
i\partial_t Z + H(t,\sigma) Z &= V(t,x) Z + F, \\
\big\la Z(t), \tilde \xi_j^m(t,\cdot) \big\ra &= 0 \label{eq:weird_orth}
\end{align}
for all $t\ge 0$, where $V(t,x)$ is a smooth function that satisfies 
$\sup_{|\gamma|\le s}\|\partial_x^\gamma V(t,\cdot)\|_{L^1\cap L^\infty} <\delta(1+t)^{-1}$,
with a nonnegative integer $s$ for all $t>0$. Then
\begin{equation}
\label{eq:decay10}
\|Z\|_{\Xl_s} \les \sum_{k=0}^{s}\|\nabla^{k}
Z(0,\cdot)\|_{L^{1}\cap L^{2}} + \|F\|_{\Yl_{s}}.
\end{equation}
\end{cor}
\begin{proof}
By \eqref{eq:weird_orth},
\begin{align}
 \big|\big\la \calM_j(\sigma,t) \calG_{v_j, D_j}(t) Z(t), \xi_j^m(\cdot;\sigma)\big \ra\big| &\le \big|\big\la \calM_j(\sigma,t) \calG_{v_j, D_j}(t) Z(t), 
-\calM_j(\sigma,t) \calG_{v_j, D_j}(t) \tilde \xi_j^m(t,\cdot)+\xi_j^m(\cdot;\sigma) \big\ra\big| \nn \\
& \le \delta \|Z(t)\|_{2+\infty}. \nn
\end{align}
Therefore, for all $0\le t\le T$,
\begin{align}
 & \big|\big\la \calM_j(\sigma,t) \calG_{v_j, D_j}(t)Z(t), \xi_j^m(t,\cdot;\sigma)\big \ra\big|
\le \delta (1+t)^{-\frac{n}{2}} \sup_{0\le \tau\le T}(1+\tau)^{\frac{n}{2}} \|Z(\tau)\|_{2+\infty} \nn \\
& \le \delta (1+t)^{-\frac{n}{2}} \|Z\|_{\Xl_s(T)} =:  (1+t)^{-\frac{n}{2}} B_T. \nn
\end{align}
Hence, using \eqref{eq:decZ_loc} one sees that 
\begin{align}
 \|Z\|_{\Xl_s(T)} &\les \sum_{k=0}^{s}\|\nabla^{k}
Z(0,\cdot)\|_{L^{1}\cap L^{2}} + \|F\|_{\Yl_{s}(T)} +\|V Z\|_{\Yl_{s}(T)} + B_T  \nn \\
&\les \sum_{k=0}^{s}\|\nabla^{k}
Z(0,\cdot)\|_{L^{1}\cap L^{2}} + \|F\|_{\Yl_{s}(T)} + \delta \|Z\|_{\Xl_s(T)},\nn
\end{align}
and the desired conclusion follows.
\end{proof}

\section{Modulation equations}

In their analysis of the stability relative to one soliton,
Buslaev and Perelman~\cite{BP1}, \cite{BP2},
 and Cuccagna~\cite{Cuc} derive
the equations for $\dot\sigma$ by imposing an 
orthogonality condition
on the perturbation~$Z$ for all times. More precisely, 
they make the ansatz
\be
\label{eq:ansatz}
 \psi=e^{i\theta(t,\sigma(t))}(w(\sigma(t))+R)
\ee
where $e^{i\theta(t,\sigma(t))}w(\sigma(t))$ is a 
single soliton evolving along
a nonlinear set of parameters. The removal of the phase 
from the perturbation~$R$
leads to an equation which is simply the translation 
of the equation involving
the stationary Hamiltonian~\eqref{eq:statHam} to 
the point $vt+D$. This in turn
makes it very easy to formulate the 
orthogonality conditions: At time~$t$, 
the function  $R(\cdot+vt+D)$  in~\eqref{eq:ansatz} 
needs to be perpendicular to all elements of
the generalized eigenspaces of all~$H_j(\sigma)^*$ as 
in~\eqref{eq:statHam}, where $\sigma$ is equal to 
the parameters~$\sigma(t)$ at time~$t$. 

\noindent In the multi-soliton case the removal of the phases
by means of this ansatz is not available, since distinct
solitons carry distinct phases. As already indicated above, 
we work with the representation
\[ \psi(t) = \sum_{j=1}^N w_j(t,\sigma(t)) + R,\]
which forces us to formulate the orthogonality 
condition in terms of a set of 
functions that is moving along with the~$w_j(t,\sigma(t))$. 
We now define these functions.

\begin{defi}
\label{def:movbas}
 Let $\sigma(t)$ be an admissible path and define 
$\theta_j(t,x;\sigma(t))$ and~$x_j(t;\sigma(t))$
as in~\eqref{eq:thetanew} and~\eqref{eq:xtnew}. Also, set 
$\phi_j(t,x;\sigma(t))=\phi(x-x_j(t;\sigma(t));\alpha_j(t))$.
Then we let
\[ \xi_j^m(t,x;\sigma(t)) = 
\left( \begin{matrix} u_j^m(t,x;\sigma(t)) \\ \bar 
u_j^m(t,x;\sigma(t))
				\end{matrix}
			 \right)
\]
with 
\begin{align}
\label{eq:xitsigma}
   u_j^1(t,x;\sigma(t))  &= w_j(t,x;\sigma(t)) = 
e^{i\theta_j(t,x;\sigma(t))}\,
 \phi_j(t,x;\sigma(t)) \\
   u_j^2(t,x;\sigma(t))&= 
\frac{2i}{\alpha_j}\,e^{i\theta_j(t,x;\sigma(t))}\,\partial_{\alpha}
\phi_j(t,x;\sigma(t)) \nn \\
    u_j^m(t,x;\sigma(t))  &= ie^{i\theta_j(t,x;\sigma(t))}\,
\partial_{x_{m-2}}\phi_j(t,x;\sigma(t)) 
\text{\ \ for\ \ }3\le m\le n+2\nn \\
    u_j^m(t,x;\sigma(t)) &= 
e^{i\theta_j(t,x;\sigma(t))}\,(x^{m-n-2}-x_j^{m-n-2}(t;\sigma(t)))
\phi_j(t,x;\sigma(t)),\text{\ for\ }\,\,n+3\le m\le 2n+2. \nn
\end{align}
\end{defi} 

The following proposition should be thought of as a time-dependent 
version of 
Proposition~\ref{prop:null}. More precisely, if $\sigma$ is a {\em 
fixed} set
of parameters, then one can define an alternate set of vectors, 
$\tilde{\xi}_j^m$, say,
by applying appropriate Galilean transforms to the stationary 
vectors in Proposition~\ref{prop:null}. For example, take some $\xi_j^m$
so that $H_j^*(\sigma)\xi_j^m=0$. Then the corresponding~$\tilde{\xi}_j^m$
satisfies 
\[ i\partial_t \tilde{\xi}_j^m + H_j(t,\sigma) \tilde{\xi}_j^m = 0,\]
with $H_j(t,\sigma)$ as in \eqref{eq:timedep}.
Naturally, one would therefore expect that
\[ i\partial_t {\xi}_j^m + H(\sigma(t)) {\xi}_j^m = O(\dot \sigma_j) 
+ O(e^{-ct}), \]
where $H(\sigma(t))$ is as in~\eqref{eq:matrix} 
(the exponentially decaying term
appears because of interactions between solitons). 
The following proposition shows that this indeed holds, 
but as in~\cite{Cuc} we will work
with a modified set of parameters 
$\tilde \sigma_j(t)=(v_j(t),D_j(t),\alpha_j(t),\tilde\gamma_j(t))$ where
\be
\label{eq:tildega}
\dot{\tilde \gamma}_j(t)= \dot{\gamma}_j(t) + \half \sum_{m=1}^n 
\dot{v}^m_j(t)  x^{m}_j(t,\sigma(t)).
\ee
The point of this modification is that the $\dot \Sigma W(\sigma(t))$ 
term in~\eqref{eq:main} and~\eqref{eq:system} can be rewritten as
\bea
\dot \Sigma W(\sigma(t))&=&\sum_{j=1}^k \Big[ \dot{\tilde\gamma}_j(t) 
 J \xi^1_j(t,x;\sigma(t))- \frac{\alpha_j}2\dot{\alpha}_j(t) 
J\xi^2_j(t,x;\sigma(t))\Big ]+ \nn\\
&&\sum_{j=1}^k \sum_{m=1}^n \Big [\dot{D}^m_j(t) J 
\xi^{m+2}_j(t,x;\sigma(t)) +
\half \dot{v}^m_j(t) J \xi^{m+n+2}_j(t,x;\sigma(t)) \Big 
],\label{eq:dotS}
\eea
where $\xi_j^m$ are as in Definition~\ref{def:movbas}. 
This is of course due to the fact that
passing to $\tilde \gamma_j$ allows us to change from 
$x$ to $x-x_j(t;\sigma(t))$ in~\eqref{eq:fdef}.

\begin{prop}
\label{prop:movingxi} Let $\sigma(t)$ be an admissible path and define
$\xi_j^m(t,x;\sigma(t))$ as in Definition~\ref{def:movbas}. Then 
\bea
i\partial_t \xi_j^1 + H^*_j(\sigma(t)) \xi_j^1 &=& 
O\bigl(\dot{ \tilde \sigma} (|\phi_j|+|D \phi_j|)\bigr) \label{eq:uzhas1} \\
i\partial_t \xi_j^2 + H^*_j(\sigma(t)) \xi_j^2 &=& i\xi_j^1 
+ O\bigl(\dot{ \tilde \sigma }
(|\phi_j|+|D \phi_j|+|D^2 \phi_j|)\bigr) \label{eq:uzhas2} \\
i\partial_t \xi_j^m + H^*_j(\sigma(t)) \xi_j^m &=&  
O\bigl(\dot{\tilde \sigma }
(|\phi_j|+|D \phi_j|+|D^2 \phi_j|)\bigr)
\text{\ \ for\ \ }3\le m\le n+2 \label{eq:uzhas3} \\
i\partial_t \xi_j^m + H^*_j(\sigma(t)) \xi_j^m &=& -2i\xi_j^{m-n} 
+ O\bigl(\dot{\tilde \sigma} 
(|\phi_j|+|D \phi_j|+|D^2 \phi_j|)\bigr) \text{\ for\ }\,\,n+3\le 
m\le 2n+2. \label{eq:uzhas4} 
\eea
Here $D$ refers to either spatial derivatives 
$\partial_{x^\ell}$ or derivatives $\partial_{\alpha}$.
Moreover, as in Definition~\ref{def:movbas}, the function 
$\phi_j$ needs to be evaluated at
$x-x_j(t;\sigma(t))$, $\alpha_j(t)$. 
\end{prop}
\begin{proof} This is verified by direct differentiation 
of the functions in Definition~\ref{def:movbas}.
\end{proof}

\noindent The following proposition collects the modulation equations 
for the path
$\sigma(t)$ that are obtained by taking scalar products 
of~\eqref{eq:matrix} with the 
functions $\xi_j^{m}$ from Definition~\ref{def:movbas}.
This will of course use~\eqref{eq:dotS}. The modulation equations
are derived from the orthogonality assumptions, see~\eqref{eq:perp} below. 
Observe that these assumptions need not be satisfied at $t=0$. 
Nevertheless, as in Buslaev and Perelman~\cite{BP1}, one shows by means of
the implicit function theorem that one can replace the initial 
decomposition~\eqref{eq:init} by a nearby one which does satisfy the 
orthogonality condition. This uses  the smallness
of the initial perturbation $R_0$, as well as the separation conditions~\eqref{eq:separat}
or~\eqref{eq:vel_sep}. The details can be found in Section~\ref{sec:exist}, see Lemma~\ref{lem:init_orth}.
In that section it is also shown that, conversely, given the modulation of Proposition~\ref{prop:modeq}
the orthogonality condition will propagate if satisfied initially.

\begin{prop}
\label{prop:modeq} 
Let $Z$ satisfy the system~\eqref{eq:system}. Suppose that for all 
$t\ge0$, 
\be
\label{eq:perp} 
\la Z(t), \xi_j^m(t,\cdot;\sigma(t)) \ra  = 0 \text{\ \ for all\ \ } j,m
\ee
where $\xi_j^m$ is as in Definition~\ref{def:movbas}. Then
the path $\tilde\sigma(t):= (v_j(t), D_j(t), \tilde {\gamma}_j(t), 
\alpha_j(t)),\,j=1,..,n$ 
satisfies the following system of equations with matrix potentials $V_r(t,x;\sigma(t))$
as in \eqref{eq:matrixV}:
\begin{align}
 -2i\dot \alpha_j(t) \Bigl \la \phi_j(\sigma(t)),\, 
\partial_{\alpha} \phi_j(\sigma(t)) \Bigr \ra + &
 O(\dot{\tilde\sigma} \|Z(t)\|_{L^{2}+L^\infty}) =  
\sum_{r\ne j}\Bigl \la   V_r(t,\cdot;\sigma(t)) Z,\,\xi^1_j(t,\cdot;\sigma(t))
\Bigr \ra +\nn \\ 
& \Bigl \la  \Big( O(w_1w_2)Z + O(w_1w_2) + O(|w|^{p-2}|Z|^2) 
+ O(|Z|^p)\Big ),\, \xi^1_j(t,\cdot;\sigma(t))\Bigr \ra 
,\label{eq:sigmadoteq}\\
 2i \dot {\tilde\gamma}_j(t) \Bigl \la \phi_j(\sigma(t)),\, 
\partial_{\alpha} \phi_j(\sigma(t)) \Bigr \ra + 
 & O(\dot{\tilde \sigma} \|Z(t)\|_{L^{2}+L^\infty}) = 
\sum_{r\ne j}\Bigl \la   V_r(t,\cdot;\sigma(t)) 
Z,\,\xi^2_j(t,\cdot;\sigma(t))\Bigr \ra +\nn \\ 
& \Bigl \la  \Big( O(w_1w_2)Z + O(w_1w_2) + O(|w|^{p-2}|Z|^2) 
+ O(|Z|^p)\Big ),\, \xi^2_j(t,\cdot;\sigma(t))\Bigr \ra ,\nn\\
 \dot v_j^m(t) \|\phi_j(\sigma(t))\|_2^2 + & O(\dot{\tilde \sigma}
\|Z(t)\|_{L^{2}+L^\infty}) = \sum_{r\ne j}\Bigl \la  
 V_r(t,\cdot;\sigma(t)) Z,\,\xi^{m+2}_j(t,\cdot;\sigma(t))\Bigr \ra +\nn \\ 
& \Bigl \la  \Big( O(w_1w_2)Z + O(w_1w_2) + O(|w|^{p-2}|Z|^2) 
+ O(|Z|^p)\Big ),\, \xi^{m+2}_j(t,\cdot;\sigma(t))\Big \ra,\nn\\
 \dot D_j^m(t) \|\phi_j(\sigma(t))\|_2^2  + & O(\dot{\tilde \sigma}
\|Z(t)\|_{L^{2}+L^\infty}) = \sum_{r\ne j}\Bigl 
\la   V_r(t,\cdot;\sigma(t)) Z,\,\xi^{n+m+2}_j(t,\cdot;\sigma(t))\Bigr \ra +\nn 
\\ 
& \!\!\!\!\!\!\!\!\!\! \Bigl \la  \Big( O(w_1w_2)Z + O(w_1w_2) 
+ O(|w|^{p-2}|Z|^2) + O(|Z|^p)\Big ),\, \xi^{n+m+2}_j(t,\cdot;\sigma(t))
\Bigr \ra.\nn
\end{align}
\end{prop}
\begin{proof} 
Differentiating \eqref{eq:perp} yields
\[ \la i\partial_t Z, \xi_j^m(t,\cdot;\sigma(t)) \ra = \la  Z, i\partial_t\xi_j^m(t,\cdot;\sigma(t)) 
\ra.\]
Taking scalar products of \eqref{eq:system} thus leads to 
\[
\la Z,i\partial_t \xi_j^m \ra  + \la Z,H^*(\sigma(t))\xi_j^m \ra 
= \Bigl \la \dot\Sigma W(\sigma(t)),\xi_j^m\Bigr \ra + 
 \Bigl \la O(w_1w_2)Z + O(w_1w_2) + O(|w|^{p-2}|Z|^2) + O(|Z|^p), 
\xi_j^m \Bigr \ra.
\]
In view of the explicit expressions \eqref{eq:xitsigma} one has
\bea 
\la J\xi_j^2(t,\cdot;\sigma(t)),\,\xi_j^1(t,\cdot;\sigma(t)) \ra &=& -2i 
\la \phi_j(\sigma(t)),\, \partial_{\alpha} \phi_j(\sigma(t)) \ra \nn \\
\la J\xi_j^m(t,\cdot;\sigma(t)),\,\xi_j^1(t,\cdot;\sigma(t)) \ra &=& 0 
\text{\ \ for\ \ } m\ne 2 \nn \\
\la J\xi_j^m(t,\cdot;\sigma(t)),\,\xi_j^2(t,\cdot;\sigma(t)) \ra &=& 0 
\text{\ \ for\ \ } m\ne 1 \nn \\
\la J\xi_j^{m+2}(t,\cdot;\sigma(t)),\,\xi_j^{m+n+2}(t,\cdot;\sigma(t)) 
\ra &=& -2i\|\phi_j(\sigma(t))\|_2^2 \text{\ \ for\ \ } 3\le m\le n+2. 
\nn
\eea
Therefore, the proposition follows by taking inner products 
in~\eqref{eq:dotS}. 
Note that the terms containing 
$\dot{\tilde\sigma}\|Z(t)\|_{L^{2}+L^\infty}$ appear 
from 
Proposition \ref{prop:movingxi}.
\end{proof}

\section{Bootstrap assumptions}

The proof of our main theorem relies on the bootstrap assumptions on 
the 
admissible path $\sigma(t)$ 
and the size of the perturbation $Z(t,x)=\binom{R(t,x)}{\bar R(t,x)}$.
in the norms of the spaces $\Xl_{s}$ defined in \eqref{eq:Xl}.

\underline{Bootstrap assumptions} 

There exists  a small constant $\delta=\delta(\epsilon)$ dependent on 
the size of the initial 
data $R_0$ and the initial separation of the solitons 
$w_j(0,x;\sigma(0))$, see \eqref{eq:separat}, and 
a sufficiently large constant $C_0$ such that for some integer 
$s>\frac n2$
\begin{align}
& |\dot{\tilde\sigma}(t)|\le \delta^2 (1+t)^{-n}, \quad \forall t\ge 
0, \label{eq:bootsigma} \\ 
& \|Z\|_{\Xl_{s}} \le \delta C_0^{-1} \label{eq:bootZ} 
\end{align}

\begin{remark}
The bootstrap assumption \eqref{eq:bootsigma} together with the 
definition 
\eqref{eq:tildega} 
implies that 
\be
\label{eq:bootga}
|\dot \ga(t)|\le \delta^2 (1+t)^{-n+1}
\ee
\end{remark}
\begin{remark}
The bootstrap assumption \eqref{eq:bootZ} together with Lemma 
\ref{le:Sobolev} implies that 
\begin{align}
&\|Z(t)\|_{L^{\infty}}\les \delta C_{0}^{-1} (1+t)^{-\frac n2},
\label{eq:sdecay}\\
&\|Z(t)\|_{H^{s}}\les \delta C_{0}^{-1} \label{eq:senergy}
\end{align}
\end{remark}

  The bootstrap assumption \eqref{eq:bootsigma} strengthens the 
notion of the admissible path.
In particular, it allows us to estimate the deviation between the 
path $x_j(t;\sigma(t))$ corresponding
to the path $\sigma(t)$ and the straight line $x_j(t,\sigma^\infty)$ 
determined by the constant parameter
$\sigma^\infty$ which was defined from $\sigma(t)$ in \eqref{eq:vecnew} and 
\eqref{eq:scalnew}.
This estimate will play an important role in our analysis.

\begin{lemma}
\label{lem:pathdiff}
Let $\sigma(t)$ be an admissible path satisfying the bootstrap 
assumption \eqref{eq:bootsigma}
and let $\sigma^\infty$ be a constant parameter vector as 
in \eqref{eq:vecnew} and \eqref{eq:scalnew}. Then
\be
\label{eq:pathdiff}
|x_j(t;\sigma(t))-tv^\infty_j-D^\infty_j| \les \delta^2 (1+t)^{-n+2}
\ee
\end{lemma}
\begin{proof}
By our choice of $v^\infty_j$ and $D^\infty_j$ one has that
\[ |x_j(t;\sigma(t))-tv^\infty_j-D^\infty_j|\les 
\int_t^\infty\int_s^\infty |\dot{v}_j(\tau)|\,d\tau + \int_t^\infty 
|\dot{D}_j(s)|\,ds \]
and the lemma follows from \eqref{eq:bootsigma}.
\end{proof}

We then have the following corollary. To formulate it, we need the 
localizing functions
\bea
\chi_0(x) &=& \exp\Big(-\half\alpha_{\rm min}(1+|x|^2)^{\half}\Big) 
\nn\\
\chi(t,x;\siginf) &=& \sum_{j=1}^k \chi_0(x-x_j(t;\siginf)). 
\label{eq:loc} 
\eea
Here $\alpha_{\rm min}>0$ satisfies $\inf_{t\ge0,1\le j\le 
k}\alpha_j(t)> \alpha_{\rm min}$ for 
any admissible path $\sigma(t)$ starting at $\sigma_0$. The exponent 
$\alpha_{\rm min}$ arises
because of the decay rate of the ground state of~\eqref{eq:ell}. 

\begin{cor}
\label{cor:Hdiff}
Let $\sigma(t)$ be an admissible path satisfying the bootstrap 
assumption \eqref{eq:bootsigma}.
With the parameters $\sigma^\infty$ as in \eqref{eq:vecnew} and 
\eqref{eq:scalnew} one has
\be
 \Big|H(t,\sigma^\infty)-H(\sigma(t))\Big| \les 
\delta^2(1+t)^{2-n}\,\chi(t,x;\sigma^\infty),
\label{eq:Hdiff}
\ee
where $H(t,\sigma^\infty)$ and $H(\sigma(t))$ are the Hamiltonians 
from~\eqref{eq:Htsigma} and~\eqref{eq:matrix}.
\end{cor}
\begin{proof} The difference
\[ H(t,\sigma^\infty)-H(\sigma(t)) \]
is a sum of matrix valued potentials that are exponentially localized 
around the solitons $w_j(\sigma(t))$
or~$w_j(t,\sigma^\infty)$, respectively. By the previous lemma, we can 
assume that all the potentials are localized
near the straight path~$x(t,\siginf)=(x_1(t;\siginf),\ldots,x_N(t;\sigma^\infty))$. 
Since $v^\infty_j=v_j(\infty),\alpha^\infty_j=\alpha_j(\infty)$, 
\begin{align} 
& \Big|[H(t,\sigma^\infty)-H(\sigma(t))]\Big| \les \sum_{j=1}^k 
\Big|tv^\infty_j+D^\infty_j-
\int_0^t v_j(s)\,ds - D_j(t)\Big|\chi_0(x-x_j(t;\siginf)) 
\label{eq:odin}\\
& + \sum_{j=1}^k \Bigl| \half \int_t^\infty \dot{v}_j(s)\cdot x\,ds 
- \half \int_0^t\int_s^\infty (\dot{v}_j(s)\cdot 
v_j(s)-\dot{\alpha}_j(s)\alpha_j(s))\,ds + \gamma_j - \gamma_j(t) 
\Bigr|\chi_0(x-x_j(t;\siginf)). \label{eq:dwa} 
\end{align}
The term \eqref{eq:odin} arises as the difference of two paths, 
whereas~\eqref{eq:dwa} is the difference
of the phases, i.e.,
\[ |e^{i\theta_j(t,x;\siginf)}-e^{i\theta_j(t,x;\sigma(t))}|.\]
In view of the definitions of $D_j,\gamma_j$ from \eqref{eq:vecnew} 
and \eqref{eq:scalnew} one has
\begin{align}
& \Big|H(t,\siginf)-H(\sigma(t))\Big| \les \sum_{j=1}^k \Bigl( 
\int_t^\infty\int_s^\infty |\dot{v}_j(\tau)|\,d\tau + \int_t^\infty 
|\dot{D}_j(s)|\,ds \Bigr)\chi_0(x-x_j(t,\sigma))  \\
& +  \sum_{j=1}^k \Bigl( \int_t^\infty\int_s^\infty 
|\dot{v}_j(\tau)\cdot 
v_j(\tau)-\dot{\alpha}_j(s)\alpha_j(s)|\,d\tau\,ds + \int_t^\infty 
|\dot{\gamma}_j(s)|\,ds +\int_t^\infty |\dot{v}_j(s)|\,ds|x| 
\Bigr)\chi_0(x-x_j(t,\siginf)) \nn \\
& \les \delta^2(1+t)^{2-n}\,\chi(t,x;\siginf).\nn 
\end{align} 
For the final inequality one uses \eqref{eq:bootga} and the fact that
\[ |x|\chi_0(x-x_j(t;\siginf)) \les t. \]
The corollary follows. 
\end{proof}

\section{Solving the modulation equations}

Our goal is to show that the system in Proposition~\ref{prop:modeq} 
has a solution $\dot{\tilde\sigma}(t)$
that satisfies the bootstrap assumptions~\eqref{eq:bootsigma}. This 
requires some care, as the right-hand
side in Proposition~\ref{prop:modeq} involves the perturbation~$Z$.
We will therefore first verify that the system of modulation 
equations is consistent with the
bootstrap assumptions~\eqref{eq:bootsigma} and~\eqref{eq:bootZ}. In 
what follows, we will use both paths
$\tilde{\sigma}(t)$ and~$\sigma(t)$. By definition, 
see~\eqref{eq:tildega},
\[ {\tilde \gamma}_j(t) = -\int_{t}^\infty \Bigl[\dot{\gamma}_j(s) + 
\half \sum_{m=1}^n \dot{v}^m_j(s)  x^{m}_j(s;\sigma(s)) \Bigr] \, 
ds. \]
The integration is well-defined provided $\tilde{\sigma}$ satisfies 
the bootstrap assumption. Indeed, in that
case $|v_j(t)|\les (1+t)^{-n}$ and since $|x_j(t;\sigma(t))|\les 1+t$, the 
integral is absolutely convergent. 
Finally, recall the property~\eqref{eq:bootga} of the derivatives.

\begin{lemma} 
\label{lem:modul}
Suppose the separation and convexity conditions hold, 
see~\eqref{eq:separat} and~\eqref{eq:stabil_2}. 
Let $\tilde{\sigma}, Z$ be {\em any} choice of functions that satisfy 
the bootstrap assumptions for sufficiently
small $\delta>0$.
If the inhomogeneous terms of the system~\eqref{eq:sigmadoteq}  are 
defined by means
of these functions, then this system has a solution~$\dot{\tilde 
\sigma}$ that satisfies~\eqref{eq:bootsigma} 
with $\delta/2$ for all times.
\end{lemma}
\begin{proof} By the nonlinear 
stability condition~\eqref{eq:stabil_2}, the left-hand side 
of~\eqref{eq:sigmadoteq}
is of the form $B_j(t)\dot{\tilde \sigma}_j(t)$ with an invertible 
matrix~$B_j(t)$. The $O$-term is
a harmless perturbation of the matrix given by the main terms on the 
left-hand side, provided $\delta$
is chosen sufficiently small. This easily follows from the smallness 
of $Z$ 
given by 
\eqref{eq:bootZ}. We need to verify that the right-hand side 
of~\eqref{eq:sigmadoteq}
decays like $\delta^2(1+t)^{-n}$. 
We consider only the first equation in \eqref{eq:sigmadoteq}, the 
others being the same.
The terms $\la V_r(t,\sigma)Z, \xi^1_j(t,\cdot;\sigma(t))\ra $ for $r\ne 
j$ and $w_1w_2$ are governed 
by the interaction of two {\it different} solitons. In view of the 
separation condition~\eqref{eq:separat}
and the exponential localization of the solitons, we have  
\begin{align}
|\dot \alpha_j(t)|\les & e^{-\alpha_{\min}(L+ct)} (1+\|Z(t)\|_{L^{2}+L^\infty}) 
+ \|Z(t)\|^2_{L^{2}+L^\infty}+ \|Z(t)\|^p_{L^{2}+L^\infty} \label{eq:8.1}\\ 
\les & 
\delta C_{0}^{-1}(1+t)^{-n} \Big (\epsilon + \delta C_{0}^{-1} + 
\delta^{p-1} 
C_{0}^{-(p-1)}\Big ) \le \Big (\frac \delta 2\Big )^{2} (1+t)^{-n}
\label{eq:alphaboot}
\end{align}
where we have used the estimate \eqref{eq:sdecay}, the condition 
\eqref{eq:quantsep}, $L\alpha_{\min}\ge |\log \epsilon|$, and that 
$p\ge 2$.
\end{proof}
More generally, the estimates leading up to~\eqref{eq:8.1} also yield 
the following result. The proof is implicit in the preceding one and is
therefore omitted.
\begin{lemma}
\label{lem:siglem}
Let $\sigma(t)$ be an admissible path satisfying the bootstrap assumption 
\eqref{eq:bootsigma} and $Z(t)$ be an arbitrary function in $\Xl_s$. Define 
the function $\dot \Sigma$ as a solution of the equation 
\begin{equation}
\label{eq:Sigmad}
\la \dot\Sigma  W(\sigma(t)), \xi^m_j(t,\cdot;\sigma(t))\ra = \la G(Z(t),\sigma(t)),
\xi^m_j(t,\cdot;\sigma(t))\ra + \la \Omega^m_j(t,\cdot;\sigma(t)), Z(t,\cdot)\ra,
\end{equation}
where 
\begin{align}
&G(Z(t),\sigma(t)) =  O(w_1(\sigma(t))w_2(\sigma(t))Z + O(w_1(\sigma(t)w_2(\sigma(t)) + 
O(|w(\sigma(t))|^{p-2}|Z|^2) + O(|Z|^p),\label{eq:defG}\\
&\Omega_j^m(t,x;\sigma (t))= O\bigl(\dot{ \tilde \sigma} 
(|\phi_j|+|D \phi_j|+|D^2 \phi_j|)\bigr) + 
\sum_{r\ne j} V_r(t,x;\sigma (t)) \xi^m_j(t,x;\sigma(t)).\label{eq:defOmega} 
\end{align}
Then 
\begin{equation}
\label{eq:estSigma}
|\dot\Sigma(t)|\le (1+t)^{-n} \Big (\frac 14\delta^2 + C\|Z\|_{\Xl_s}^2 + C\|Z\|_{\Xl_s}^p\Big )
\end{equation}
\end{lemma}
\begin{remark}
The functions $G(Z(t),\sigma(t))$ and $\Omega_j^m(t,x;\sigma (t))$ arise as follows.
In Section~\ref{sec:exist} we will rewrite the $Z$ equation in the form
\[ i\partial_t Z + H(\sigma(t))Z = \dot\Sigma W(\sigma(t)) + G(Z(t),\sigma(t)),\]
The quantity
$\Omega_j^m(t,x;\sigma (t))$ is defined via the equation
 \[ i\partial_t \xi_j^m(t,\cdot;\sigma (t)) + H^*(\sigma (t)) \xi_j^m(t,\cdot;\sigma (t)) = 
 {\cal S}_k^m \xi_j^k(t,\cdot;\sigma (t)) + \Omega_j^m(t,\cdot;\sigma (t)), \] 
 where the matrix ${\cal S}$ collects the terms $  \xi_j^m(t,\cdot;\sigma (t))$ on the 
 right-hand sides of \eqref{eq:uzhas1}-\eqref{eq:uzhas4}. 
However, the previous proof does not require the explicit form of $G$ or~$\Omega_j^m$.
\end{remark} 

\section{Solving the $Z$ equation}
In this section we verify the bootstrap assumptions \eqref{eq:bootZ} 
for the 
perturbation $Z$. This together with the already verified bootstrap 
estimates for $\dot{\tilde\sigma}$ 
will also lead to the existence of the function $Z(t)$ asserted in 
our main result.
At this point we recall the imposed orthogonality conditions 
\eqref{eq:perp}
\be
\label{eq:perpZ}
\la Z(t), \xi_j^m(t,\cdot;\sigma(t)) \ra  = 0 \text{\ \ for all\ \ } j,m
\ee
with $\xi_j^m(t,x;\sigma(t))$ is as in Definition~\ref{def:movbas}. We 
next rewrite the equation \eqref{eq:system} for $Z$ in  the form
\begin{align}
&i\partial_t Z + H(t,\sigma^\infty) Z = F,\label{eq:equatZ}\\
&F=\big (H(t,\sigma^\infty) - H(\sigma (t))\big )Z + \dot\Sigma W(\sigma(t)) 
+ O(w_1w_2)Z + O(w_1w_2) + O(|w|^{p-2}|Z|^2) + O(|Z|^p)
\label{eq:Fagain}
\end{align}
with the reference Hamiltonian $H(t,\sigma^\infty)$ as defined in 
\eqref{eq:Htsigma}.
To verify the bootstrap assumption \eqref{eq:bootZ} we need to apply the the 
dispersive estimate for the inhomogeneous 
charge transfer problem stated in Theorem \ref{thm:charge}. 
The following lemma shows that the orthogonality conditions 
\eqref{eq:perp} and the bootstrap assumptions \eqref{eq:bootsigma} , 
\eqref{eq:bootZ}
imply that $Z$ is 
asymptotically orthogonal to the bound states of $H_j^*(\siginf)$, as 
required 
in Theorem \ref{thm:charge}.
\begin{lemma}
\label{lem:orthZ}
Let $Z$ be an arbitrary function in $\Xl_s$  satisfying the 
orthogonality conditions \eqref{eq:perpZ} with respect to an admissible path
$\sigma(t)$ obeying the bootstrap assumption \eqref{eq:bootsigma}, and  so that $Z$ verifies
the bootstrap assumption \eqref{eq:bootZ}.
Then $Z$ is asymptotically orthogonal to the null spaces of the 
Hamiltonians $H^*_j(\siginf)$ in the sense 
of \eqref{eq:asorthZ}. In fact,
\be
\label{eq:charor}
\|P_{N_j}(\siginf) \calG_{{v_j^\infty}, {D_j^\infty}}(t) Z(t,\cdot)\|_{L^2} 
\les \delta^3 (1+t)^{-\frac n2-1},
\quad \forall j=1,..,k  
\ee
\end{lemma}
\begin{proof}
By the assumption $Z(t)$ is orthogonal to the vectors 
$\xi_j^m(t,x;\sigma(t))$ introduced in Definition \ref{def:movbas},
while \eqref{eq:charor} is equivalent to  the estimates 
$$
|\la \calG^*_{v^\infty_j,D^\infty_j}(t)\xi_j^m(\cdot;\sigma^\infty), Z(t)\ra |\les \delta^3 
(1+t)^{-\frac n2-1}, \quad \forall j,m
$$
Here $\xi_j^m(x;\sigma^\infty)$, defined in Proposition \ref{prop:null}, refer to the 
elements of the null spaces $\rootsp_j(\sigma^\infty)$ of the stationary Hamiltonians 
$H^*_j(\sigma^\infty)$.
The desired estimate would then follow from the bootstrap assumption 
\eqref{eq:bootZ}, in particular \eqref{eq:sdecay}, and 
the inequality
\be
\label{eq:diffxi}
\|\calG^*_{v^\infty_j,D^\infty_j}(t)\xi_j^m(\cdot;\sigma^\infty) - \xi_j^m(t,\cdot;\sigma(t))\|_{L^1}\les 
\delta^2 (1+t)^{-1}
\ee
The vectors $\xi_j^m$ are composed of the functions derived from the 
bound state $\phi$. In particular,
$\xi_j^1=\binom{\phi}{\phi}$. Therefore,
\begin{align}
|\calG^*_{v^\infty_j,D^\infty_j}(t)\xi_j^1(x;\sigma^\infty) - 
\xi_j^1(t,x;\sigma(t))|= &2 |e^{i(\frac 12 
v^\infty_j\cdot x - \frac 14
(|v^\infty_j|^2-|\alpha^\infty_j|^2) t + \gamma^\infty_j )
} \phi (x-v^\infty_j t - D^\infty_j) -\nn \\ & e^{i(\frac 12 v_j(t)\cdot x - \frac 
14
\int\limits_0^t(|v_j(\tau)|^2-\alpha_j(\tau)^2)\,d\tau + \gamma_j(t) 
)} \phi (x- x_j(t,x;\sigma(t))) |\label{eq:difort}
\end{align}
According to Lemma \ref{lem:pathdiff}, $|x_j(t,x;\sigma(t)) - v^\infty_j t - D^\infty_j|\les 
\delta^2 (1+t)^{-n+2}$. Similarly, \eqref{eq:dwa} of Corollary
\ref{cor:Hdiff} gives the estimate for the difference of the phases 
appearing in \eqref{eq:difort}
$$
|e^{i\theta_j(t,x;\sigma^\infty)} - e^{i\theta_j(t,x;\sigma(t)}|\les 
\delta^2 (1+t)^{-n+2}
$$
The estimate \eqref{eq:diffxi} follows immediately since $n\ge 3$.
\end{proof}
We now in the position to apply Theorem \ref{thm:charge} to 
establish the 
improved $\Xl_{s}$ estimates for $Z(t)$.
\begin{lemma}
\label{lem:bootin}
Let $Z$ be a solution of the equation \eqref{eq:equatZ} satisfying 
the 
bootstrap assumption \eqref{eq:bootZ} with some sufficiently small 
constants $\de$ and $C_{0}^{-1}$. 
We also assume (due to Lemma \ref{lem:modul}) that the admissible 
path $\sigma(t)$ obeys the estimate \eqref{eq:bootsigma}.
Then we have the following  estimate
\be
\label{eq:verdecZ}
\|Z(t)\|_{\Xl_{s}} \le \frac {\delta}2 C_0^{-1} 
\ee
\end{lemma} 
\begin{proof}
Perturbation $Z$ is a solution of the inhomogeneous charge transfer 
problem 
\eqref{eq:equatZ}
\begin{align}
&i \partial_t Z + H(t,\sigma^\infty) Z= F,\nn\\ 
&F:= \big (H(t,\sigma^\infty) - H(\sigma (t))\big )Z + \dot\Sigma 
W(\sigma(t)) + 
O(w_1w_2)Z + O(w_1w_2) + O(|w|^{p-2}|Z|^2) + O(|Z|^p)
\label{eq:defF}
\end{align}
Lemma \ref{lem:orthZ} shows that $Z$ is asymptotically orthogonal 
(with the constant $\delta^{3}$) to 
the null 
spaces of the Hamiltonians $H_j^*(\siginf)$.
Therefore, Theorem \ref{thm:charge} gives the estimate 
\be
\label{eq:repdecayZ}
\|Z(\cdot)\|_{\Xl_{s}} \les \sum_{k=0}^{s}
\|\nabla^{k} Z_0\|_{L^1\cap 
L^2} + 
\|F\|_{\Yl_{s}} + \delta^{3}
\ee
with 
\be
\label{eq:trip}
\| F\|_{\Yl_{s}}=  \sup_{t\ge 0}\Big ( \sum_{k=0}^{s}
\int_{0}^{t}\|\nabla^{k} F(\tau,\cdot)\|_{L^{1}}\,d\tau + 
(1+t)^{\frac n2 +1}
\|F(t,\cdot)\|_{H^{s}}\Big )
\ee
By the assumptions on the initial data 
$\sum_{k=0}^{s}
\|\nabla^{k} Z_0\|_{L^1\cap 
L^2} \le \epsilon <<\delta$. 
Therefore, to obtain the conclusion of 
Lemma \ref{lem:bootin} it would suffice to verify that 
\be
\label{eq:tripF}
\|F\|_{\Yl_{s}} \les \delta^{2}
\ee
with $F$ defined as in \eqref{eq:defF}.
The estimate \eqref{eq:tripF} relies on the following lemma and the bootstrap assumptions
\eqref{eq:bootZ} on $Z(t)$.
\end{proof}
\begin{lemma}\label{lem:FZ}
Let $\sigma(t)$ be an arbitrary admissible path satisfying the bootstrap assumption
\eqref{eq:bootsigma} and $Z(t)$ be an arbitrary function in $\Xl_s$. 
Then the nonlinear expression $F$ defined in terms of the path $\sigma(t)$ and 
$Z(t)$ as in \eqref{eq:defF} obeys the estimate
\begin{equation}
\label{eq:EstF}
\|F\|_{\Yl_s} \les \de^2 +  \|Z\|_{\Xl_s}^2 + \|Z\|_{\Xl_s}^p
\end{equation}
\end{lemma}
\subsection{Algebra estimates}
In this section we establish several simple lemmas designed to 
ease the task of estimating the $\Yl_{s}$ norm of $F=F(Z,w,\sigma)$
in connection with the $\Xl_{s}$ norms of $Z$.

We start by formulating a version of the Sobolev estimate 
tailored to the use of the space $L^{2}+L^{\infty}$.
\begin{lemma}
Let $s$ be a positive integer. Then for any nonnegative
integer $k\le s$ and any $q\in [2,q_{k}]$, where 
\begin{align}
&\frac 1{q_{k}} = \frac 12 - \frac {s-k}n,\qquad {\text{if}}\,\,\,
k> s-\frac n2\label{eq:defq1}\\
&q_{k}=\infty,\qquad {\text{if}}\,\,\,k<s-\frac n2\nn
\end{align}
and $q\in [2,\infty)$ if $k=s-\frac n2$
the following estimates hold true 
\be
\label{eq:sobolev1} 
\|\nabla^{k} f\|_{L^{q}+L^{\infty}} \les 
\sum_{l=0}^{s}\|\nabla^{l} f\|_{L^{2}+L^{\infty}}\le
(1+t)^{-\frac n2} \|f\|_{\Xl_{s}}
\ee
In particular, if $s>\frac n2$
\be
\label{eq:sobolev}
\|f\|_{L^{\infty}} \les 
(1+t)^{-\frac n2} \|f\|_{\Xl_{s}}
\ee
\label{le:Sobolev}
\end{lemma}
\begin{proof}
By duality and density it suffices to show that 
$$
\|f\|_{L^{1}\cap L^{2}}\les \sum_{l=0}^{s-k}\|\nabla^{l} 
f\|_{L^{1}\cap L^{q'}} 
$$
The $L^{1}$ estimate is trivial while the the estimate for the 
$L^{2}$ norm follows from the standard Sobolev embedding 
$W^{k-l,q'}\subset L^{2}$, which holds for the range of parameters
$(k,l,q)$ described in the Lemma.
\end{proof}
Next are the estimates of the nonlinear quantities arising in 
\eqref{eq:defF} in terms of the $\Xl_{s}$ norm.
\begin{lemma}
Let $\gamma(\tau)$ be a smooth function which obeys the estimates 
\be
\label{eq:assumga}
|\gamma^{(\ell)}(\tau)|\les \tau^{(\frac{p-1}2-\ell)_{+}}
\ee 
for some $p\ge 2+\frac 2n$ and all non-negative integers $\ell$. 
Here $r_{+}=r$ if $r\ge 0$ and $r_{+}=0$ if $r<0$.
Then for any $s>\frac n2$ and any non-negative integer 
$k\le s$
\begin{align}
&\|\nabla^{k}\Big (\gamma(|f|^{2})f\Big )\|_{L^{1}}\les 
(1+t)^{-1} (\|f\|_{\Xl_{s}}^{p} + 
\|f\|_{\Xl_{s}}^{2k+1+(p-1-2k)_{+}})\label{eq:L1gamma},\\
&\|\nabla^{k}\Big (\gamma(|f|^{2}) f\Big)\|_{L^{2}}\les 
(1+t)^{-\frac n2 -1} (\|f\|_{\Xl_{s}}^{p} + 
\|f\|_{\Xl_{s}}^{2k+1+(p-1-2k)_{+}})
\label{eq:L2gamma}
\end{align}
In addition, if $\gamma$ is a smooth function obeying 
\eqref{eq:assumga} for some $p\ge 2$ and $\zeta(x)$ is 
an exponentially localized smooth function, then for any 
$q\in [1,2]$
\be
\label{eq:locgamma}
\|\nabla^{k}\Big (\zeta \,\gamma(|f|^{2})f\Big )\|_{L^{q}}\les 
(1+t)^{-n} (\|f\|_{\Xl_{s}}^{p} + 
\|f\|_{\Xl_{s}}^{2k+1+(p-1-2k)_{+}}).
\ee
\label{le:L1leibn}
\end{lemma}
\begin{remark}
\label{re:l1time}
It will become clear from the proof below that if the function $\gamma$ satisfies
\eqref{eq:assumga} for some $p>2+\frac 2n$, then the 
estimate \eqref{eq:L1gamma} holds with a better rate of decay in $t$. In 
particular,
\be
\label{eq:L1tga}
\int_{0}^{\infty} \Big\|\nabla^{k}\Big (\gamma(|f|^{2})f\Big )\Big\|_{L^{1}}dt\les 
\|f\|_{\Xl_{s}}^{p} + 
\|f\|_{\Xl_{s}}^{2k+1}.
\ee
\end{remark}
\begin{proof}
By Leibnitz's rule
$$
\nabla^{k} \Big (\gamma(|f|^{2}) f\Big ) =
\sum_{\ell=0}^{k} \sum_{m_{1}+\ldots+m_{2\ell+1}=k} 
C_{\ell\vec m}\gamma^{(\ell)}(|f|^{2}) \nabla^{m_{1}}f \nabla^{m_{2}} f\ldots
\nabla^{m_{2\ell+1}} f
$$
with some positive integer constants $C_{\ell\vec m}$ and non-negative 
vectors $\vec m=(m_{1},\ldots,m_{2\ell+1})$. We may assume that 
$m_{2\ell+1}\ge m_{2\ell}\ge\ldots\ge m_{1}$. 
Define 
\be
\label{eq:defq}
q_{m_{r}}= (\frac 12 - \frac {s-m_{r}}n)^{-1} 
\ee
for $m_{r}\ge s-\frac n2 $ and $q_{m_{r}}=\infty$ otherwise. 
With the above definition the Sobolev embeddings 
\begin{equation}
\label{eq:sob_imbed}
H^{s}\subset W^{m_{r},q_{m_{r}}}, \qquad W^{s,2}+W^{s,\infty}\subset 
W^{m_{r},q_{m_{r}}}+W^{m_{r},\infty}
\end{equation}
(recall that $m_r\le k\le s$) 
\footnote{In the case of $m_{r}=s-\frac n2$ the value of $q_{m_{r}}$ can 
be set arbitrarily large, but the embedding fails for $q_{m_{r}}=\infty$. 
However, since the following argument has some ``slack'', we can allow ourselves to 
still set $q_{m_{r}}=\infty$ for simplicity.} 
hold true by Lemma \ref{le:Sobolev}. Then
\begin{align}
\|\nabla^{k}& \Big (\gamma(|f|^{2}) f\Big )\|_{L^{1}}\les 
 \|\gamma(|f|^{2}) |\nabla^{k} f|^{1-\frac 2n}\|_{L^{1}\cap L^{\frac 
 n{n-1}}}
\|\nabla^{k} f\|_{L^{2}+L^{\infty}}^{\frac 2n}+ \nn \\
&\sum_{\ell=1}^{k} \sum_{m_{1}+\ldots+m_{2\ell+1}=k} 
\|\gamma^{(l)}(|f|^{2}) \nabla^{m_{1}}f\ldots 
\nabla^{m_{2\ell}}f\|_{L^{1}\cap L^{q'_{m_{2\ell+1}}}}
\big (\|f\|_{L^{2}+L^{\infty}}+
\|\nabla^{s} f\|_{L^{2}+L^{\infty}}\big ),\label{eq:leibn}
\end{align}
where the final term inside the parentheses in~\eqref{eq:leibn} arises
via the second embedding in~\eqref{eq:sob_imbed}. 
We claim for $\ell>0$ that there exist two sets of parameters 
$q^{1}_{m_{r}}$ and $q^{2}_{m_{r}}$ for $r=1,\ldots,2\ell $ such that 
\begin{align}
&2\le q^{1,2}_{m_{r}}\le q_{m_{r}}, \qquad \forall r=1,\ldots,2\ell,
\label{eq:subcrit}\\
&\sum_{r=1}^{2\ell} \frac 1{q^{1}_{m_{r}}} =1,\qquad
\sum_{r=1}^{2\ell} \frac 1{q^{2}_{m_{r}}} =\frac 1{q'_{m_{2\ell+1}}}.\nn
\end{align}
To prove the claim we 
let $\tau$ be 
the number of $m_{r}$ for $r=1,\ldots,2\ell$ such that $m_{r}\ge s-\frac n2$. 
Observe that 
$$
\sum_{r=1:\, m_{r}\ge s-\frac n2}^{2\ell} m_{r} \le k-m_{2\ell+1}
$$
Therefore,
\begin{align*}
\sum_{r=1}^{2\ell} \frac 1{q_{m_{r}}} &\le  
\frac \tau 2 - \frac {\tau s-k+m_{2\ell+1}}n \le
 - \tau (\frac {s}n-\frac 12) + 
\frac {k-m_{2\ell+1}}n\\ &\le 
- \tau (\frac {s}n-\frac 12) - \frac {s-k}n + \frac 12
\le \frac 12 \le \frac 1{q'_{m_{2\ell+1}}}.
\end{align*}
The inequality in the second line above follows since 
$m_{2\ell+1} \ge s-\frac n2$, which holds if 
$\tau > 0$. 
On the other hand,
$$
\sum_{r=1}^{2\ell} \frac 12 =  \ell \ge 1
$$
and the claim immediately follows, provided that $\ell >0$.
Thus using the sequence $q^1_{m_{r}}$ to handle the $L^{1}$
norm in \eqref{eq:leibn} and $q^{2}_{m_{r}}$ for the 
$L^{q'_{m_{2\ell+1}}}$ norm, we obtain 
\begin{align*}
& \|\nabla^{k} \Big (\gamma(|f|^{2}) f\Big )\|_{L^{1}}\les 
\|\gamma(|f|^{2}) |\nabla^{k} f|^{1-\frac 2n}\|_{L^{1}\cap 
L^{\frac n{n-1} }}
\|\nabla^{k} f\|_{L^{2}+L^{\infty}}^{\frac 2n}+\\
&\sum_{\ell=1}^{k} \sum_{m_{1}+\ldots+m_{2\ell+1}=k} 
\|\gamma^{(\ell)}(|f|^{2}) \|_{L^{\infty}}
\|\nabla^{m_{1}}f\|_{L^{q^{1}_{m_{1}}}\cap L^{q^{2}_{m_{1}}}} 
\ldots
 \|\nabla^{m_{2\ell}} f\|_{L^{q^{1}_{m_{2\ell}}}\cap 
L^{q^{2}_{m_{2\ell}}}} 
\|\nabla^{s} f\|_{L^{2}+L^{\infty}}.
\end{align*}
By H\"older's inequality
$$
\|\gamma(|f|^{2}) |\nabla^{k} f|^{1-\frac 2n}\|_{L^{1}\cap L^{\frac n{n-1}}}
\les  \|\nabla^{k} f\|_{L^{2}}^{1-\frac 2n} 
\|f\|^{p-1}_{L^{(p-1)\frac {2n}{n+2}}\cap L^{2(p-1)}}
\les \|f\|_{H^{s}}^{p-\frac 2n}
$$
provided that $p\ge 2+\frac 2n$, which is dictated by the condition
that $(p-1)\frac {2n}{n+2}\ge 2$.
Finally, using the property \eqref{eq:subcrit} together with the
estimate~\eqref{eq:assumga} we obtain
\begin{align*}
\|\nabla^{k} \Big (\gamma(|f|^{2}) f\Big )\|_{L^{1}}&\les 
\|f\|^{p-\frac 2n}_{H^{s}} \|\nabla^{k} f\|^{\frac 2n}_{L^{2}+L^{\infty}}+ 
(\|f\|^{p-1}_{H^{s}} + \|f\|_{H^{s}}^{2k+(p-1-2k)_{+}})
\|\nabla^{s} f\|_{L^{2}+L^{\infty}}
\\ &\les t^{-1} \Big(\|f\|_{\Xl_{s}}^{p} + 
\|f\|_{\Xl_{s}}^{2k+1-(p-1-2k)_{+}}\Big). 
\end{align*}
Similarly, we estimate 
\begin{align}
\|\nabla^{k}& \Big (\gamma(|f|^{2}) f\Big )\|_{L^{2}}\les 
\|\gamma(|f|^{2})\|_{L^{2}\cap L^{\infty}}
\|\nabla^{k} f\|_{L^{2}+L^{\infty}}+\label{eq:leibnl2}\\
&\sum_{\ell=1}^{k} \sum_{m_{1}+\ldots+m_{2\ell+1}=k} 
\|\gamma^{(\ell)}(|f|^{2}) \nabla^{m_{1}}f\ldots\nabla^{m_{2\ell}}f\|
_{L^2\cap L^{\frac{2q_{m_{2\ell+1}}}{q_{m_{2\ell+1}}-2}}} 
\big (\|f\|_{L^{2}+L^{\infty}}+\|\nabla^{s}f\|_{L^{2}+L^{\infty}}\big )\nn
\end{align}
To estimate the first term in \eqref{eq:leibnl2} we note that
$$
\|\gamma(|f|^{2})\|_{L^{2}\cap L^{\infty}}\le 
\|f\|^{p-1}_{L^{\infty}\cap L^{2(p-1)}}\le 
\|f\|_{L^{\infty}}^{{p-2}}\|f\|_{L^{2}\cap L^{\infty}}
\les (1+t)^{-\frac n2(p-2)} \|f\|_{\Xl_{s}}^{p-1}, 
$$
where the second inequality follows from interpolating 
$L^{2(p-1)}$ between $L^{2}$ and $L^{\infty}$ and 
the last inequality is a consequence of Lemma~\ref{le:Sobolev}
and the definition of the space $\Xl_{s}$. Thus
\be
\label{eq:leibnfirst}
\|\gamma(|f|^{2})\|_{L^{2}\cap L^{\infty}}
\|\nabla^{k} f\|_{L^{2}+L^{\infty}}\les 
(1+t)^{-\frac n2(p-1)} \|f\|_{\Xl_{s}}^{p}. 
\ee
Furthermore, using definition~\eqref{eq:defq} we have that 
$$
\frac{2q_{m_{2\ell+1}}}{q_{m_{2\ell+1}}-2} = \frac n{s-m_{2\ell+1}}.
$$
Then
\begin{align}
\|\gamma^{(\ell)}(|f|^{2}) \nabla^{m_{1}}f\ldots\nabla^{m_{2\ell}}f\|&
_{L^2\cap L^{\frac n{s-m_{2l+1}}}}\les 
\|\nabla^{m_{1}}f\|_{L^{q_{m_{1}}}+L^{\infty}}\times\label{eq:leibnint}\\
&\|\gamma^{(\ell)}(|f|^{2}) 
\nabla^{m_{2}}f\ldots \nabla^{m_{2\ell}}f\|_{L^{2}\cap 
L^{\frac n{s-m_{2\ell+1}}}
\cap L^{\frac {2q_{m_{1}}}{q_{m_{1}}-2}}\cap L^{\frac 
{nq_{m_{1}}}{q_{m_{1}}(s-m_{2\ell+1})-n}}}\nn
\end{align}
Using the definition of $q_{m_1}$ from \eqref{eq:defq} and the
assumption that $m_{1}\le m_{2\ell+1}$ we infer that the last norm 
reduces to the one of the space 
\begin{align*}
&L^{2}\cap L^{\frac n{s-m_{2\ell+1}}},\qquad {\text{for}}
\,\,\, m_{1}\le s-\frac n2,\\
&L^{2}\cap L^{\frac n{s-m_{2\ell+1}+ (s-m_{1}-\frac n2)}},\qquad
{\text{for}}\,\,\, m_{1}> s-\frac n2
\end{align*}
We now
let $\tau$ be 
the number of $m_{r}$ for $r=2,..,2\ell$ such that $m_{r}\ge s-\frac n2$. 
Observe that since $s> \frac n2$ and $s\ge k$
\begin{align*} 
\sum_{r=2}^{2\ell} \frac 1{q_{m_{r}}} &\le  
\frac \tau 2 - \frac {s \tau-k+m_{2\ell+1}+m_{1}}n \\ &=
- \tau(\frac {s}n-\frac 12) + \frac {k-m_{2\ell+1}-m_{1}}n\\ &\le 
\min\{\frac {s-m_{2\ell+1}}n, \frac {(s-m_{2\ell+1}-m_{1})+s-\frac n2}n\}
\end{align*}
On the other hand
$$
\sum_{r=2}^{2\ell} \frac 12 =  \ell-\frac 12 \ge \frac 12.
$$
It therefore follows that there exist 2 sets of parameters 
$q_{m_{r}}^{1}$ and $q_{m_{r}}^{2}$ for $r=1,\ldots,2\ell$ such that
\begin{align}
&2\le q^{1,2}_{m_{r}}\le q_{m_{r}}, \qquad \forall r=2,\ldots,2\ell
\label{eq:subcrit2}\\
&\sum_{r=2}^{2\ell} \frac 1{q^{1}_{m_{r}}} =2,\nn\\
&\sum_{r=2}^{2\ell} \frac 1{q^{2}_{m_{r}}} = \frac 
{(s-m_{2\ell+1}-m_{1})+s-\frac n2}n\,\,\,{\text{or}}\,\,\, \frac 
{s-m_{2\ell+1}}n.\nn
\end{align}
In either case, with the help of Lemma \ref{le:Sobolev}, we can estimate
$$
\|\gamma^{(\ell)}(|f|^{2}) 
\nabla^{m_{2}}f\ldots \nabla^{m_{2\ell}}f\|_{L^{2}\cap 
L^{\frac n{s-m_{2\ell+1}}}
\cap L^{\frac {2q_{m_{1}}}{q_{m_{1}}-2}}\cap L^{\frac 
{nq_{m_{1}}}{q_{m_{1}}(s-m_{2\ell+1})-n}}}\les 
\|f\|_{H^{s}}^{(p-1-2\ell)_{+}+2\ell-1}
$$
It therefore follows that the second term in \eqref{eq:leibnl2} 
is
\be
\label{eq:leibnsec}
\les (1+t)^{-n} \sum_{\ell=1}^{k}\|f\|_{\Xl_{s}}^{(p-1-2\ell)_{+}+2\ell+1}.
\ee
Now combining this with 
\eqref{eq:leibnfirst}, and using the 
condition that $p\ge 2+\frac 2n$, we infer that 
$$
\|\nabla^{k}\Big (\gamma(|f|^{2}) f\Big)\|_{L^{2}}\les 
(1+t)^{-\frac n2 -1} (\|f\|_{\Xl_{s}}^{p} + 
\|f\|_{\Xl_{s}}^{2k+1+(p-1-2k)_{+}}).
$$
The proof of \eqref{eq:locgamma} proceeds along the lines of 
the argument  for the $L^{2}$ estimate \eqref{eq:L2gamma}. 
We first observe that since $\zeta(x)$ is an exponentially 
localized function, the  $L^{q}$ estimate for $1\le q\le 2$ 
can be reduced to the $L^{2}$ estimate. We then note that the
condition that $p\ge 2+\frac 2n$ was only used in the estimate
\eqref{eq:leibnfirst} which now takes the form
\begin{align*}
\|\gamma(|f|^{2})\,\zeta\|_{L^{2}\cap L^{\infty}}
\|\nabla^{k} f\|_{L^{2}+L^{\infty}}&\les 
\|\gamma(|f|^{2})\|_{L^{\infty}}
\|\nabla^{k} f\|_{L^{2}+L^{\infty}}\\ &\les 
\||f|^{p-1}\|_{L^{\infty}}
\|\nabla^{k} f\|_{L^{2}+L^{\infty}}\\
&\les (1+t)^{-\frac n2 p} \|f\|_{\Xl_{s}}^{p}
\end{align*}
The remaining estimates already have the desired form 
\eqref{eq:leibnsec}. 
\end{proof}

\subsection{$L^{1}$ estimates}
In this section and the following we prove Lemma~\ref{lem:FZ}.
We start with the verification of 
$$
\sum_{k=0}^{s}\int_{0}^{\infty} \|\nabla^{k} F(t,\cdot)\|_{L^{1}}\,dt 
\les \delta^{2} + \|Z\|_{\Xl_s}^2 + \|Z\|_{\Xl_s}^p 
$$
with $F$ as in \eqref{eq:defF}. By Corollary \ref{cor:Hdiff} we have
\be
\label{eq:difH2}
 \Big|H(t,\siginf)-H(\sigma(t))\Big| \les 
\delta^2(1+t)^{2-n}\,\chi(t,x;\siginf),
\ee
where $\chi(t,x;\siginf)$ is a smooth cut-off function localized 
around the union of the paths $x_j(t;\siginf)=v^\infty_j t + D^\infty_j$.
Moreover, the spatial derivatives of the above difference also 
satisfy the same estimates. Using the bootstrap assumptions \eqref{eq:bootZ}  we obtain
\begin{align}
\sum_{k=0}^{s}\int_0^\infty \Big\|\nabla^{k} \Big (\Big [H(\tau,\siginf)-H(\sigma(\tau))\Big ] 
Z(\tau)\Big )\Big\|_{L^1}\,d\tau &\les 
\delta^{2} \sum_{k=0}^{s} \int_{0}^{\infty} \|\nabla^{k} 
Z(\tau)\|_{L^{2}+L^{\infty}} (1+\tau)^{2-n}\,d\tau\nn \\ &\les
\delta^2 \|Z\|_{\Xl_s} \int_0^\infty 
(1+\tau)^{2-n-\frac n2}\,d\tau \les 
\delta^2 \|Z\|_{\Xl_s}. \label{eq:Hdiffest2}
\end{align}
The term $\dot\Sigma W(\sigma(t))$ obeys the point wise bound 
$$
|\dot\Sigma W(\sigma(t)) |\les \max_j |\dot{\tilde \sigma}_j(t)| 
\chi(t,x;\siginf)
$$
This can be easily seen from the equation \eqref{eq:dotS} and Lemma 
\ref{lem:pathdiff}. 
The same estimate also holds for the spatial derivatives of the 
quantity above.
Thus, with the help of the already
verified estimate \eqref{eq:bootsigma} we infer that
\be
\label{eq:sigmaest2}
\sum_{k=0}^{s}\int_0^\infty 
\Big\|\nabla^{k} \Big (\dot\Sigma W(\sigma(\tau))\Big )\Big\|_{L^1}\,
d\tau \les \delta^2 \int_0^\infty (1+\tau)^{-n}\, d\tau\les \delta^2 
\ee
The estimates for the $O(w_1w_2) Z$ and 
$O(w_1w_2)$ terms in \eqref{eq:defF} are straightforward 
due to the separation 
and the exponential localization of the solitons $w_1$ and $w_2$, e.g.,
\be
\label{eq:intsol2}  
\sum_{k=0}^{s}\int_0^\infty \|O\Big (\nabla^{k}( w_1w_2)\Big )\|_{L^{1}}\les 
\int_0^\infty e^{-\alpha_{\min} 
(L+c\tau)}\,d\tau 
\le \frac {e^{-\alpha_{\min} L}}{c\alpha_{\min} }\les
\frac {\epsilon}{\alpha_{\min}} \les \delta^2
\ee
Here we have used the separation assumption \eqref{eq:separat} and 
the condition \eqref{eq:quantsep}, 
$\alpha_{\min} L\ge |\log \epsilon |$.

The  exponential localization of the multi-soliton state $w$, the 
bootstrap assumptions \eqref{eq:bootZ} and the estimate 
\eqref{eq:locgamma} of Lemma \ref{le:L1leibn} yield the estimate
\be
\label{eq:quadrest2}
\sum_{k=0}^{s}\int_0^\infty \| O\Big (\nabla^{k}(|w|^{p-2} Z^2)\Big 
)\|_{L^1}\,d\tau \les \|Z\|_{\Xl_s}^2 \int_0^\infty (1+\tau)^{-n}\,d\tau \les \|Z\|_{\Xl_s}^2. 
\ee
Finally, with the help of \eqref{eq:L1gamma} (more specifically 
using the improvement \eqref{eq:L1tga} of Remark~\ref{re:l1time}), we obtain
\be
\label{eq:Zpest2}
\sum_{k=0}^{s}
\int_0^\infty\|\nabla^{k} (Z^p(\tau))\|_{L^1}\,d\tau \les 
\|Z\|_{\Xl_{s}}^{p}
\ee

\subsection {$L^2$ estimates}
In this subsection we establish the estimate 
$$
\|F(t,\cdot)\|_{H^{s}}\les (1+t)^{-\frac n2-1}\Big(\delta^{2} + \|Z\|_{\Xl_s}^2 + \|Z\|_{\Xl_s}^p \Big).
$$ 
The arguments follows closely those of the previous section.
Using the estimates~\eqref{eq:difH2}, \eqref{eq:locgamma} and  the 
bootstrap assumptions~\eqref{eq:bootZ}  we obtain
\begin{align}
\Big\|\Big (H(\tau,\siginf)-H(\sigma(\tau))\Big )Z(t)\Big\|_{H^{s}}\,
& \les \delta^2 (1+t)^{2-n} \sum_{k=0}^{s} 
 \|\nabla^{k}Z(t,\cdot)\|_{L^{2}+L^{\infty}}\nn \\ &\les 
\delta^2 \|Z\|_{\Xl_s}\,(1+t)^{2-n-\frac n2} \les \delta^{2}\|Z\|_{\Xl_s}\, (1+t)^{-\frac n2-1}
\label{eq:Hdiffest}
\end{align}
where the last inequality follows since $n\ge 3$. Similar to \eqref{eq:sigmaest2}
\be
\label{eq:sigmaest}
\|\dot\Sigma W(\sigma(t))\|_{H^s} \les \delta^2 (1+t)^{-n} \les
\delta^2 (1+t)^{-\frac n2-1} 
\ee
The estimates for the $O(w_1w_2) Z$ and $O(w_1w_2)$ terms again
follow from the separation 
and the exponential localization of the solitons $w_1$ and $w_2$, 
\be
\label{eq:intsol}  
\|O(w_1w_2)\|_{H^{s}}\les 
e^{-\alpha_{\min}(L+ct)}\les \delta^{2}(1+t)^{-\frac n2-1}.
\ee
The exponential localization of the multi-soliton state $w$ together 
with the estimate \eqref{eq:locgamma} of Lemma \ref{le:L1leibn} and the 
bootstrap assumption \eqref{eq:bootZ}, also give the estimate
\be
\label{eq:quadrest}
\| O(|w|^{p-2} Z^2)\|_{H^{s}}\les \|Z\|_{\Xl_s}^2  
(1+t)^{-n}\les \|Z\|_{\Xl_s}^2 
(1+t)^{-\frac n2-1}.
\ee
Finally, using the estimate \eqref{eq:L2gamma} of Lemma 
\ref{le:L1leibn},
we obtain
\be
\|Z^p(t)\|_{H^s} \les (1+t)^{-\frac n2-1} \|Z\|_{\Xl_{s}}^{p}. 
\label{eq:Zpest}
\ee
This completes the proof of Lemma~\ref{lem:FZ}. 
\section{Scattering}
In this section we intend to prove the last part of Theorem \ref{thm:main1}.
More precisely, we shall show that 
there exists $u_0\in L^2$ so that
\[ \Big\|\psi(t,\cdot)-\sum_{j=1}^N w_j(t,x;\sigma^\infty) - e^{i\frac{t}{2}\Laplace} u_0 \Big\|_{L^2} 
\to0\]
as $t\to\infty$. 
Observe that the equation \eqref{eq:main} for $Z=(R,\bar R)$ and the 
Corollary \ref{cor:Hdiff} allow us to write the solution 
$\psi(t)$ in the form 
$$
\psi(t)=\sum_{j=1}^N w_j(t,x;\sigma^\infty) + R(t,x) +Q(t,x),
$$
where the function $Q$ is spatially localized and decays in time with the rate 
$t^{-n+2}$ and thus $\|Q(t)\|_{L^2_x}\to 0$ as $t\to\infty$. The function 
$R$ verifies an inhomogeneous Schr\"odinger  equation 
$$
i\partial_t R + \frac 12 \Delta R = V(t,x) R + F, \qquad R(0,x)=R_0(x)
$$
with spatially exponentially localized
 potential $V$. According to Lemma \ref{lem:bootin} the inhomogeneous term $F$ satisfies the
  estimate 
 \begin{equation}\label{eq:bound-F}
  \|F(t,\cdot)\|_{L^2}\les (1+t)^{-\frac n2-1} \delta^2,
  \end{equation}
while the solution $R$ has been shown to obey the 
dispersive estimate
\begin{equation}\label{eq:bound-R}
\|R(t,\cdot)\|_{L^2+L^\infty}\les (1+t)^{-\frac n2} \delta.
\end{equation}
The standard scattering theory argument show that the desired function 
$u_0$ can be constructed as the limit 
$$
u_0= R_0 - i \lim_{t\to\infty}\int_0^t e^{-i\frac s2 \Delta } 
\big (V(s,\cdot ) R(s,\cdot ) + F(s,\cdot )\big )\, ds.
$$
The estimates \eqref{eq:bound-F} and \eqref{eq:bound-R} together
with the localization of $V$ guarantee the existence of the limit.
\section{Existence}
\label{sec:exist}
In Lemmas \ref{lem:modul} and \ref{lem:bootin} we established the 
estimates 
\be
\label{eq:siges}
|\dot {\tilde\sigma}|\le \frac 14 \delta^{2} (1+t)^{-n}
\ee
for the admissible path $\sigma(t)$ and 
\be
\label{eq:Zest}
\|Z\|_{\Xl_{s}}\le \delta^{2}
\ee
for the solution $Z(t,x)$ of 
the nonlinear inhomogeneous matrix charge transfer problem~\eqref{eq:equatZ},
under the bootstrap assumptions \eqref{eq:bootsigma}, \eqref{eq:bootZ}
\begin{align}
&|\dot {\tilde\sigma}|\le  \delta^{2} (1+t)^{-n}\label{eq:bos},\\
&\|Z\|_{\Xl_{s}}\le \delta C_{0}^{-1}\label{eq:boZ}
\end{align}
and the condition that $Z$ is asymptotically orthogonal to the 
null spaces of the Hamiltonians $H_{j}^{*}(\sigma)$ with the constant
$\delta^{3}$. In this section we shall show that these are sufficient 
to establish the existence of the desired admissible path and the 
perturbation $R$. We prove existence by iteration.
We shall define a sequence of admissible paths $\sigma^{(n)}(t)$ and 
approximate solutions $Z^{(n)}(t)$ for $n=1,\ldots$ according to the 
following rules.
First, we write the $Z$ equation~\eqref{eq:equatZ} and~\eqref{eq:Fagain}
in the form
\begin{equation}
\label{eq:Zequ1} i\partial_t Z + H(\sigma(t))Z= \dot\Sigma W(\sigma(t)) + G(Z,\sigma(t)) 
\end{equation}
where $G(Z,\sigma(t)) = O(w_1w_2)Z + O(w_1w_2) + O(|w|^{p-2}|Z|^2) + O(|Z|^p)$. 
Set 
$$
\sigma^{(1)}(t)=\sigma(0),\qquad Z^{(0)}\equiv 0
$$
where $\sigma^{(1)}$ is to be understood as 
the constant path coinciding with the initial data $\sigma(0)$
common to all admissible paths. 
We now define functions $Z^{1}(t,x)$ and $\sigma^{2}(t)$ to be a solution of the 
following {\textit{linear}} system
\begin{align}
&i \pa_{t} Z^{(1)} + H(\sigma^{(1)}(t)) Z^{(1)} = \dot\Sigma^{(2)} 
W(\sigma^{(1)}(t)) + G(Z^{(0)},\sigma^{(1)}(t)),\label{eq:Z1s} \\
&Z^{(1)}(0,x) = Z_{0}(x),\nn\\
&\Bigl \la \dot\Sigma^{(2)} W(\sigma^{(1)}(t)),\xi_j^m(t,\cdot;\sigma^{(1)}(t))\Bigr \ra = 
\Bigl \la G(Z^{(0)},\sigma^{(1)}(t)), 
\xi_j^m(t,\cdot;\sigma^{(1)}(t)) \Bigr \ra + \Big\la \Omega_j^m(t,\cdot;\sigma^{(1)}(t)), Z^{(1)} \Big\ra.
 \label{eq:sigma2}
\end{align}
Here $\Omega_j^m(t,x;\sigma^{(1)}(t))$ is defined via the equation
\[ i\partial_t \xi_j^m(t,\cdot;\sigma^{(1)}(t)) + H^*(\sigma^{(1)}(t)) \xi_j^m(t,\cdot;\sigma^{(1)}(t)) = 
{\cal S}_k^m \xi_j^k(t,\cdot;\sigma^{(1)}(t)) + \Omega_j^m(t,\cdot;\sigma^{(1)}(t)), \]
where the matrix ${\cal S}$ collects the terms $  \xi_j^m(t,\cdot;\sigma^{(1)}(t))$ on the 
right-hand sides of \eqref{eq:uzhas1}-\eqref{eq:uzhas4}.  
Thus, using Propositions \ref{prop:movingxi} and \ref{prop:modeq}, we have  
\[ \Omega_j^m(t,x;\sigma^{(1)}(t))= O\bigl(\dot{ \tilde \sigma}^{(1)} 
(|\phi_j|+|D \phi_j|+|D^2 \phi_j|)\bigr) + 
\sum_{r\ne j} V_r(t,x;\sigma^{(1)}(t)) \xi^m_j(t,x;\sigma^{(1)}(t)). 
\]
Observe that \eqref{eq:Z1s} arises from the nonlinear equation \eqref{eq:Zequ1} by 
replacing  $\sigma(t)$   with the already defined path  $\sigma^{(1)}(t)$ 
as well as $Z$ on the right-hand side with~$Z^{(0)}=0$. 
The equation \eqref{eq:sigma2} determining $\sigma^{(2)}(t)$ 
ensures that $\la Z^{(1)}(t),\xi_{j}^{m}(t,\cdot;\sigma^{(1)}(t))\ra = 0$.
Indeed, taking scalar products of~\eqref{eq:Z1s} with~$\xi_j^m(t,\cdot;\sigma^{(1)}(t))$ 
and using Propositions~\ref{prop:movingxi} and \ref{prop:modeq} yields
\[ \frac{d}{dt} \Xi(t;\sigma^{(1)}(t)) + {\mathcal S} \Xi(t;\sigma^{(1)}(t)) = 0,\]
where $\Xi(t;\sigma^{(1)}(t))$ is the vector of 
$\la \xi_{j}^{m}(t,\cdot;\sigma^{(1)}(t)), Z^{(1)}(t,\cdot)\ra$
and where $\cal S$ is the constant matrix defined above. 
As long as 
\begin{equation}
\label{eq:Xi_init}
\Xi(0;\sigma^{(1)}(0))=\Xi(0;\sigma(0))=0
\end{equation}
one therefore has $\Xi(t;\sigma^{(1)}(t))=0$ for all $t\ge0$. 
Generally speaking, \eqref{eq:Xi_init} need not be satisfied.
However, using the fact that the initial perturbation~$R_0$ is small,
we proceed as in~\cite{BP1} Proposition~1.3.1 to show that one can 
modify the initial splitting in such a way that it does hold. More precisely,
one has the following lemma.

\begin{lemma}
\label{lem:init_orth}
Let $\psi(0)=\sum_{j=1}^N w_j(0,x;\sigma(0)) + R_0$ with $R_0$ small as in~\eqref{eq:smalldata}
and assume either one of the separation conditions~\eqref{eq:separat} or~\eqref{eq:vel_sep}
as well as the convexity condition~\eqref{eq:stabil_2} in a small neighborhood $\calU$ of $\sigma(0)$.
Then there exist $\wt{\sigma(0)}\in\calU$ such that in the decomposition
\[ \psi(0)=\sum_{j=1}^N w_j(0,x;\wt{\sigma(0)}) + \wt{R_0} \]
the new perturbation 
$\left( \begin{array}{c} \wt{R_0} \\  \overline{\wt{R_0}}
 \end{array} \right)
$ 
is orthogonal to the root spaces $\rootsp_j^*$ of $H_j^*(\wt{\sigma(0)})$ and satisfies the 
smallness condition \eqref{eq:smalldata}.
\end{lemma}
\begin{proof} We need to solve the equation (with the solution being $\wt{\sigma(0)}$) 
\begin{equation}
\label{eq:init_cond}
\Big\la \Psi(0)-\sum_{j=1}^N W_j(0,\cdot;\sigma(0)), J \partial_{\sigma_r} W_\ell(0,\cdot;\sigma(0))
 \Big\ra = 0 \text{\ \ for all\ \ }1\le r\le 2n+2,
\end{equation}
where $J$ is the matrix from Proposition~\ref{prop:null} and $\Psi$, $W_j$
are the complexified versions of $\psi,w_j$, respectively. 
One solves~\eqref{eq:init_cond} by means of the implicit function theorem.
Indeed, the derivative of the left-hand side of~\eqref{eq:init_cond} is given by
\[ 
\Big\la -\sum_{j=1}^N \partial_{\sigma_k} W_j(0,\cdot;\sigma(0)), J \partial_{\sigma_r} 
W_\ell(0,\cdot;\sigma(0)) \Big\ra + \Big\la \Psi(0)-\sum_{j=1}^N W_j(0,\cdot;\sigma(0)), 
J \partial^2_{\sigma_r\sigma_k} W_\ell(0,\cdot;\sigma(0))
 \Big\ra.
\]
The second term is $O(\eps)$ where $\eps$ controls the size of the
initial perturbation~$R_0$ in $L^1$, say. On the other hand, the first
term is separated from zero by virtue of the convexity condition, and 
either the separation condition~\eqref{eq:separat}
or the assumption of large relative velocities of the solitons~\eqref{eq:vel_sep}. 
Indeed,  as in Proposition~1.3.1 from~\cite{BP1} one sees that 
\begin{equation}
\label{eq:det_BP}
 \Big| \det \Big\{ \Big\la  \partial_{\sigma_k} W_\ell(0,\cdot;\sigma(0)), J \partial_{\sigma_r} W_\ell(0,\cdot;\sigma(0)) \Big\ra \Big\}_{1\le k,\ell\le 2n+2} \Big| \gtrsim \|\phi(\cdot;\alpha_\ell(0))\|_2^4
\big(\partial_\alpha \| \phi(\cdot;\alpha_\ell(0))\|_2^2 \big)^2 >0
\end{equation}  
uniformly in the small neighborhood of $\alpha_\ell(0)$ that we are allowing, see the convexity
condition~\eqref{eq:stabil_2}. The remaining entries of the derivative matrix, which
involve inner products with $W_j, W_k$ for $j\ne k$, are small because of either the 
(physical)  separation condition~\eqref{eq:separat} or the velocity condition~\eqref{eq:vel_sep}.
The latter ensures that we are taking scalar products of quantities that are
almost orthogonal by virtue of the large distances of their Fourier transforms. 
Hence the determinant of the derivative is essentially bounded below by the product of
the matrices with $j=k$, see~\eqref{eq:det_BP}. This proves that the derivatives are invertible,
and since the original perturbation~$R_0$ is small, the image of the diffeomorphism given by the
left-hand side of~\eqref{eq:init_cond} contains zero, as claimed. 
\end{proof}
In general, we define
\begin{align}
&i \pa_{t} Z^{(n)} + H(\sigma^{(n)}(t)) Z^{(n)}  = \dot\Sigma^{(n+1)} 
W(\sigma^{(n)}(t)) + G(Z^{(n-1)}, \sigma^{(n)}(t))  \label{eq:equatA}\\
&Z^{(n)}(0,x) =Z_{ 0}(x)\nn\\
&\Bigl \la \dot\Sigma^{(n+1)} W(\sigma^{(n)}(t)),\xi_j^m(t,\cdot;\sigma^{(n)}(t))\Bigr \ra =
 \Bigl \la G(Z^{(n-1)}, \sigma^{(n)}(t)) , 
\xi_j^m(t,\cdot;\sigma^{(n)}(t)) \Bigr \ra + \Big\la \Omega_j^m(t,\cdot; \sigma^{(n)}(t)),
 Z^{(n)}(t) \Big \ra. 
\label{eq:orthn}
\end{align}
Here 
$$
G(Z^{(n-1)}, \sigma^{(n)}(t))=O(w_{1}^{(n)} w_{2}^{(n)}) + 
O(w^{(n)}_1w^{(n)}_2)Z^{(n-1)}+ O(|w^{(n)}|^{p-2}|Z^{(n-1)}|^2) + 
O(|Z^{(n-1)}|^p)
$$  and 
$\Omega_j^m(t,x;\sigma^{(n)}(t))$ once again is defined via the equation
\[ i\partial_t \xi_j^m(t,\cdot;\sigma^{(n)}(t)) + H^*(\sigma^{(n)}(t)) \xi_j^m(t,\cdot;\sigma^{(n)}(t)) = 
{\cal S}_k^m \xi_j^k(t,\cdot;\sigma^{(n)}(t)) + \Omega_j^m(t,\cdot;\sigma^{(n)}(t)), \]
and has the form
\[ \Omega_j^m(t,x;\sigma^{(n)}(t))= O\bigl(\dot{ \tilde \sigma}^{(n)} 
(|\phi_j|+|D \phi_j|+|D^2 \phi_j|)\bigr) + 
\sum_{r\ne j} V_r(t,x;\sigma^{(n)}(t)) \xi^m_j(t,x;\sigma^{(n)}(t)). 
\]
Observe that by the same argument  as in the case of  $Z^{(1)}$, the perturbation 
$Z^{(n)}$ is orthogonal 
to the functions $\xi^m_j(t,x;\sigma^{(n)}(t))$. 
We shall assume that solutions $Z^{(n)}, \sigma^{(n)}$ of \eqref{eq:equatA}, \eqref{eq:orthn}
have already been constructed and we now proceed to estimate them. We will
return to the issue of constructing those solutions at the end of this section.
We shall now assume that $\sigma^{(n)}$ and $Z^{(n-1)}$ satisfy 
\eqref{eq:bos} and \eqref{eq:boZ} and  prove the estimates  
\eqref{eq:siges} and \eqref{eq:Zest} for $\sigma^{(n+1)}$ and 
$Z^{(n)}$.
First we estimate $\dot{\tilde\sigma}^{(n+1)}$ in terms of $Z^{(n)}$.
Observe that $\dot\sigma^{(n+1)}$ verifies the system of ODE's described in Lemma 
\ref{lem:siglem} with the function $G(Z^{(n)}(t),\sigma^{(n)}(t)$ and $\Omega^m_j(t,x;\sigma^{(n)}(t))$. 
Therefore, 
\begin{equation}
\label{eq:tildesigma}
|\dot{\tilde\sigma}^{(n+1)}(t)|\le  (1+t)^{-n}(\frac 14 \delta^2 +
C \|Z^{(n)}\|_{\Xl_s}^2 + C \|Z^{(n)}\|_{\Xl_s}^p)
\end{equation}
We now  
consider the $Z^{(n)}$ equation. First, by construction $Z^{(n)}$ 
is orthogonal to 
$\xi_{j}^{m}(t,x;\sigma^{(n)}(t))$.
Observe also that $Z^{(n)}$ satisfies the equation 
$$
i\pa_t Z^{(n)} + H(t,{\sigma^{(n)}}^\infty) Z^{(n)} = \Big ( H(t,{\sigma^{(n)}}^\infty)- 
H(\sigma^{(n)}(t))\Big ) Z^{(n)} + \dot\Sigma^{(n+1)} W(\sigma^{(n)}(t)) + G(Z^{(n-1)}, \sigma^{(n)}(t))
$$
Corollary \ref{cor:Hdiff} implies that 
$$
H(t,{\sigma^{(n)}}^\infty)- 
H(\sigma^{(n)}(t)) = V(t,x),
$$
where a smooth localized potential $V$ has the property that 
$\sup_{0\le |\gamma |\le s} \|\pa^\gamma V(t,\cdot)\|_{L^1\cap L^\infty}\les \delta^2 (1+t)^{-n+2}$. 
Moreover, the calculation leading to \eqref{eq:difort} of Lemma \ref{lem:orthZ} shows that  
$$
\| \calM_j(\sigma,t) \calG_{v_j, D_j}(t) \xi_j^m(t,\cdot;\sigma^{(n)}(t))- 
\xi_j^m(\cdot;{\sigma^{(n)}}^\infty)\|_{L^1\cap L^2}\les 
\delta^2 (1+t)^{-n+2}
$$
since $\sigma^{(n)}$ is an admissible path satisfying the bootstrap assumptions.
Therefore, by the results of Corollary~\ref{cor:charge2}, taking into account smallness
of the initial data $Z^{(n)}(0)=Z_0$,
\begin{align}
\|Z^{(n)}\|_{\Xl_{s}}&\les \sum_{k=0}^{s}\|\nabla^k Z_0\|_{L^1\cap L^2}  + 
\|G(Z^{(n-1)},\sigma^{(n)})\|_{\Yl_s} + 
\|\dot\Sigma^{(n+1)} W(\sigma^{(n)}(t))\|_{\Yl_s}\nn \\
&\les  \delta^{2} + \|Z^{(n-1)}\|^2_{\Xl_{s}} + \|Z^{(n-1)}\|_{\Xl_s}^p + 
 \|Z^{(n)}\|_{\Xl_s}^2 +  \|Z^{(n)}\|_{\Xl_s}^p, \label{eq:Zn}\\ &\les  
\delta^2 +  \|Z^{(n)}\|_{\Xl_s}^2 +  \|Z^{(n)}\|_{\Xl_s}^p\label{eq:Sigman},
\end{align}
where the inequality leading to \eqref{eq:Zn} follows from the estimates on the
nonlinear term $G(Z^{(n-1)}(t),\sigma^{(n)}(t))$ and $\dot\Sigma^{(n+1)} W(\sigma^{(n)}(t))$
obtained in Lemma \ref{lem:FZ}. The bound on $\dot\Sigma^{(n+1)} W(\sigma^{(n)}(t))$ 
also uses the inequality \eqref{eq:tildesigma}. To pass to \eqref{eq:Sigman} we used 
the assumption that $Z^{(n-1)}$ satisfies  \eqref{eq:Zest}. 
In the same way one obtains a local in time version of equation \eqref{eq:Sigman}.
\begin{equation}
\label{eq:localZ}
\|Z^{(n)}\|_{\Xl_{s}(T)}\les \delta^2 +  
\|Z^{(n)}\|_{\Xl_s(T)}^2 +  \|Z^{(n)}\|_{\Xl_s(T)}^p,
\end{equation}
which by continuity in $T$ implies the desired estimate $\|Z^{(n)}\|_{\Xl_s} \les \delta^2$.
 
Thus, the sequence $Z^{(n)}$ is uniformly bounded and small in the space $\Xl_s$ while 
$ (1+t)^n \dot{\tilde \sigma}^{(n)}$ is uniformly small point wise in time. 
Therefore, we can choose a convergent subsequence of the 
paths $\sigma^{(k)}(t)\to \sigma(t)$ and a weekly convergent 
in $H^{s}(\R^{n})$ subsequence $Z^{(k)}\to Z$. We multiply the 
equation \eqref{eq:equatA} by a smooth compactly supported function
$\zeta(x)$, integrate over the entire space and pass to the limit using 
that on any compact set $Z^{(n)}\to Z$ strongly in $H^{s'}$ for any 
$s'<s$. In particular, since $s>\frac n2 $, $Z^{(n)}\to Z$ point wise.
It will follow that $Z$ is a solution of the equation  
\begin{align}
&i \pa_{t} Z+ H(\sigma(t)) Z = \dot\Sigma 
W(\sigma (t)) + G(Z(t),\sigma(t))\label{eq:limitA},\\
& Z(0,x)  = Z_{0}(x)\nn
\end{align}
We also pass to the limit in the equation \eqref{eq:orthn}
to obtain 
\be
\Bigl \la \dot\Sigma W(\sigma (t)),\xi_j^m(t,\cdot;\sigma(t))\Bigr \ra = 
 \Bigl \la G(Z(t),\sigma(t)),\xi_j^m(t,\cdot;\sigma(t))\Bigr \ra  +
\Bigl \la \Omega^m_j(t,\cdot;\sigma(t)), Z\Bigr \ra.\label{eq:orthnZ}
\ee
Comparing equations \eqref{eq:limitA} and \eqref{eq:orthnZ} 
we conclude that 
$$
\la Z(t),\xi_{j}^{m}(t,\cdot;\sigma(t))\ra =0
$$
for all $j, m$. Therefore, the function 
$\psi = R + w_{1} + w_{2}$ solves the original NLS and by 
uniqueness, say in $L^{2}$, $\psi$ is our original solution.

To show existence of the solution $Z^{(n)}, \sigma^{(n+1)}$ of the 
linear system \eqref{eq:equatA}, \eqref{eq:orthn} we first 
construct the solution on a small time interval.
We note that the "system" \eqref{eq:orthn} for 
$\dot{\tilde\sigma}^{(n+1)}$ can be resolved algebraically due 
to the spatial separation of the paths $\sigma^{(n)}_{j}(t)$.
Therefore, for simplicity we can replace the system 
\eqref{eq:equatA}, \eqref{eq:orthn} by the following caricature:
\begin{align*}
&i\pa_{t} z + \frac 12 \Laplace z = V(t,x) z +  \omega(t) a(t,x) + g(t,x),\\
& \omega (t) = \la z, b(t,\cdot)\ra + f(t)
\end{align*}
Here $V, a, b, f$ are sufficiently smooth given functions and $g(t,\cdot)\in H^s$
uniformly in $t$.
We eliminate $\omega(t)$ and infer that 
$$
i\pa_{t} z + \frac 12 \Laplace z = V(t,x) z +  \la z, c(t,\cdot)\ra 
a(t,x) + F(t,x)
$$
with some new smooth functions $c, a$ and an $H^s$ function $F(t,\cdot)$.
Using the standard energy estimates we obtain that 
$$
\|z(t)\|_{H^{s}}\le \|z_{0}\|_{H^{s}} + C_{1} t \sup_{\tau\le 
t}\|z(\tau)\|_{H^{s}} + C_2 t \sup_{\tau\le 
t}\|F(\tau,\cdot)\|_{H^{s}},
$$
where the constants $C_{1}$ and $C_{2}$ depend on  $V, a, b, f$.
Therefore, we can establish the existence of the solution 
on the time interval of size $\frac 12 C_{1}^{-1}$ by means of the 
standard contraction argument. Then we can repeat this argument 
indefinitely thus constructing a global {\it {classical}} solution. 

\section{The linearized problem}
\label{ap:lin}

\subsection{Estimates for matrix charge transfer models}

In this section we recall some of the estimates from
Sections~7 and~8 from our companion paper~\cite{RSS}.
First, consider the case of a system with a single matrix
potential:
\be
\label{eq:sys_sing}
i\partial_t \binom{\psi_1}{\psi_2} + \bm H+U & -W \\ W & 
-H-U \endm \binom{\psi_1}{\psi_2}
=0
\ee
with $U,W$ real-valued and $H=\Lapl-\mu$, $\mu>0$.  
We say that $A:=\bm H+U 
& -W \\ W & -H-U \endm$ is {\em admissible} iff
the conditions of the following Definition~\ref{def:spec_ass}
hold. 

\begin{defi}
\label{def:spec_ass} Let $A$ be as above  with $U,W$ real-valued 
and exponentially decaying.
The operator $A$ on $\Dom(A)=H^2(\R^n)\times H^2(\R^n)\subset
\Hil:=L^2(\R^3)\times L^2(\R^3)$ is {\em admissible} provided
\begin{itemize}
\item $\spec(A)\subset \R$ and $\spec(A)\cap (-\mu,\mu) = 
\{\omega_\ell\:|\: 0\le\ell \le M\}$, for some $M<\infty$ 
      where $\omega_0=0$ and all $\omega_j$ are distinct eigenvalues. 
There are no eigenvalues in $\spec_{ess}(A)=(-\infty,-\mu]\cup[\mu,\infty)$.  
\item For $1\le \ell\le M$, 
$L_\ell:=\ker(A-\omega_\ell)^2=\ker(A-\omega_\ell)$, and 
$\ker(A)\subsetneq \ker(A^2)=\ker(A^3)=:L_0$. Moreover, these spaces 
are finite dimensional.
\item The ranges $\Ran(A-\omega_\ell)$ for $1\le\ell\le M$ and 
$\Ran(A^2)$ are closed.
\item The spaces $L_\ell$ are spanned by exponentially decreasing 
    functions in $\Hil$ (say with bound $e^{-\eps_0|x|}$).
\item The points $\pm \mu$ are not resonances of $A$.
\item All these assumptions hold as well for the adjoint $A^*$. We 
denote the corresponding (generalized)
eigenspaces by $L_\ell^*$.
\end{itemize}
\end{defi}

\noindent 
We will discuss these conditions in detail in the following
Subsection~\ref{subsec:specass}. It is possible to establish some of these properties
by means of ``abstract'' methods (for example, the exponential decay
of elements of generalized eigenspaces
via a variant of Agmon's argument, or the closedness of 
$\Ran(A-\omega_\ell)$ from Fredholm's theory), whereas others
can be reduced to statements concerning certain semi-linear 
elliptic operators $L_{+}, L_{-}$, see~\eqref{eq:L+L-}
(for example, that the spectrum is real or that only $0$  
can have a generalized eigenspace). In a later section we will 
prove for a particular model that $L_{+},L_{-}$ have the required properties.
\noindent One condition that we will not deal with
in this paper is the absence of embedded eigenvalues
in the essential spectrum. This property will remain an
assumption.

\smallskip

It is shown in \cite{RSS}, Lemma~7.2 that under these
conditions there is a direct sum decomposition
\be
\label{eq:split}
 \Hil = \sum_{j=0}^M L_j + \Bigl(\sum_{j=0}^M L_j^*\Bigr)^\perp 
\ee
and we denote by $P_s$ the induced projection 
onto~$\Bigl(\sum_{j=0}^M L_j^*\Bigr)^\perp$.
In general, $P_s$ is non-orthogonal. The letter ``s''
here stands for ``scattering'' (subspace). It is known
that $\Ran(P_s)$ plays the role of the scattering
states for the evolution~$e^{it A}$. Indeed, 
the main result from
Section~7 in~\cite{RSS} is that 
if~$A$ is admissible and the {\em linear stability property} 
\be
\label{eq:stab_cond}
\sup_t \|e^{it A} P_s \|_{2\to2} < \infty
\ee
holds,  then one has the dispersive bound
\be
\| e^{it A} P_s \psi_0 \|_{L^2+L^\infty} \lesssim
 |t|^{-\frac32} \| \psi_0 \|_{L^1 \cap L^2}
\label{eq:Adecay}
\ee
(if in addition $\|\hat{V}\|_1<\infty$, then the $L^2$ norm 
can be removed on the left-hand side). 
Next, we recall the notion of matrix charge transfer models
from Section~8 in~\cite{RSS}. 

\begin{defi}
\label{def:chargetransm} 
By a {\em matrix charge transfer model} we mean a system
\bea
&& i \partial_t \vpsi + \bm \Lapl & 0 \\ 0 & -\Lapl 
\endm\vpsi + \sum^\nu_{j =1}
V_j(\cdot - \vec{v_j} t) \vpsi = 0 \label{eq:transferm} \\
&& \vpsi |_{t=0} = \vpsi_0,  \nn
\eea
where $\vec v_j$ are distinct vectors in $\R^3$, and $V_j$ are matrix 
potentials of the form
\[ V_j(t,x) = \bm U_j(x) & -e^{i\theta_j(t,x)}\,W_j(x) \\ 
e^{-i\theta_j(t,x)}\,W_j(x) & -U_j(x) \endm, \]
where $\theta_j(t,x)=(|\vec v_j\,|^2+\alpha_j^2)t+2x\cdot \vec v_j + 
\gamma_j$, $\alpha_j,\gamma_j\in\R$,
$\alpha_j\not=0$. Furthermore, we require that
each
\[ H_j = \bm \Lapl - \half\alpha_j^2 + U_j & -W_j \\ W_j & -\Lapl + 
\half\alpha_j^2 - U_j \endm \]
be admissible in the sense of Definition~\ref{def:spec_ass} and  
that it satisfy the linear stability condition~\eqref{eq:stab_cond}. 
\end{defi}

\noindent It is clear that the Hamiltonian in~\eqref{eq:ref_hamil}
is of this form. As in Lemma~\ref{lem:conjug} above one
now verifies the following. The Galilean transforms $\calG_{\vec v}:=\calG_{\vec v,0}$
are defined as in~\eqref{eq:matrGal}, i.e.,  
\[ 
\calG_{\vec v}(t) \binom{f_1}{f_2} =
\binom{\calg_{\vec v,0}(t)f_1}{\overline{\calg_{\vec v,0}(t)\bar f_2}} 
\]
where $\calg_{\vec v,0}(t)= e^{-i\frac{|\vec v|^2}2 t} e^{-i{x\cdot \vec v}} e^{it\vec v\cdot p}$.

\begin{lemma}
\label{lem:trans_law} 
Let $\alpha\in\R$ and set
\[ A := \bm \Lapl - \half\alpha^2 + U & -W \\ W & -\Lapl +
\half\alpha^2 - U \endm \]
with real-valued $U,W$. Moreover, let 
$\vec v\in\R^3$,  $\theta(t,x)=(|\vec v\,|^2+\alpha^2)t+2x\cdot \vec 
v + \gamma$, $\gamma\in\R$, and define
\[ H(t) := \bm \Lapl+U(\cdot-\vec vt) & -e^{i\theta(t,\cdot-\vec 
vt)}W(\cdot-\vec vt) \\
 e^{-i\theta(t,\cdot-\vec vt)}W(\cdot-\vec vt) & -\Lapl-U(\cdot-\vec 
vt) \endm.
\]
Let $\calS(t)$, $\calS(0)=Id$, denote the propagator of the system
\[ i\partial_t \calS(t) + H(t)\calS(t) =0.\]
Finally, let 
\be
\label{eq:M_def}
\calM(t)=\calM_{\alpha,\gamma}(t)=\bm e^{-i\omega(t)/2} & 0 \\ 0 & 
e^{i\omega(t)/2} \endm 
\ee
where $\omega(t)=\alpha^2 t+\gamma$.  Then 
\be
\label{eq:trans_law}
 \calS(t) = \calG_{\vec v,0}(t)^{-1} \calM(t)^{-1} e^{itA} 
\calM(0)\calG_{\vec v,0}(0).
\ee
\end{lemma}
\begin{proof} One has
\begin{equation}
i \partial_t \calM(t)\calG_{\vec v}(t)\calS(t) = \bm 
\half\dot{\omega} & 0\\ 0 & - \half\dot{\omega} \endm 
\calM(t)\calG_{\vec v}(t)\calS(t) 
 + \calM(t) {i}\dot\calG_{\vec v}(t) \calS(t) - 
\calM(t)\calG_{\vec v}(t)H(t)\calS(t). \label{eq:big_dot}
\end{equation}
Let $\rho(t,x)=t|\vec v\,|^2+2x\cdot\vec v$. One now checks the 
following properties by differentiation:
\bea
\calM(t){i}\dot\calG_{\vec v}(t) &=& - \bm \half|\vec 
v\,|^2+\vec v\cdot\vec p & 0 \\
0 & -\half|\vec v\,|^2+\vec v\cdot\vec p \endm \calM(t)\calG_{\vec 
v}(t) \nn \\
\calM(t)\calG_{\vec v}(t)H(t) &=& \bm \Lapl+U & 
-e^{i(\theta-\rho-\omega)}W \\
 e^{-i(\theta-\rho-\omega)}W & - \Lapl-U \endm \calM(t)\calG_{\vec 
v}(t) \nn \\
&& \quad - \bm \half|\vec v\,|^2+\vec v\cdot\vec p & 0 \\
0 & -\half|\vec v\,|^2+\vec v\cdot\vec p \endm \calM(t)\calG_{\vec 
v}(t). \label{eq:diff_com}
\eea
The right-hand side of~\eqref{eq:diff_com} arises as follows.
First, the Galilei transform introduces a factor of $e^{-ix\cdot \vec v}$,
which needs to be commuted with $\Lapl$. Since
\bea
\Lapl \Big(e^{-ix\cdot \vec v} f\Big) &=& - \half|\vec v|^2 e^{-ix\cdot \vec v} f
- e^{-ix\cdot \vec v} i \vec v\cdot \vec \nabla f +
\half e^{-ix\cdot \vec v} \triangle f \nn \\
&=& \half|\vec v|^2 e^{-ix\cdot \vec v} f
-  i \vec v\cdot \vec \nabla \Bigl( f e^{-ix\cdot \vec v}\Bigr) +
\half e^{-ix\cdot \vec v} \triangle f 
\nn \\
&=& \Big(\half|\vec v|^2 
+  \vec v\cdot \vec p \Big) \Bigl( f e^{-ix\cdot \vec v}\Bigr) +
\half e^{-ix\cdot \vec v} \triangle f, \nn
\eea
one obtains the final term on the right-hand side of~\eqref{eq:diff_com}.
It remains to check the terms involving the potentials 
(for simplicity $\theta(\cdot-t\vec v)=\theta(t,\cdot-\vec v t))$:
\bea
&& 
\calM(t)\calG_{\vec v}(t) \bm U(\cdot-\vec vt) & -e^{i\theta(t,\cdot-\vec 
vt)}W(\cdot-\vec vt) \\
 e^{-i\theta(t,\cdot-\vec vt)}W(\cdot-\vec vt) & -U(\cdot-\vec 
vt) \endm \binom{f_1}{f_2} \nn \\
 &=&  
\bm e^{-i\omega(t)/2} & 0 \\ 0 & 
e^{i\omega(t)/2} \endm \binom{\calg_{\vec v}(t) U(\cdot-\vec v t)f_1 - 
\calg_{\vec v}(t) e^{i\theta(\cdot-\vec vt)} W(\cdot -\vec vt)f_2}
{\overline{\calg_{\vec v}(t) 
e^{i\theta(\cdot-\vec vt)} W(\cdot-\vec v t)\overline{f_1}} - 
\overline{\calg_{\vec v}(t)  U(\cdot -\vec vt)\overline{f_2}}} \nn \\
 &=&   
\binom{U\calg_{\vec v}(t) (e^{-i\omega(t)/2} f_1) 
- W e^{-i(v^2 t+2 x\cdot \vec v)} e^{i(\theta-\omega)} 
\overline{\calg_{\vec v}(t) \overline{e^{i\omega(t)/2} f_2}}}
{We^{i(v^2t+2x\cdot\vec v)} e^{i(\omega-\theta)} \;\calg_{\vec v}(t) (e^{-i\omega(t)/2} f_1) 
- U \;\overline{\calg_{\vec v}(t) \overline{e^{i\omega(t)/2} f_2}}} \nn \\
&=& \bm U & - e^{i(\theta-\omega-\rho)} W \\ e^{-i(\theta-\omega-\rho)} W & -U \endm
\bm e^{-i\omega(t)/2} & 0 \\ 0 & e^{i\omega(t)/2} \endm \binom{\calg_{\vec v}(t) f_1}
{\overline{\calg_{\vec v}(t)\overline{f_2}}}\nn,
\eea
as claimed.
In view of our definitions, $\theta-\rho-\omega=0$. Since 
$\dot{\omega}=\alpha^2$, the lemma
follows by inserting~\eqref{eq:diff_com} into~\eqref{eq:big_dot}.
\end{proof}

In order to prove our main dispersive estimates for such
matrix charge transfer problems we need to formulate a
condition which ensures  that the initial condition belongs 
to the stable subspace. 
To do so, let $P_s(H_j)$ and~$P_b(H_j)$ be the projectors 
induced by the decomposition~\eqref{eq:split} for the
operator~$H_j$. 
Abusing terminology somewhat, we refer to $\Ran(P_b(H_j))$ as the 
{\em bound states} of~$H_j$. 

\begin{defi} 
\label{def:asympm}  Let $U(t) \vpsi_0 = \vpsi(t, \cdot)$ be the
solution of~\eqref{eq:transferm}. We say that $\vpsi_0$ is a
{\em scattering state} relative to~$H_j$ if 
$$
\|P_b(H_j,t)U(t) \vpsi_0 \|_{L^2} \to 0\text{ as }t\to +\infty.
$$
Here
\be
\label{eq:Proj2m}
P_b(H_j,t) := \calG_{\vec v_j}(t)^{-1}\calM_j(t)^{-1} P_b(H_j)\, 
\calM_j(t)\calG_{\vec v_j}(t)
\ee
with $\calM_j(t)=\calM_{\alpha_j,\gamma_j}(t)$ as in \eqref{eq:M_def}.
\end{defi}

The formula~\eqref{eq:Proj2m} is of course motivated 
by~\eqref{eq:trans_law}. 
Clearly,  $P_b(H_j,t)$ is the projection onto the bound states of 
$H_j$ that have been 
translated to the
position of the matrix potential $V_j(\cdot-t\vec{v}_j)$. 
Equivalently, one 
can think of it as translating the solution of~\eqref{eq:transferm} 
from that position to the 
origin, projecting onto the bound states of~$H_j$, and then 
translating back. 

\noindent We now formulate our decay estimate for matrix charge 
transfer models, see Theorem~8.6 in~\cite{RSS}.

\begin{theorem}
\label{thm:mainm}
Consider the matrix charge transfer model as in 
Definition~\ref{def:chargetransm}.
Let $U(t)$ denote the propagator of the 
equation~\eqref{eq:transferm}. Then 
for any initial data $\vpsi_0 \in L^1\cap L^2$, which is 
a scattering state relative to each~$H_j$
in the sense of Definition~\ref{def:asympm}, 
one has the decay estimates
\be
\| U(t) \vpsi_0 \|_{L^\infty} \lesssim
\langle t\rangle^{-\frac32}\|\vpsi_0\|_{L^1\cap L^2}. \label{eq:mainm}
\ee
\end{theorem}

For technical reasons, we need the estimate~\eqref{eq:mainm}
for perturbed matrix charge transfer equations, as described
in the following corollary. This is discussed in Remark~8.6
in~\cite{RSS}. 

\begin{cor}
\label{cor:per_chtr}
Let $\vec \psi$ be a solution of the equation 
\bea
&& i \partial_t \vpsi + \bm \Lapl & 0 \\ 0 & -\Lapl 
\endm\vpsi + \sum^\nu_{j =1}
V_j(\cdot - \vec{v_j} t) \vpsi +  V_0(t,x)\vpsi = 0 \label{eq:transferm_per} \\
&& \vpsi |_{t=0} = \vpsi_0,  \nn
\eea
where everything is the same as in Definition~\ref{def:chargetransm} 
up to the perturbation $V_0(t,x)$ which satisfies
\[ \sup_{t} \|V_0(t,\cdot)\|_{L^1\cap L^\infty} < \eps.\]
Let $\tilde{U}(t)$ denote the propagator of the 
equation~\eqref{eq:transferm_per}. Then 
for any initial data $\vpsi_0 \in L^1\cap L^2$, which is 
a scattering state relative to each~$H_j$
in the sense of Definition~\ref{def:asympm} (with $U(t)$ replaced by~$\tilde U(t)$),  
one has the decay estimates
\be
\| \tilde U(t) \vpsi_0 \|_{L^\infty} \lesssim
\langle t\rangle^{-\frac32}\|\vpsi_0\|_{L^1\cap L^2}\label{eq:per_dec}
\ee
provided $\eps$ is sufficiently small.
\end{cor}
 
The corresponding inhomogeneous bound is stated in Section~\ref{ap:lin}. 

\subsection{The spectral properties I: general arguments}\label{subsec:specass}

In order for the linear estimates to apply, we need
to impose the conditions in Definition~\ref{def:spec_ass}
as well as the linear stability condition~\eqref{eq:stab_cond}
on the operators from~\eqref{eq:statHam}. 
The admissibility conditions of 
Definition~\ref{def:spec_ass} were motivated to a large
extent by Buslaev and Perelman~\cite{BP1}, who built on
earlier work of Weinstein~\cite{W1}. We now analyze these
conditions in detail. As before, 
\begin{equation}
\label{eq:HUV}
A :=  \bm H + U & - W\\ W & -H -U \endm = B + V 
\end{equation}
where $U,W$ are real-valued,  $H=\Lapl-\mu$ with $\mu>0$,
and $V$ is the matrix potential consisting of $U,W$. 
In this subsection we deal with those properties that 
can be dealt with by means of general arguments, that make no
use of any special structure of the operator. 

\begin{lemma}
\label{lem:easy_spec} 
Let the matrix potential $V$ be bounded and go to zero at infinity. Then
$(A-z)^{-1}$ is a meromorphic function 
in $\Omega:=\Compl\setminus (-\infty,-\mu]\cup [\mu,\infty)$.
The poles are eigenvalues of $A$ of finite
multiplicity and $\Ran(A-z)$ is closed for all $z\in\Omega$. 
Finally, the complement of $\Omega$ agrees with the essential spectrum
of~$A$, i.e., $\spec_{\rm ess}(A)  = (-\infty,-\mu]\cup [\mu,\infty)$.
\end{lemma}
\begin{proof}
Suppose that $z\in\Omega$. 
Then $B-z$ is invertible, and $A-z=\Big(1+V(B-z)^{-1}\Big)(B-z)$. 
Since $V(B-z)^{-1}$ is analytic and compact in that region of~$z$'s, 
the analytic Fredholm theorem implies that $1+V(B-z)^{-1}$
is invertible for all but a discrete set of $z$'s in 
$\Omega$. Furthermore, the poles are precisely eigenvalues of 
$A$ of finite multiplicity.
It is also a general property that the ranges $\Ran(1+V(B-z)^{-1})$ 
are closed. Indeed, if $K$ is any compact operator on a Banach space, then
it is well-known and also easy to see that $\Ran(I-K)$ is closed.
Since $B-z$ has a bounded inverse for all $z\in\Omega$, 
this implies that $\Ran(A-z)$ is closed, as claimed.
Conjugating by the matrix $P=\bm 1 & i\\1 & -i\endm$
leads to the Hamiltonians
\be
\tilde A := P^{-1}AP= i\left (\begin{array}{cc}
0 & H+V_1 \\
-H-V_2 & 0 \end{array}\right ) = \tilde B+V, \quad 
\tilde B = i\bm 0 & H \\ -H & 0 \endm,\quad V= i\bm 0 & V_1 \\ -V_2 & 0 
\endm
\label{eq:sys2} 
\ee 
where $V_1=U+W$ and $V_2=U-W$. 
The system~\eqref{eq:sys2} corresponds to writing a vector
in terms of real and imaginary parts, whereas~\eqref{eq:sys}
corresponds to working with the solution itself and its conjugate.
By means of the matrix $J=\bm 0 & i \\ -i & 0 \endm$ one can also 
write
\[ \tilde B = \bm H & 0 \\ 0 & H \endm J,\quad \tilde A = \bm H+V_1 & 0 \\ 0 & 
H+V_2 \endm J.\]
Since $\tilde B^*=\tilde B$ it follows that $\spec(\tilde B)\subset \R$. 
One checks that for $\Re z\ne0$
\bea
(\tilde B-z)^{-1} &=& (\tilde B+z)\bm (H^2-z^2)^{-1} & 0 \\ 0 & (H^2-z^2)^{-1} 
\endm  \nn \\
&=& \bm (H^2-z^2)^{-1} & 0 \\ 0 & (H^2-z^2)^{-1} \endm (\tilde B+z) 
\label{eq:Bspec} \\
(\tilde A-z)^{-1} &=& (\tilde B-z)^{-1} - (\tilde B-z)^{-1} W_1 \Bigl[1+W_2 J (\tilde B-z)^{-1} 
W_1\Bigr]^{-1} W_2 J (\tilde B-z)^{-1}
\label{eq:grill}
\eea
where $W_1$ and $W_2$ are the following matrix potentials that go to zero at infinity: 
\[
 W_1 = \bm |V_1|^{\half} & 0 \\ 0 & |V_2|^{\half} \endm,\quad W_2 = 
\bm |V_1|^{\half}\sign(V_1) & 0 \\ 0 & |V_2|^{\half}\sign(V_2) \endm.
\]
The inverse of the operator in brackets exists if $z=it$ with $t$ large, for example. 
Moreover, by the assumed decay of the potential the entire operator
that is being subtracted from the right-hand side is compact in that case. 
One is therefore in a position to apply Weyl's criterion, 
see Theorem~XIII.14 in~\cite{RS4},  whence
\be
\label{eq:ess_spec}
 \spec_{\rm ess}(A)=\spec_{\rm ess}(\tilde A) = (-\infty,-\mu]\cup [\mu,\infty). 
\ee
The identity~\eqref{eq:grill} goes back to Grillakis~\cite{Grill}.
\end{proof}

Next, we need  to locate possible eigenvalues of~$A$ 
or equivalently, $\tilde A$.  This will not be done on 
the same general level, but require 
analysis of $L_{+}, L_{-}$ from~\eqref{eq:L+L-}. But 
we first discuss another general property of the matrix
operator~$A$.

\begin{lemma}
\label{lem:agmon}
Let $A$ be as in \eqref{eq:HUV} with $U,W$
continuous and $W$ exponentially decaying, whereas $U$
is only required to tend to zero.  
If $f\in\ker(A-E)^k$ for some $-\mu<E<\mu$ and
some positive integer~$k$, then $f$ decays exponentially.
\end{lemma}
\begin{proof} We want to emphasize
that the following result is ``abstract'' and does not rely on 
any special structure of the matrix potential or on any properties of $L_{+}$
or $L_{-}$. We will use a variant of Agmon's argument~\cite{Ag}.
More precisely, 
suppose that for some $-\mu < E< \mu$, there are 
$\psi_1,\psi_2\in H^2(\R^n)$ so that 
\bea
(\Laplace - \mu + U)\psi_1 - W\psi_2 &=& E\psi_1 \nn \\
W\psi_1 + (-\Laplace + \mu - U)\psi_2 &=& E\psi_2. \label{eq:eig_sys}
\eea
As usual, $U,W$ are real-valued and exponentially
decaying, $\mu>0$. Suppose $|W(x)|\les e^{-b|x|}$. 
Then define the Agmon metrics
\bea 
\rho_E^{\pm}(x) &=& \inf_{\gamma:0\to x} L^\pm_{\rm Ag}(\gamma) \nn \\
L^\pm_{\rm Ag}(\gamma) &=& \int_0^1 \min\Big( \sqrt{(\mu\pm E-U(\gamma(t)))_{+}}\,,\,
b/2 \Big) \|\dot \gamma(t)\|\,dt \label{eq:Ag_def}
\eea
where $\gamma(t)$ is a $C^1$-curve with $t\in[0,1]$,
and the infimum is to be taken over such curves that 
connect $0,x$.  These functions satisfy
\be
\label{eq:ag_met}
 |\nabla \rho_E^{\pm}(x)|\le \sqrt{(\mu\pm E-U(x))_{+}}.
\ee
Moreover, one has $\rho_E^{\pm}(x)\le b|x|/2$ by construction.
Now fix some small $\eps>0$ and 
set $\omega^\pm(x):= e^{2(1-\eps)\rho^{\pm}_E(x)}$. 
Our goal is to show that
\be
\label{eq:Agmon}
\int \Big[ \omega^+(x) |\psi_1(x)|^2 + \omega^{-}(x)|\psi_2(x)|^2\Big] \, dx < \infty.
\ee
Not only does this exponential decay in the mean suffice for our applications 
(cf.~Section~7 in~\cite{RSS}), but it can also be improved
to point wise decay using regularity estimates for $\psi_1,\psi_2$. 
We do not elaborate on this, see for example \cite{Ag} and 
Hislop, Sigal~\cite{HS}.

\noindent Fix $R$ arbitrary and large. For technical reasons, we set
\[ \rho_{E,R}^{\pm}(x):= \min\Big(2(1-\eps)\rho^{\pm}_E(x),R\Big),\quad  
\omega_R^\pm(x):= e^{\rho_{E,R}^{\pm}(x)}. \]
Notice that \eqref{eq:ag_met} remains valid in this case, and also that 
$\rho_E^{\pm}(x)\le \min(b|x|/2,R)$. Furthermore, 
by choice of $E$ there is a smooth functions $\phi$ that is 
equal to one for large $x$ so that
\[ 
\supp(\phi)  \subset \{\mu+E-U>0\} \cap \{\mu-E-U>0\}.
\]
It will therefore suffice to prove the following modified form
of~\eqref{eq:Agmon}:
\be
\label{eq:Agmon'}
\sup_{R}\int \Big[ \omega_R^+(x) |\psi_1(x)|^2 + 
\omega_R^{-}(x)|\psi_2(x)|^2 \Big]\phi^2(x)\, dx < \infty.
\ee
All constants in the following argument will be independent of~$R$. 
By construction, there is $\delta>0$ such that
\bea
&& \delta \int  \omega_R^+(x) |\psi_1(x)|^2 \phi^2(x)\,dx \le 
\int \omega_R^+(x) ( \mu+E-U(x)) |\psi_1(x)|^2 \phi^2(x)\,dx  \label{eq:delta_int} \\
&& = 
\int \omega_R^+(x) (\Laplace \psi_1 - W\psi_2)(x) \bar\psi_1(x) \phi^2(x)\,dx  \nn \\
&& = -  \int \nabla(\omega_R^{+}(x) \phi^2(x)) \nabla \psi_1(x) \bar\psi_1(x)\,dx 
 -  \int \omega_R^{+}(x) \phi^2(x) |\nabla \psi_1(x)|^2\,dx  \label{eq:nabla_rho} \\
&&  -  \int \omega_R^+(x) W(x) \psi_1(x)\bar\psi_2(x)\phi^2(x) \, dx. \label{eq:cross}
\eea
As far as the final term \eqref{eq:cross} is concerned, notice that
$\sup_{x,R} |\omega_R^+(x) \phi^2(x) W(x)|\les 1$ by construction, whence
$ |\eqref{eq:cross}| \les \|\psi_1\|_2\|\psi_2\|_2$. 
Furthermore, by \eqref{eq:ag_met} and Cauchy-Schwarz, 
the first integral in~\eqref{eq:nabla_rho} satisfies 
\bea 
&& \left| \int \nabla(\omega_R^{+}(x) \phi^2(x)) \nabla \psi_1(x) \bar\psi_1(x)\,dx 
\right| \nn \\
&& \le 2(1-\eps) \Bigl( \int \omega_R^+(x)(\mu+E-U(x)) \phi^2(x)|\psi_1(x)|^2 \, dx \Bigr)^{\half}
\Bigl( \int \omega_R^+(x) \phi(x)^2\, |\nabla\psi_1(x)|^2 \, dx \Bigr)^{\half} \label{eq:int_again} \\
&& \qquad + 2 \Bigl( \int \omega_R^+(x) \phi^2(x)|\nabla \psi_1(x)|^2 \, dx \Bigr)^{\half}
\Bigl( \int \omega_R^+(x) |\nabla\phi(x)|^2\, |\psi_1(x)|^2 \, dx \Bigr)^{\half} \nn 
\eea
Since the first integral in~\eqref{eq:int_again} 
is the same as that in~\eqref{eq:delta_int}, inserting~\eqref{eq:int_again}
into~\eqref{eq:nabla_rho} yields after some simple manipulations
\bea
 \eps \int \omega_R^+(x) ( \mu+E-U(x)) |\psi_1(x)|^2 \phi^2(x)\,dx 
&\le& \eps^{-1} \int \omega_R^+(x) |\nabla\phi(x)|^2\, |\psi_1(x)|^2\, dx \nn \\
&&  - \int \omega_R^+(x) \phi(x)^2 W(x) \psi_2(x)\,\bar\psi_1(x)\, dx. \nn
\eea
Since $\nabla\phi$ has compact support, and by our previous considerations 
involving $\omega_R^+ W$, the entire right-hand side is bounded independently of~$R$, 
and thus also~\eqref{eq:delta_int}. A symmetric argument applies to
the integral with $\psi_2$, and~\eqref{eq:Agmon'}, \eqref{eq:Agmon} hold.
This method also shows that functions belonging to generalized eigenspaces
 decay exponentially. Indeed, suppose $(A-E)\vec g=0$ and 
$(A-E)\vec f=\vec g$. Then 
\bea
(\Laplace - \mu + U)f_1 - W f_2 &=& E f_1 + g_1\nn \\
W f_1 + (-\Laplace + \mu - U) f_2 &=& E f_2 + g_2 \nn
\eea
with $g_1,g_2$ exponentially decaying. 
Decreasing the value of $b$ in~\eqref{eq:Ag_def} 
if necessary allows one to use the same argument as before 
to prove~\eqref{eq:Agmon} for~$\vec f$. By induction, one then
deals with all values of~$k$ as in the statement of the lemma.
\end{proof}

\subsection{The spectral properties II: reduction to $L_+,L_{-}$ }\label{subsec:specass2}

We now need to specialize $A$ from~\eqref{eq:HUV}
to the form~\eqref{eq:statHam}, i.e., 
\begin{equation}
A= \left (\begin{array}{cc}
\half\Laplace -\frac{\alpha^2}2 + \beta(\phi^2) 
+\beta'(\phi^2)\phi^2 &  
\beta'(\phi^2)\phi^2  \\
-\beta'(\phi^2)\phi^2 & -\half\Laplace 
+ \frac{\alpha^2}2 -\beta(\phi^2) -
\beta'(\phi^2)\phi^2 
	\end{array}\right ). \label{eq:sys}  
\end{equation}
As shown in Section~\ref{sec:results}, these are the stationary Hamiltonians derived from
the linearization of NLS, see~\eqref{eq:statHam} and Lemma~\ref{lem:trans_law}. 
Let $\alpha>0$ and~$\phi$ be a nonzero solution of 
\be
\label{eq:ground} \Lapl \phi -\frac{\alpha^2}{2}\phi +\beta(\phi^2)\phi=0,
\ee
i.e., $\phi=\phi(\cdot;\alpha)$ for all $\alpha\in(\alpha_0-c_0,\alpha_0+c_0)$.
Moreover, $\phi$ is smooth in both variables and 
\begin{equation}
\label{eq:smooth_phi}
\|\partial_\alpha\phi\|_{H^1(\R^n)}+\|\partial^2_\alpha\phi\|_{H^1(\R^n)} <\infty.
\end{equation}
Finally, we require that $\phi$ is sufficiently rapidly decaying.
A particular case would be a $\phi$ which is positive and radially symmetric. 
Such a solution is known to exist and to be unique if $\beta(u)=|u|^\sigma$ 
provided $0<\sigma<\frac{2}{d-2}$ and is referred to as the ``ground state''. 
It decays exponentially. 
However, in order to keep this section as general as possible, we do not
require~$\phi$ to be the ground state. 
Let
\be
\label{eq:L+L-}
 L_- := - \Lapl  + \frac{\alpha^2}{2} - \beta(\phi^2), \qquad 
L_+ := - \Lapl + \frac{\alpha^2}2 - \beta(\phi^2) -
2\beta'(\phi^2)\phi^2
\ee
with domains $\Dom(L_+)=\Dom(L_{-})=H^2(\R^n)$ so that
\be
\label{eq:tildeAdef}
 \tilde A = \bm 0 & -i L_{-} \\ i L_{+} & 0 \endm 
\ee
with $\Dom(\tilde A) = H^2(\R^n)\times H^2(\R^n)$.
Here $\tilde A$ is obtained by conjugating $A$ with the 
matrix~$P$, see~\eqref{eq:sys2}.
For simplicity, however, 
we no longer distinguish between $A$ and~$\tilde A$, i.e., we set
$A=\tilde A$. The spectrum of $L_{\pm}$ on $[\mu,\infty)$
is purely absolutely continuous, and below~$\mu=\frac{\alpha^2}{2}>0$ there 
are at most a finite number of eigenvalues of finite multiplicity 
(by Birman-Schwinger, the assumed decay of $\phi$ as well as $\beta(0)=0$).
Clearly,
\be
\label{eq:bas_rel}
L_-\phi=0,\quad L_+(\partial_j \phi)=0,\;1\le j\le n,\quad
L_+(\partial_\alpha \phi)=-\alpha\phi,
\ee
where the final property is formal. We now collect some crucial 
properties discovered by M.~Weinstein.

\medskip

\begin{defi}\label{def:L+L-}
\noindent 
\fbox{{\bf Needed properties of the scalar elliptic operators $L_{+}$ and $L_{-}\;$}:}

\smallskip\noindent Let $\phi(\cdot;\alpha)$ be as above,
in particular assume that~\eqref{eq:smooth_phi} holds.
 The kernels have the following explicit form:
\[\ker(L_{-})=\spa\{\phi\} \text{\ \ and\ \ }
\ker(L_{+})=\spa\{\partial_j\phi\:|\: 1\le j\le n\}.
\] 
 The operator $L_{+}$ has a single negative eigenvalue~$E_1$ with a
 unique ground state $\psi>0$, whereas $L_{-}$ is
 nonnegative, with $0$ as an isolated eigenvalue. 
Furthermore, the {\em convexity condition}  $\la
 \partial_\alpha \phi(\cdot;\alpha),
\phi(\cdot;\alpha)\ra > 0$ holds,
see~\eqref{eq:stabil}.
\end{defi}

\medskip
In the following section, these conditions will be verified for
a particular choice of nonlinearity. 
Note that according to this definition, $\phi$ is the ground state
of the linear operator $L_{-}$ and as such is positive.   
In case of $\phi(\cdot;\alpha)$ being the ground state of the nonlinear problem~\eqref{eq:ground}, 
these properties have been shown to hold by Weinstein~\cite{W1} 
and~\cite{W2} in case of power nonlinearities, i.e., $\beta(u)=|u|^\sigma$, 
$0<\sigma<\frac{2}{d-2}$. 

The purpose of this subsection is to reduce some of the admissibility 
conditions from Definition~\ref{def:spec_ass} to the properties
 of $L_{+}, L_{-}$ from Definition~\ref{def:L+L-}. Recall from the previous
subsection that several other properties hold in greater generality.
We now collect those properties that follow from Definition~\ref{def:L+L-}
into a single proposition.

\begin{proposition}
\label{prop:koschmar}
Impose the spectral assumption on $L_{+}$ and $L_{-}$ 
from Definition~\ref{def:L+L-}. Then
\begin{itemize}
\item  $\spec(A)\subset\R$, the only eigenvalue
that admits a generalized eigenspace is~$0$, and $\Ran(A^2)$ is closed. 
\item the linear stability condition holds.
\item equality holds in the relation concerning the root spaces~\eqref{eq:rootA} 
and~\eqref{eq:rootA*}.
In particular, one has $\ker(A^2)=\ker(A^3)$ and $\ker((A^*)^2)=\ker((A^*)^3)$
 and the rootspace ${\cal N}(A^*)$ has dimension $2n+2$.
\end{itemize}
\end{proposition}

The proof of this proposition is split into several lemmas below.

\begin{lemma}
\label{lem:fine_spec}
Impose the spectral assumption on $L_{+}$ and $L_{-}$ 
from Definition~\ref{def:L+L-}. Then $\spec(A)\subset\R$, the only eigenvalue
that admits a generalized eigenspace is~$0$, and $\Ran(A^2)$ is closed. 
\end{lemma}
\begin{proof}
Consider 
\be
\label{eq:TTstar}
A ^2 = \bm T^* & 0 \\ 0 & T \endm, \quad T = L_{+}L_{-} 
\ee
with domain $H^4(\R^n)=W^{4,2}(\R^n)$. 
Following~\cite{BP1}, we
first show that any eigenvalue of~$T$, and therefore also of 
$\spec(A ^2)$ is real,
and then under the assumption~\eqref{eq:stabil}, that it is
nonnegative. Because of~\eqref{eq:ess_spec}, the latter then 
implies that $\spec(A )$
is real, as required in Definition~\ref{def:spec_ass}. 
Clearly, $T\phi=0$. Let $\psi\not\in\spa\{\phi\}$, $T\psi=E\psi$. 
Let $\psi=\psi_1+c\phi$, $\psi_1\perp\phi$. Then 
\[ L_{-}^{\half} L_{+} L_{-}^{\half} L_{-}^\half \psi_1 = EL_{-}^\half\psi_1,\]
so that $L_{-}^\half\psi_1\ne0$ is an eigenfunction of the
symmetric operator $ L_{-}^{\half} L_{+} L_{-}^{\half}$ (with domain $H^4(\R^n)$),
and thus $E$ is real. Hence any eigenvalue of $A$ can only be real
or purely imaginary. 
Since $\phi\perp \ker(L_{+})$ by our assumption concerning $L_{+}$, 
the function
\[  g(E) := \la (L_{+}-E)^{-1}\phi,\phi \ra \]
is well-defined on an interval of the form $(E_1,E_2)$ for some $E_2>0$.
Moreover, 
\[ g'(E)=\|(L_{+}-E)^{-1}\phi\|^2>0\]
so that $g(E)$ is strictly increasing on the interval. Finally, 
\be
\label{eq:gnull} g(0) = -\frac{1}{\alpha} \la \partial_\alpha \phi,\phi \ra < 0
\ee
in view of \eqref{eq:bas_rel} and \eqref{eq:stabil}. Now suppose that $A^2$ has
a negative eigenvalue. Then by the preceding, so does~$T$,
and therefore also $ L_{-}^{\half} L_{+} L_{-}^{\half}$. 
More precisely, the argument from before implies that 
there is $\chi\in \ker(L_{-})^\perp$, $\chi\ne0$,  so that 
\[ \la L_{-}^\half L_{+} L_{-}^\half \chi,\chi \ra = 
\la L_{+}\psi,\psi \ra < 0\]
with $\psi=L_{-}^\half \chi$. Let $P_{-}^\perp$ denote the projection
onto the orthogonal complement of~$\ker(L_{-})=\spa(\phi)$. 
By the Rayleigh principle this implies that the self-adjoint operator 
$P_{-}^\perp L_{+}P_{-}^\perp$ has a negative eigenvalue, say $E_3<0$.
Thus $L_{+}\psi=E_3\psi+c\phi$ for some $\psi\perp\phi$. 
If $c=0$, then $E_3=E_1$ so 
that $\psi>0$ as the ground state of~$L_{+}$. But then $\la \phi,\psi \ra>0$, which is
impossible. So $c\ne0$, and one therefore obtains
\[ (L_{+}-E_3)^{-1}\phi=\frac{1}{c}\psi \quad\Longrightarrow\quad g(E_3)=0.\]
But this contradicts \eqref{eq:gnull} by strict monotonicity of $g$.
Thus $A^2$ does not have any negative eigenvalues, which implies that
$A$ does not have imaginary eigenvalues. 
Hence all eigenvalues of $A$ are real, as desired. 

We now turn to generalized eigenspaces.
Suppose $A\psi=E\psi+\chi$, where $E\ne0$, $(A-E)\chi=0$ and $\chi\ne0$.
This is equivalent to saying that $A$ has a generalized
eigenspace at~$E$. Then $\psi,\chi\in\Dom(A^2)$, and moreover 
\[ (A^2-E^2)\chi=0,\quad (A^2-E^2)\psi=(A-E)\chi+2E\chi = 2E\chi,\]
so that $A^2$ would have a generalized eigenspace at~$E$, and therefore
also~$T$. 
Hence, suppose $T\psi=E\psi$, with $E\ne0$, $\psi\ne0$.
If $(T-E)\chi=c\psi$ with $c\ne0$, then
\[ (L_{-}^{\half} L_{+} L_{-}^{\half} -E)L_{-}^{\half}
\chi_1=cL_{-}^{\half}\psi_1\ne0,
\qquad (L_{-}^{\half} L_{+} L_{-}^{\half} -E)^2L_{-}^{\half}
\chi_1=cL_{-}^{\half}(L_{+}L_{-}-E)\psi=0\]
where $\psi_1,\chi_1$ denote the projections of $\psi,\chi$ onto
the orthogonal complement of $\phi$. But  $L_{-}^{\half}\psi_1\ne0$ since
$E\ne0$ and thus $E$ would have to be a generalized eigenvalue of 
$L_{-}^{\half} L_{+} L_{-}^{\half}$, which is impossible.
So only $E=0$ can have a generalized eigenspace. Here we used the 
property that $L_{-}^{\half} L_{+} L_{-}^{\half}$ is self-adjoint
on its domain~$H^4(\R^n)$. While symmetry is obvious, self-adjointness on~$H^4(\R^n)$
requires a bit more care. 
Suppose $\la L_{-}^{\half} L_{+} L_{-}^{\half}f, g \ra = \la f, h\ra$
for all $f\in H^4(\R^n)$, and some fixed $g,h\in L^2(\R^n)$. Taking $f\in\ker(L_{-})$
shows that $P_{-}^\perp h=h$, i.e., that $h\in (\ker (L_{-}^\half))^\perp$. 
By the Fredholm alternative applied to the self-adjoint operator $L_{-}^\half$, 
one can write $h=L_{-}^\half h_1$ with some $h_1\in \Dom(L_{-}^\half)=H^1(\R^n)$. 
Note that $h_1$ is defined only up to an element in~$\ker(L_{-}^\half)$,
i.e., $h_1+c\phi$ has the same property for any constant $c$. Thus
\[  \la L_{-}^{\half} L_{+} L_{-}^{\half} f,g \ra = \la f, L_{-}^\half (h_1+c\phi)\ra
= \la L_{-}^\half f, h_1 + c\phi \ra\]
for all $f\in H^4(\R^n)$. Equivalently, setting $f_1=L_{-}^{\half} f$, one has 
\[ \la L_{-}^{\half} L_{+} f_1 ,g \ra = \la  f_1, h_1 + c\phi \ra. \]
Note that the class of $f_1$ are all functions in $H^3(\R^n)$ with $f_1\perp \phi$.
We now want to remove the latter restriction, which can be achieved by a suitable 
choice of~$c$. Indeed, in order to achieve
\[ \la L_{-}^{\half} L_{+} (f_1+\lambda \phi) ,g \ra = \la  f_1+\lambda\phi, h_1 + c\phi \ra \]
for all $f_1\in H^3(\R^n)$, $f_1\perp \phi$, $\lambda\in\C$
one chooses $c$ such that
\[ \la L_{-}^{\half} L_{+}  \phi ,g \ra = \la \phi, h_1 + c\phi \ra, \]
which can be done since $\la \phi,\phi \ra >0$. 
Renaming $h_1+c\phi$ into $h_1$, one thus arrives at 
\begin{equation}
\label{eq:zw}
 \la L_{-}^{\half} L_{+} f_1 ,g \ra = \la  f_1, h_1  \ra \text{\ \ for all\ \ }f_1\in H^3(\R^n). 
\end{equation}
Recall that $h=L_{-}^\half h_1$. One can now continue this procedure. Indeed,
since~\eqref{eq:zw} implies that $h_1\perp\ker(L_{+})$, one can write
$h_1=L_{+} h_2= L_{+}(h_2+\sum_{j=1}^k c_j\psi_j)$, where $\ker(L_{+})=\spa\{\psi_j\}_{j=1}^k$
and $h_2\in H^3(\R^n)$ (in fact, $\psi_j=\partial_j\phi$ by our assumption). 
As before, the constants $\{c_j\}$ are chosen in such a way that 
\[ \Big\la L_{-}^{\half} (L_{+} f_1 + \sum_{j=1}^k \lambda_j\psi_j),g \Big\ra = 
\Big\la L_{+} f_1 + \sum_{j=1}^k \lambda_j\psi_j, h_2+\sum_{j=1}^k c_j\psi_j \Big\ra\]
for all $\lambda_j$. 
This can be done because of the invertibility of the Gram matrix of $\{\psi_j\}_{j=1}^k$.
Hence  
\[
 \la L_{-}^{\half} f_2 ,g \ra = \la  f_2, h_2  \ra \text{\ \ for all\ \ }f_2\in H^1(\R^n). 
\]
Moreover, $h=L_{-}^\half L_{+}h_2$ with $h_2\in H^3(\R^n)$. 
By the self-adjointness of $L_{-}^\half$ this implies that $h_2=L_{-}^\half g$.
It follows that $g\in H^4(\R^n)$ and $h= L_{-}^\half L_{+} L_{-}^\half g$ as desired. 

Finally, we show that $\Ran(A^2)$ is closed. By~\eqref{eq:TTstar} it suffices to 
show that the ranges of both $T=L_{+}L_{-}$ and $T^*=L_{-}L_{+}$ 
are closed with domain $H^4(\R^n)$. We will first verify that these
operators are closed on this domain. Indeed, they each can be written in the form
$\triangle^2+F_1\triangle+\triangle F_2+F_3$. Since for $M$ large
\[
\| (\triangle^2+M+F_1\triangle+\triangle F_2+F_3)f\|_2 \ge 
\| (\triangle^2+M)f\|_2 - C\|f\|_{W^{2,2}} 
\ge \half \| (\triangle^2+M)f\|_2\gtrsim \|f\|_{W^{4,2}}, 
\]
one concludes that $T+M$, $T^*+M$ are closed, and therefore
also $T,T^*$. Next, let $P_{-}$ and $P_{+}$ be the
projections onto $\ker(L_{-})$ and $\ker(L_{+})$, respectively.
Then $L_{-}L_{+}=L_{-}P_{-}^\perp L_{+}P_{+}^\perp$. Now
\be
\label{eq:switch}
 P_{-}^\perp L_{+}f = L_{+}f - \|\phi\|^{-2}\la L_{+}f,\phi\ra \phi = 
L_{+}\Big(f-\|\phi\|^{-2}\la f,L_{+}\phi\ra \tilde\phi\Big),
\ee
where we have written $\phi=L_{+}\tilde\phi$, 
$\tilde\phi\in \ker(L_{+})^\perp= \Ran(P_{+}^\perp)$ by virtue of
the fact that 
\[\phi\in\Ran(L_{+})=\overline{\Ran(L_{+})}=\ker(L_{+})^\perp 
= \spa\{\partial_j\phi\}^\perp.\]
The last equality here is Weinstein's characterization, more precisely,
our assumption on~$L_{+}$. Define 
\[ \tilde \phi := P_{+}^\perp\tilde\phi_1,\quad 
Qf:= f-\|\phi\|^{-2}\la f,L_{+}\phi\ra \tilde\phi,
\quad \tilde Qf:= f-\|\phi\|^{-2}\la f,P_{+}^\perp L_{+}\phi\ra \tilde\phi_1 \]
so that \eqref{eq:switch} gives 
\[ P_{-}^\perp L_{+}P_{+}^\perp = L_{+}QP_{+}^\perp=L_{+}P_{+}^\perp\tilde Q
\quad \Longrightarrow\quad L_{-}L_{+} = L_{-}L_{+}P_{+}^\perp\tilde Q. \]
In particular, 
\bea
\|T^* f\|_2 &=& \|L_{-}P_{-}^\perp L_{+}P_{+}^\perp f\|_2 \ge 
c_1\, \|P_{-}^\perp L_{+}P_{+}^\perp f\|_2 \nn \\
&=& c_1\, \|L_{+} QP_{+}^\perp f\|_2 = c_1\, \|L_{+} P_{+}^\perp \tilde Q f\|_2
\ge c_1c_2\, \|P_{+}^\perp \tilde Q f\|_2, \label{eq:c1c2}
\eea
where the existence of $c_1,c_2>0$ follows from the 
self-adjointness of $L_{-},L_{+}$. Hence, if $T^* f_n=T^* P_{+}^\perp \tilde Q f_n
 \to h$ in $L^2$, then 
by~\eqref{eq:c1c2} and linearity  $P_{+}^\perp \tilde Q f_n \to g$ in~$L^2$. 
Since $T^*$ was shown to be closed, it follows that $h=T^* g$, and $\Ran(T^*)$
is closed. A similar argument shows that $\Ran(T)$ is closed. 
\end{proof}

Next, we derive the {\em linear stability assumption} as well as the
structure of the generalized eigenspaces of $A$ and $A^*$
from our spectral assumptions on $L_{-},L_{+}$. 
 From the spectral assumptions in Definition~\ref{def:L+L-} 
as well as 
\[ L_{-}\phi=0,\; L_{-}(x_j\phi)=-\partial_j\phi,\;
 L_{+}(\partial_j\phi) = 0,\; L_{+}(\partial_\alpha\phi)=-\alpha\phi
\]
it follows that
\bea
\ker(A) &=& \spa\Big\{ \binom{0}{\phi},\;\binom{\partial_j\phi}{0}\;: 1\le j\le n \Big\} \nn \\
\ker(A^*) &=& \spa\Big\{ \binom{\phi}{0},\;\binom{0}{\partial_j\phi}\;:1\le j\le n \Big\}\nn \\ 
\calN(A) &:=& \bigcup_{k=1}^\infty \ker(A^k) \supset \spa\Big\{ \binom{0}{\phi},\;
\binom{\partial_\alpha\phi}{0},\;\binom{0}{x_j\phi},\;\binom{\partial_j\phi}{0}\;: 1\le j\le n \Big\}
=:\calM \label{eq:rootA} \\
\qquad\calN(A^*) &:=&\bigcup_{k=1}^\infty \ker((A^*)^k) \supset \spa\Big\{ \binom{\phi}{0},\;
\binom{0}{\partial_\alpha\phi},\;\binom{x_j\phi}{0},\;\binom{0}{\partial_j\phi}\;:1\le j\le n \Big\}
=:\calM_*. \label{eq:rootA*} 
\eea
One of our goals is to show that equality holds in the last two relations.
This is the same as the structure statement made in Proposition~\ref{prop:null},
but one needs to apply the matrix~$P=\bm 1 & i\\1 & -i\endm$ to pass between 
these two representations. 

Now suppose that $i\partial_t\vec \psi+A\vec \psi=0$. This can be written as
$\partial_t \vec\psi+ JM\vec\psi=0$ where $J=\bm 0& -1\\ 1 & 0 \endm$ and 
$M=\bm L_{+} & 0 \\ 0 & L_{-}\endm$. Therefore,
\[
\frac{d}{dt} \la \vec\psi, M \vec\psi \ra = 2\Re \la \partial_t\vec\psi , M\vec\psi \ra 
= -2\Re \la JM\vec\psi , M\vec\psi \ra =0 
\]
by anti-selfadjointness of $J$. In other words,
\[ Q(\vec\psi) := \la L_{+}\psi_1, \psi_1\ra + \la L_{-}\psi_2,\psi_2 \ra \]
is constant in time if $\vec \psi(t) = e^{itA}\vec\psi(0)$ 
(here $\vec\psi=\binom{\psi_1}{\psi_2}$). Although the previous calculation 
basically required classical solutions, it is clear that its
natural setting are $H^1(\R^n)$-solutions. In that case one needs to interpret
the form $Q(\vec \psi)$ via
\bea
\la L_{+}\psi_1, \psi_1\ra &=& \half\|\nabla \psi_1\|_2^2 + \frac{\alpha^2}{2}\|\psi_1\|_2^2 
-\la \beta(\phi^2)\psi_1,\psi_1\ra \nn \\
\la L_{-}\psi_2, \psi_2\ra &=& \half\|\nabla \psi_2\|_2^2 + \frac{\alpha^2}{2}\|\psi_2\|_2^2 
-\la (\beta(\phi^2)+2\beta'(\phi^2)\phi^2)\psi_2,\psi_2\ra. \label{eq:formH1}
\eea
In what follows, we will tacitly make this interpretation whenever it is needed. 
The following lemmas are due to Weinstein~\cite{W1}.

\begin{lemma}
\label{lem:E1}
Impose the spectral assumptions on $L_{+},L_{-}$ from Definition~\ref{def:L+L-}.
Then $\la L_{+}f,f\ra \ge0$ for all $f\in H^1(\R^n)$, $f\perp \phi$.
\end{lemma}

This is a special case of Lemma~E.1 in~\cite{W1}, and we refer the reader
to that paper for the proof. 

\begin{lemma}
\label{lem:stabil}
Impose the spectral assumptions on $L_{+},L_{-}$ from Definition~\ref{def:L+L-}.
Then there exist constants $c=c(\alpha,\beta)>0$ such that
for all $\vec\psi\in H^1(\R^n)$, 
\begin{enumerate}
\item $\la L_{-}\psi_2,\psi_2\ra \ge c\|\psi_2\|_2^2$ \ \  if\ \  $\psi_2\perp\partial_\alpha\phi$, 
$\psi_2\perp \partial_j\phi$
\item $\la L_{+}\psi_1,\psi_1\ra \ge c\|\psi_1\|_2^2$ \ \ if\ \  $\psi_1\perp\phi$, 
$\psi_1\perp x_j\phi$.
\end{enumerate}
The constant $c(\alpha,\beta)$ can be taken to be uniform in $\alpha$ in
the following sense: If $\alpha_0$ satisfies Definition~\ref{def:L+L-},
then there exists $\delta>0$ so that 1.\ and 2.\ above hold for all 
$|\alpha-\alpha_0|<\delta$ with $c(\alpha,\beta)>\half c(\alpha_0,\beta)$.
\end{lemma}
\begin{proof} Consider the minimization problems
\bea 
\inf_{f\in H^1} \la L_{-}f,f\ra && \text{\ subject to constraints\ \ } \|f\|_2=1,\;f\perp\partial_\alpha \phi,\; f\perp \partial_j\phi \label{eq:L-min} \\
\inf_{g\in H^1} \la L_{+}g,g\ra && \text{\ subject to constraints\ \ } \|g\|_2=1,\;g\perp \phi,\;
g\perp x_j\phi. \label{eq:L+min}
\eea
As usual, one would like to establish the existence of minimizers by means of passing
to weak limits in minimizing sequences. While such sequences are bounded in $H^1(\R^n)$,
this is not enough to guarantee strong convergence in $L^2(\R^n)$ because some
(or all) of the $L^2$-mass might escape to infinity. Using the fact that the quadratic
forms in question are perturbations of $\half\|\nabla f\|_2^2 + \frac{\alpha^2}{2}\|f\|_2^2$
by a potential that decays at infinity, one can easily exclude that {\em all} the $L^2$-mass escapes
to infinity. One then proceeds to show that the remaining piece of the limit, 
normalized to have $L^2$-norm one,  is a minimizer. This, however, 
 is a simple consequence of the non-negativity of $L_{-}$ and $L_{+}$,
the latter under the constraint $f\perp \phi$, see Lemma~\ref{lem:E1} above.
 This argument is presented in all details in~\cite{W1}, page~478 for the case of
power nonlinearities. But the same argument also applies to the general nonlinearities
 considered here, and we do not write it out. 

Assume therefore that $f_0$ is a minimizer of~\eqref{eq:L-min} with $\|f_0\|_2=1$,
$f_0\perp \partial_\alpha\phi$, $f_0\perp \partial_j\phi$. 
Then 
\be
\label{eq:euler1}
 L_{-}f_0=\lambda_0 f_0 + c_0\partial_\alpha\phi + \sum_{j=1}^n c_j \partial_j\phi
\ee
for some Lagrange multipliers $\lambda_0,c_0,\ldots,c_n$. Clearly, $\lambda_0$
agrees with the minimum sought, and therefore it suffices to show that $\lambda_0>0$.
If $\lambda_0=0$, then taking the scalar product of~\eqref{eq:euler1} with~$\phi$
implies that $c_0=0$ (using that $\la \partial_\alpha \phi,\phi\ra >0$). Taking
scalar products with $x_k\phi$ shows that also $c_k=0$ for $1\le k\le n$. 
Thus $L_{-}f_0=0$, which would imply that $f_0=\gamma \phi$ for some $\gamma\ne0$.
However, this is impossible because of $f_0\perp\partial_\alpha\phi$.

Proceeding in the same manner for $L_{+}$, one arrives at the Euler-Lagrange equation
\[ L_{+}g_0 = \lambda_0 g_0 + c_0\phi + \sum_{j=1}^n c_j x_j \phi.\]
 As before, $\lambda_0$ is the minimum on the left-hand side of~\eqref{eq:L+min} 
and thus $\lambda_0\ge0$ by Lemma~\ref{lem:E1}.
If $\lambda_0=0$, then taking scalar products with $\partial_k\phi$ leads to
$c_k=0$ for all $1\le k\le n$. Hence $L_{+}g_0=c_0\phi$ which implies that
\[  g_0=-\frac{c_0}{\alpha}\partial_\alpha\phi+\sum_{\ell=1}^n b_\ell \partial_\ell\phi.\]
Taking scalar products of this line with $\phi$ and $x_j\phi$ shows
that $c_0=0$ and $b_\ell=0$ for all $1\le\ell\le n$, respectively. 
But then $g_0=0$ which is impossible. 

Since the constants $c(\alpha,\beta)>0$ were obtained by contradiction,
one has no control on their dependence on~$\alpha$. However, let 
$|\alpha-\alpha_0|<\delta$ be as in Definition~\ref{def:L+L-}. 
Suppose $\|f\|_2=1$ satisfies $f\perp \phi(\cdot,\alpha)$, 
$f\perp x_j\phi(\cdot,\alpha)$. Then there is 
\[
h\in\spa\Big\{\phi(\cdot,\alpha),\,\phi(\cdot,\alpha_0),\,
x_j\phi(\cdot,\alpha),\, x_j\phi(\cdot,\alpha_0)\::\: 1\le j\le n\Big\}
\]
so that $f+h\perp\phi(\cdot,\alpha_0)$ and 
$f+h\perp x_j\phi(\cdot,\alpha_0)$. Moreover, since 
$\|\partial_\alpha\phi\|_{H^1(\R^n)}+\|\partial^2_\alpha\phi\|_{H^1(\R^n)} <\infty$
one can take $\|h\|_{H^1(\R^n)}$ as small as desired provided $\delta$ 
is chosen small enough. One can therefore use inequality 2.\ from
this lemma at $\alpha_0$ for $f+h$ to obtain a similar bound for $f$ at~$\alpha$. 
\end{proof}

The following corollary proves the crucial {\em linear stability assumption}
contingent upon the spectral assumptions on $L_{+},L_{-}$ from above
(and thus, in particular, contingent upon the convexity condition).
Strictly speaking, the following corollary gives a stronger statement
than~\eqref{eq:stab_cond}, since the range of $P_s$ is potentially 
smaller than needed for the stability to hold. 

\begin{cor}
\label{cor:wein}
Impose the spectral assumptions on $L_{+},L_{-}$ from Definition~\ref{def:L+L-}.
Then there exist constants $C=C(\alpha,\beta)<\infty$ so that 
for all $\vec\psi_0\in H^1(\R^n)$
\be 
\label{eq:H^1stab}
\|e^{itA}\vec\psi_0\|_{H^1(\R^n)} \le C\|\vec\psi_0\|_{H^1(\R^n)}
\text{\ \ provided\ \ } \vec\psi\in\calM_*^\perp.
\ee
Here $\calM_*$ is the $A^*$-invariant subspace from~\eqref{eq:rootA*}. 
Moreover, the same bound
holds for $H^s(\R^n)$-norms for any real~$s$ with $s$-dependent constants
(and thus in particular for $L^2(\R^n)$).
Analogous statements hold for $e^{itA^*}$. Finally, the constants $C(\alpha,\beta)$
can be taken to be uniform in $\alpha$ in
the following sense: If $\alpha_0$ satisfies Definition~\ref{def:L+L-},
then there exists $\delta>0$ so that~\eqref{eq:H^1stab} holds for all 
$|\alpha-\alpha_0|<\delta$ with $C(\alpha,\beta)<2 C(\alpha_0,\beta)$.
\end{cor}
\begin{proof} Let $\vec\psi(t)=e^{itA}\vec\psi_0$. 
By Lemma~\ref{lem:stabil} one has
\[ Q(\vec\psi_0)=Q(\vec\psi(t))\ge c\,\|\vec\psi(t)\|^2_{L^2(\R^n)}\]
provided that $\vec\psi_0\in\calM_*^\perp$.
Since clearly $Q(\vec\psi_0)\le C\|\vec\psi_0\|^2_{H^1(\R^n)}$ one concludes 
that~\eqref{eq:H^1stab} holds with $L^2(\R^n)$ on the left-hand side.
In order to pass to~$H^1(\R^n)$ write
\bea
 \la L_{-}f,f\ra &=& (1-\eps)\la L_{-}f,f\ra + \frac{\eps}{2}\|\nabla f\|_2^2 
+\eps\frac{\alpha^2}{2} \|f\|_2^2 -  
\eps \int_{\R^n} \beta(\phi^2(x)) |f(x)|^2\,dx \nn \\
&\ge& \frac{\eps}{2}\|\nabla f\|_2^2 + c(1-\eps)\|f\|_2^2 + 
\eps(\alpha^2/2-\|\beta\|_\infty)\|f\|_2^2,
\label{eq:H1split}
\eea
where the constant $c$ in~\eqref{eq:H1split} is the one from Lemma~\ref{lem:stabil}. 
Taking $\eps$ small enough, one sees that the third term can be absorbed into
the second. Thus the entire right-hand side of~\eqref{eq:H1split} admits the lower
bound~$\frac{\eps}{2}\|f\|_{H^1(\R^n)}^2$. The same argument applies to~$L_{+}$,
and~\eqref{eq:H^1stab} follows. The uniformity statement concerning the
constants $C(\alpha,\beta)$ is an immediate consequence of the analogous statement
in Lemma~\ref{lem:stabil}.
To obtain~\eqref{eq:H^1stab} for all $H^s$ spaces
note first that 
\be
\label{eq:Mell}
C_\ell^{-1}\,\|\vec\psi\|_{H^{2\ell}(\R^n)} \le \|(A+iM)^\ell \vec\psi\|_2 
\le C_\ell\,\|\vec\psi\|_{H^{2\ell}(\R^n)} 
\ee
for all integers $\ell$ and sufficiently large  $M=M(\ell)$. 
Indeed, to check the lower bound for $\ell=1$
one can use $(A+iM)^{-1}=(B+iM)^{-1}\Big[1+V(B+iM)^{-1}\Big]^{-1}$.
The inverse of the operator in brackets exists provided $M$ is large and it is a
bounded operator on~$L^2(\R^n)$. Taking powers of this relation allows one
to deal with all $\ell\ge1$ (in our case $V$ is $C^\infty$ which is needed here). 
Since $\calM_*$ is $A^*$-invariant and therefore $\calM_*^\perp$ is $A$-invariant, 
inserting~\eqref{eq:Mell} into~\eqref{eq:H^1stab} allows one to pass 
to all odd integers~$s$. The case of general~$s$ then follows by interpolation.
Finally, since all arguments in this section apply equally well to $A^*$ as~$A$,
the corollary follows. 
\end{proof}

This corollary has an important implication concerning the 
structure of the root spaces as required in Proposition~\ref{prop:null}. 

\begin{cor}
\label{cor:null_struc}
Impose the spectral assumptions on $L_{+},L_{-}$ from Definition~\ref{def:L+L-}.
Then equality holds in the relation concerning the root spaces~\eqref{eq:rootA} 
and~\eqref{eq:rootA*}.
In particular, one has $\ker(A^2)=\ker(A^3)$ and $\ker((A^*)^2)=\ker((A^*)^3)$.
\end{cor}
\begin{proof}
Suppose $\dim(\calN(A))>2n+2$. Then there exists $\vec\psi_0\in\calN(A)$
such that $\vec\psi_0 \in\calM_*^\perp$. This is because a system of
$2n+2$ equations in $2n+3$ variables always has a nonzero solution.
Since $\partial_\alpha\phi
\not\perp \phi$ and $\partial_j\phi\not\perp x_j\phi$,
one checks that $\vec\psi_0\not\in \ker(A)$.  Therefore, 
$\vec\psi_0\in\ker(A^k)\setminus\ker(A^{k-1})$ for some $k\ge2$. 
Expanding $e^{itA}$ into a series implies that 
$\|e^{itA}\vec\psi_0\|_2>c\, t^{k-1}$ for some constant $c>0$,
which contradicts Corollary~\ref{cor:wein}. 
Therefore, $\dim(\calN(A))\le 2n+2$. Since moreover $\phi>0$
and $\la \partial_\alpha\phi,\phi\ra>0$ imply that the $2n+2$
vectors on the right-hand side of~\eqref{eq:rootA} 
are linearly independent, equality must hold as claimed.
Analogously for~\eqref{eq:rootA*}.
\end{proof}
\section {Uniqueness of the negative eigenvalue of $L_+$}
In this section we show that the assumption that the linear operator $L_+$ 
has a single simple negative eigenvalue of Definition~\ref{def:L+L-}
is automatically satisfied for a large class of nonlinearities. 

We consider the ground state of the problem
\be
\label{eq:stat2}
-\frac 12\Laplace\phi - \beta(|\phi|^{2})\phi + \frac{\alpha^{2}}2\phi = 0
\ee
constructed from the constrained minimization problem for the following functional:
\be
\label{eq:constraint2}
J[u] = \big\{\int_{\R^{n}} |\nabla u|^{2}\,\,:\,\, W[u]=\int_{\R^{n}} 
G(u) = 1\big\},  
\ee
where $
G(s) = \int_{0}^{s} \big (
\beta(s^{2}) s  - \frac{\alpha^{2}}2 s\big ).$
If $w$ is a minimum it solves the equation 
\be
\label{eq:EuL}
-\frac 12\Laplace w - \lambda (\beta(w^{2}) w - \frac{\alpha^{2}}2 w =0
\ee
where the Lagrange multiplier $\lambda$ is determined from the condition 
that $W[w]= 1$. We can then find a ground state via rescaling 
\be
\label{eq:rescale2}
\phi(x) = w(\lambda^{-\frac 12} x)
\ee 
for an appropriate choice of $\lambda$.
In this section we prove the following:
\begin{theorem}
Let $\phi$ be the ground state of \eqref{eq:stat2} constructed in the 
variational problem \eqref{eq:constraint2} with the nonlinearity 
$\beta$ verifying the condition that $|\beta(s^2) \le C |s|^{p-1}$ with 
$p\in (1,  1+\frac 4{n-2})$ and that there exists $s_0$ such that 
$C(s_0) > 0$. Moreover, assume that 
$$
\beta (s^2) s^2 \ge 2 \int_0^s \beta(\tau^2) \tau\, d\tau
$$
for all $s\ge 0$. Then the operator $L_{+}=-\frac 12 \Delta -\beta(\phi^2)-\beta'(\phi^2)\phi^2 + \frac {\alpha^2}2$ has a unique simple negative 
eigenvalue.
\end{theorem}
\begin{proof}
Let us consider a smooth family of functions $w_{z}$ satisfying
the constraint $W[w_{z}]=1$ for all $z$ and passing through a
minimizer $w$, i.e., $w_{0} = w$.
Observe that since $w$ is a minimizer of the constrained variational 
problem we have that 
$$
\frac {d}{dz} J[w_{z}] |_{z=0}=0,\qquad
\frac {d^{2}}{dz^{2}} J[w_{z}] |_{z=0}\ge 0 
$$
Denote 
$$
\dot w = \frac {d}{dz} w_{z} |_{z=0}, \qquad
\ddot w = \frac {d^{2}}{dz^{2}} w_{z}|_{z=0}
$$
Then
$$
0=\frac {d}{dz} J[w_{z}]|_{z=0} = 2 \int_{\R^{n}}\nabla w 
\nabla\dot w
 $$
and 
\begin{equation}
\label{eq:secvar}
0\le \frac {d^{2}}{d z^{2}} J[w_{z}]|_{z=0} = 
2\int_{\R^{n}} |\nabla \dot w|^{2} + \nabla w \nabla \ddot w
= - 2\int_{\R^{n}} \Laplace\dot w \dot w + \Laplace w \ddot w
\end{equation}
Moreover, since the family $w_{z}$ verifies the constraint 
$W[w_{z}] =1$ we have 
$$
\frac {d}{d z} W[w_{z}]= \frac {d^{2}}{d z^{2}} 
W[w_{z}] = 0
$$
Computing the derivatives we obtain 
\begin{align}
&\int_{\R^{n}} (\beta (w^{2}) w - \frac {\alpha^{2}}2 w) \dot w =0,
\label{eq:forth}\\
&\int_{\R^{n}} (\beta (w^{2})  + 2\beta'(w^{2}) w^{2} - 
\frac {\alpha^{2}}2 ) |\dot w |^{2} + (\beta(w^{2}) w -\frac 
{\alpha^{2}}2 w) \ddot w =0
\label{eq:sorth}
\end{align}
Now observe that the Euler-Lagrange equation \eqref{eq:EuL} implies that 
$$
-\int_{\R^{n}} \Laplace w \ddot w = 
2\lambda \int_{\R^{n}} (\beta(w^{2}) w - \frac{\alpha^{2}}2 w) \ddot w
$$
On the other hand from \eqref{eq:sorth} we have that 
$$
\int_{\R^{n}} (\beta(w^{2}) w -\frac 
{\alpha^{2}}2 w) \ddot w = - \int_{\R^{n}} 
(\beta (w^{2})  + 2\beta'(w^{2}) w^{2} - 
\frac {\alpha^{2}}2 ) |\dot w |^{2}
$$
Substituting this into \eqref{eq:secvar} we obtain that 
\be
\label{eq:sub}
- 2\int_{\R^{n}} \Big (\Laplace\dot w  + 2 \lambda 
(\beta (w^{2})  + 2\beta'(w^{2}) w^{2} - 
\frac {\alpha^{2}}2 ) \dot w \Big )\dot w \ge 0
\ee
Recall the definition of the operator $L_{+}$ associated with a 
ground state $\phi$:
$$
L_{+} = - \frac 12\Laplace  - \beta (\phi^{2})  - 2\beta'(\phi^{2}) 
\phi^{2} + \frac {\alpha^{2}}2 
$$
We know from \eqref{eq:rescale2} that
$\phi(x) = w(\lambda^{-\frac 12} x)$. Thus rescaling \eqref{eq:sub}
leads to the inequality
\be
\label{eq:L-}
\lambda^{\frac n2+1} \la L_{+} \dot w(\lambda^{-\frac 12}\cdot )\,\,,\, 
\dot w (\lambda^{-\frac 12}\cdot )\ra \ge 0
\ee
Furthermore, \eqref{eq:forth} implies that 
the function $\dot w(\lambda^{-\frac 12}\cdot )$ can be chosen arbitrarily
from the subspace orthogonal to the function
$$
\Psi = (\beta (\phi^{2}) \phi - \frac {\alpha^{2}}2 \phi)
$$
Let $\Pi$ denote the orthogonal projection on the above
subspace. Then \eqref{eq:L-} implies\footnote{It is easy to see
that the condition that ensures that $\lambda > 0$ is that 
$$
\beta(s^{2}) s^{2} \ge 2\int_{0}^{s}\beta(\tau^{2}) \tau\,d\tau
$$}
that the operator 
$\Pi L_{+} \Pi \ge 0$. 
Now let $\mu$ and $u$ be correspondingly 
a negative eigenvalue of $L_{+}$ and the associated eigenfunction.
Decompose\footnote{One can always achieve this decomposition 
by normalizing $u$} 
$$
u= v +\Psi 
$$
where $\Pi v = v$.
Then 
$$
L_{+} u = L_{+} v + L_{+}\Psi = \mu (v+\Psi)
$$
Projecting we obtain that 
$$
(\Pi L_{+}\Pi-\mu) v = -\Pi L_{+}\Psi
$$
Since the operator $\Pi L_{+} \Pi$ is non-negative 
and $\mu <0$ we have that $ (\Pi L_{+}\Pi-\mu)$ is invertible and
$$
v =- (\Pi L_{+}\Pi-\mu)^{-1} \Pi L_{+}\Psi
$$
Therefore the eigenvalue $\mu$ is simple.
We now substitute $v$ back into the eigenvalue equation for $u$
and take a scalar product with $\Psi$. We obtain that
$$
-\langle L_+ (\Pi L_{+}\Pi-\mu)^{-1} \Pi L_{+}\Psi, \Psi\rangle 
+ \langle L_+\Psi, \Psi\rangle = \mu \langle \Psi, \Psi\rangle
$$
We now assume that there is another negative eigenvalue $\lambda$.
Then  the corresponding equation also holds with $\lambda$ replacing 
$\mu$. Subtracting we obtain
$$
-\langle L_+(\Pi L_{+}\Pi-\mu)^{-1} \Pi L_{+}\Psi, \Psi\rangle 
+\langle L_+ (\Pi L_{+}\Pi-\lambda)^{-1} \Pi L_{+}\Psi, \Psi\rangle =
(\mu -\lambda) \langle \Psi, \Psi\rangle
$$
Using the resolvent identity
$$
 (\Pi L_{+}\Pi-\mu)^{-1} -  (\Pi L_{+}\Pi-\lambda)^{-1} =
 (\mu-\lambda)  (\Pi L_{+}\Pi-\mu)^{-1} (\Pi L_{+}\Pi-\lambda)^{-1}
$$
Thus,
$$
\langle  (\Pi L_{+}\Pi-\mu)^{-1}(\Pi L_{+}\Pi-\lambda)^{-1}
 \Pi L_{+}\Psi, \Pi L_+\Psi\rangle = -  \langle \Psi, \Psi\rangle
$$
Observe that $ (\Pi L_{+}\Pi-\mu)^{-1}(\Pi L_{+}\Pi-\lambda)^{-1}$ is 
a positive operator since $\Pi L_+\Pi$ is non-negative and 
$\lambda, \mu < 0$. Contradiction.
\end{proof}

\section{Verification of the needed properties of $L_+$ and $L_-$ 
for a particular class of nonlinearities}

The purpose of this section is to verify that the conditions
imposed on $L_{+}$ and $L_{-}$ in Definition~\ref{def:L+L-}
hold for  a particular class of nonlinearities. This class
is chosen so as to allow for the nonlinear analysis of the previous
sections to go through.
\be
\label{eq:beta'}
\beta_{\theta}(s^{2}) = s^{p-1} \frac {f(s^2)}{\theta + f(s^2)}
\ee
with a constant $\theta >0$ and the function $f$ satisfying the conditions
\be\label{eq:condit-f'}
C_1 s^{r+1}\le |f(s^2)\le C_2 s^{r+1}, \qquad 
s^2 |f'(s^2)|\le C_2 |f(s^2)|, \qquad p\in (1, 1+\frac 4n),\,\,
r>-1
\end{equation} 
The existence of a ground state for the problem 
\be
\label{eq:stat}
-\frac 12\Laplace\phi - \beta(|\phi|^{2})\phi + \frac{\alpha^{2}}2\phi = 0
\ee
for $\alpha\ne 0$ had been established by Berestycki and Lions under 
the following
conditions on the function $\beta$:
\begin{enumerate}
\item
$0\ge \overline{\lim}_{s\to +\infty} {\beta(s)}{s^{-\frac 2{n-2}}}
\ge +\infty$
\item There exists $s_{0}>0$ such that $G(s_{0}) = \int_{0}^{s_{0}}
\beta(s^{2}) s \,ds - \frac{\alpha^{2}}4 s^{2}_{0} > 0$
\end{enumerate}
Moreover, in the case when the function $\beta(s)$ satisfies a 
stronger condition that 
\be
\label{eq:nonlcon}
\lim_{s\to +\infty} \beta(s) s^{{-\frac 2{n-2}}} = 0
\ee
a ground state can be constructed from a solution of the 
constrained minimization problem for the following functional:
\be
\label{eq:constraint}
J[u] = \big\{\int_{\R^{n}} |\nabla u|^{2}\,\,:\,\, W[u]=\int_{\R^{n}} 
G(u) = 1\big\} 
\ee
If $w$ is a minimum it solves the equation 
$$
-\frac 12\Laplace w - \lambda (\beta(w^{2}) w - \frac{\alpha^{2}}2 w) =0
$$
where the Lagrange multiplier $\lambda$ is determined from the condition 
that $W[w]= 1$. We can then find a ground state via rescaling 
\be
\label{eq:rescale}
\phi(x) = w(\lambda^{-\frac 12} x)
\ee 
Observe that it is possible to choose $w$ a positive spherically 
symmetric function.
We now consider the case of the monomial subcritical nonlinearity
$\beta(s) = s^{\frac {p-1}2}$ with $p < \frac {n+2}{n-2}$. 
By the results of Coffman, McLeod-Serrin, and Kwong there exists 
a unique positive radial solution of the equation \eqref{eq:stat}
for $\alpha\ne 0$. Let $w$ denote the corresponding minimizer of 
the functional $J$.
\begin{defi}
\label{de:dega}
Given $\ga > 0$ and $w$ the minimizer of $J$ corresponding to the
unique ground state $\phi$ define 
\begin{align}
\sigma(\ga) = \inf \big \{\theta\,:\,\, &{\text{
for any positive non-increasing radial function}}\,\, 
u\,\, {\text{ with the}}\nn\\ &{\text{property that}}\,\, 
\|u- w\|_{H^{1}}\ge \theta 
\,\,{\text{ we have
that}}\,\, J[u]\ge J[w] + \ga\big \} \label{eq:defde}
\end{align}
\end{defi}
We now make the following claim
\begin{lemma}
\label{le:convex}
Let $\beta(s^2)=s^{p-1}$ be a monomial nonlinearity with $p\in (1, \frac 4{n-2})$. 
Then the 
function $\sigma (\ga)\to 0$ as $\ga\to 0$.
\end{lemma}
\begin{proof}
We argue by contradiction. Assume that there exists a 
sequence $\ga_{k}\to 0$, a positive constant $\sigma$, and 
positive radial functions $ u_{\ga_{k}}$ 
such that $ \|u_{\ga_{k}}- w\|_{H^{1}}\ge \theta$ but 
$J[u_{\ga_{k}}] < J[w] +\ga_{k}$.  Then the sequence 
$u_{\ga_{k},\theta}$ is minimizing for the functional $J$. 
This implies that 
$$
\|\nabla u_{\ga_{k}}\|_{L^{2}}\to \|\nabla w\|_{L^{2}} 
$$
Using the constraint $W[u_{\ga_{k}}]=1$ it is not difficult
to show that the sequence $u_{\ga_{k}}$ is uniformly bounded
in $H^{1}$, see \cite{BL}. Thus without loss of generality we assume that 
$u_{\ga_{k}}\to u$ weakly in $H^{1}$ for some radial 
non-increasing function $u$.
Therefore, $u$ is another minimizer of the functional $J$ and its
rescaled version is a non-increasing radial solution of the 
equation \eqref{eq:stat}. By the strong maximum principle it is 
positive\footnote{The minimizer $u$ cannot be identically zero 
since one can show that the minimum is attained on the function
satisfying the constraint $W[u]=1$, see \cite{BL}.}
 and therefore a ground state. Since the ground state is 
unique, after rescaling back we conclude that $u=w$.
Therefore, we have constructed a sequence $u_{\ga_{k}}$
with the properties that
\begin{align}
&u_{\ga_{k}}\to w \quad{\text{weakly in}}\quad 
H^{1},\label{eq:weak}\\
&\nabla u_{\ga_{k}}\to \nabla w \quad {\text{in}}\quad 
L^{2},\label{eq:nabstrong}\\
&\int_{\R^{n}} |u_{\ga_{k}}|^{p+1} = 1 + \frac {\alpha^{2}}4
\int_{\R^{n}} |u_{\ga_{k}}|^{2},\label{eq:Vconstr}\\
&\|u_{\ga_{k}}- w\|_{H^{1}}\ge \theta\label{eq:far}
\end{align}
Since  $2<p+1\le \frac {2n}{n-2}$, conditions \eqref{eq:weak}
and \eqref{eq:nabstrong}	 
imply that 
$$
\int_{\R^{n}} |u_{\ga_{k}}|^{p+1} \to 
\int_{\R^{n}} |w|^{p+1}
$$
Thus from \eqref{eq:Vconstr}
$$
\int_{\R^{n}} |u_{\ga_{k}}|^{2}\to 
\int_{\R^{n}} |w|^{2}
$$
and with the help of \eqref{eq:weak} and \eqref{eq:nabstrong} we
conclude that $u_{\ga_{k}}\to w$ in $H^{1}$. This contradicts
\eqref{eq:far}. 
\end{proof}
\subsection {Variational continuity}
We now consider the ground state problem 
\be
-\frac 12\Laplace\phi_{\theta} - \beta_{\theta}(|\phi_{\theta}|^{2})\phi_{\theta} + 
\frac{\alpha^{2}}2\phi_{\theta} = 0
\tag{gr$_{\theta}$}
\ee
for the nonlinearities
\be
\label{eq:betaep} 
\beta_{\theta}(s^{2}) = -|s|^{{p-1}} \frac {f(s^2)}{\theta + f(s^2)},
\ee
where $1<p<1+\frac 4n$, the function $f$ satisfies the following estimate
$$
C_1 s^q\le |f(s^2)|\le C_2 s^{q+1}
$$ 
for some positive constants $C_1, C_2$ and $q>-1$, and $\theta >0$ is a positive number.
Define 
\begin{align}
G_{\th} (\tau) = \int_{0}^{\tau} \beta_{\th}(s^{2})s\,ds - 
\frac{\alpha^{2}}4 \tau^{2},\label{eq:Gep}\\
W_{\th} [u] = \int_{\R^{n}}G_{\th}(u(x))\,dx\label{eq:Wep}
\end{align}
\begin{lemma}
We have the following estimate
\be
\label{eq:diffep}
|W_{\th}[u] - W_{0}[u]|\les \th^{\frac {p-1}{p+q}} 
\int_{\R^{n}}\big (|u|^{2} + |u|^{p+1}\big )
\ee
\end{lemma}
\begin{proof}
Estimate \eqref{eq:diffep} immediately follows from the 
inequality
\be
\label{eq:Gtep}
|G_{\th}(\tau)-G_{0}(\tau)|\le \tau^{p+1}\frac {\th}{\th + 
\tau^{q+1}} = \th\tau^{p-q} \frac {\tau^{q+1}}{\ep + 
\tau^{q+1}}
\ee
since the above expression can be bounded by
$$
\min \big \{\th\tau^{p-q}, \tau^{p+1}\big \}
$$
Thus using the first term for the values of $\tau\ge \th^{\frac 1{p+q}}$, and the second term 
when $\tau \le \th^{\frac 1{p+q}}$  we obtain \eqref{eq:Gtep}. 
\end{proof}
We now consider the variational problem 
\be
\label{eq:Jep}
J_{\th}[u] = \big\{ \int_{\R^{n}} |\nabla u|^{2}:\,\,
W_{\th}[u] =1\big\}
\ee
\begin{prop}
\label{prop:Contin}
Let $\phi$ be the  ground state 
of the problem (gr$_{0}$). Then for any 
sufficiently small $\th >0$ there exists a  positive constant $\de'=\de'(\th)
\to 0$ as $\th\to 0$, 
and a ground state 
$\phi_{\th}$ of (gr$_{\th}$) 
such that $\|\phi_{\th} -\phi\|_{H^{1}}<\de'$. 
\end{prop}
\begin{proof}
We start by choosing a sufficiently large constant $M$ such that 
for all sufficiently small $\th$ any minimizer of $J_{\th}$ is 
contained in a ball $B_{M/2}$ of radius $M/2$ in the space $H^{1}$.
In particular, using \eqref{eq:diffep} we will assume that 
for $u\in B_{M}$
\be
\label{eq:dep}
|W_{\th}[u] - W_{0}[u]|\les \ep^{\frac {p-1}{p+q}} 
\ee
We now observe the following trivial property of the constraint 
functionals $W_{\th}[u]$: for any $\th\ge 0$ and an arbitrary
$\mu\ne 0$
\be
\label{eq:rescep}
W_{\th}[u(x)] = \mu^{n} W_{\th}[u(\frac x\mu)]
\ee
We now fix a sufficiently small $\th >0$.
Let $w$ be the minimizer of the variational problem $J=J_{0}$ 
corresponding to the unique ground state $\phi$. The function 
$w$ satisfies the constraint $W_{0}[w]=1$. Therefore, using the
rescaling property \eqref{eq:rescep} and \eqref{eq:dep} we can show
that there exists $\mu=\mu(w)$ with the property that 
\begin{align}
&W_{\th}[w(\frac x\mu)]=1,\nn\\
&|\mu-1|\le \th^{\frac {p-1}{p+q}}\label{eq:muep}
\end{align} 
Moreover,
\be
\label{eq:Jeep}
J_{\th}[w(\frac x\mu)] = \mu^{n-2} J_{0}[w(x)]=
J_{0}[w] + O(\th^{\frac {p-1}{p+q}})
\ee
We now claim that there exists a small positive $\de=\de(\th)\to 0$
as $\th\to 0$,
such that for any positive non-increasing radial function 
$u$ satisfying the constraint $W_{\th}[u]=1$ and the property 
that 
\be
\label{eq:uw}
\|u-w(\frac x\mu)\|_{H^{1}}\ge \de
\ee
we have
\be
\label{eq:nonmin}
J_{\th}[u]\ge J_{\th}[w(\frac x\mu)] + \th^{\frac {p-1}{p+q}}
\ee
Assume for the moment that the claim holds. Then \eqref{eq:uw} and 
\eqref{eq:nonmin} imply that $J_{\th}$ has a minimizer in 
the $\delta$ neighborhood of the function $w(\frac x\mu)$.
We denote this minimizer by $w_{\th}$. Then 
\eqref{eq:muep} implies that 
$$
\|w_{\th}-w\|_{H^{1}}\le \|w_{\th}- w(\frac x\mu)\|_{H^{1}}
+ \|w - w(\frac x\mu)\|_{H^{1}}\le \de + \|w - w(\frac x\mu)\|_{H^{1}}
$$
Observe that $\|w - w(\frac x\mu)\|_{H^{1}}\to 0$ as $\mu\to 1$,
which follows by the 
density argument and the fact that it is easily satisfied on 
functions of compact support\footnote{In fact, the minimizer $w$ is 
smooth and localized in space and thus one could even give the precise 
dependence on $\mu$}. 
Define the function $a_{w}(\ep)$:
\be
\label{eq:modulus}
a_{w}(\ep):= \sup_{|\mu-1|\le \th^{\frac {p-1}{p+q}}}
\|w - w(\frac x\mu)\|_{H^{1}},\qquad a_{w}(\th)\to 0
\quad {\text{as}}\,\,\,\th\to 0 
\ee
Therefore,
\be
\label{eq:close}
\|w_{\th}- w \|_{H^{1}}
\les \de + a_{w}(\th)
\ee
The functions $w_{\th}$, $w$ are the solutions of the 
Euler-Lagrange equations
\begin{align}
&-\frac 12\Laplace w_{\th} - \lam_{\th}\big (\beta_{\th}(w_{\th})^{2})w_{\th}
-\frac{\alpha^{2}}2 w_{\th}\big ) = 0,\label{eq:ELep}\\
&-\frac 12\Laplace w - \lam\big (\beta_{0}(w)^{2})w
 - \frac{\alpha^{2}}2 w\big ) = 0,\label{eq:EL}
\end{align}
where the Lagrange multipliers $\lam_{\th}$, $\lam$ are determined
from the conditions that 
$W_{\th}[w_{\th}]=W_{0}[w]=1$.
We multiply the equations \eqref{eq:ELep} and \eqref{eq:EL} by
$w_{\th}$ and $w$ correspondingly, integrate by parts, and subtract
one from another.
Using the estimate
$$
\int_{\R^{n}}\big 
|\beta_{\th}(w_{\th})^{2})w^{2}_{\th} - 
\beta_0(w_{\th})^{2})w_{\th}^{2}\big |\les 
\th^{\frac {p-1}{p+q}},
$$
which is essentially the same as the estimate 
\eqref{eq:dep}, and the estimate \eqref{eq:close} we obtain that
\be
(\lam -\lam_{\th}) \int_{\R^{n}}(\beta_{0}(w)^{2})w^{2}
-\frac{\alpha^{2}}{2}w^{2}\big )  = O(\de^2) + O(\th^{\frac {2(p-1)}{p+q}})
+ a^2_{w}(\th)
\ee
Recall that $\beta_{0}(w^{2}) = w^{p-1}$. The condition that
$W[w]=1$ implies that 
$$
\int_{\R^{n}}\big ( \frac 1{p+1} |w|^{p+1} -\frac {{\alpha}^{2}}4
|w|^{2}\big) = 1
$$
Thus,
$$ 
\int_{\R^{n}}(\beta_{0}(w)^{2})w^{2}
-\frac{\alpha^{2}}2w^{2}\big ) = 2 + \frac {p-1}{p+1} 
\int_{\R^{n}} |w|^{p+1} \ge 2
$$
This allows us to conclude that 
\be
\label{eq:lep}
|\lam -\lam_{\th}|\le \de^2 + \th^{\frac {2(p-1)}{p+q}} + a^2_{w}(\th)
\ee
Finally, recall that the ground states $\phi_{\th}$ and 
$\phi$ are obtained by the rescaling of the minimizers
$w_{\th}$ and $w$.
$$
\phi_{\th} (x) = w_{\th}(\lam_{\th}^{-\frac 12} x),
\qquad \phi (x) = w(\lam^{-\frac 12} x)
 $$
Thus
\begin{align*}
\|\phi_{\th} - \phi \|_{H^{1}} & \les 
\|w_{\th}(\lam_{\th}^{-\frac 12} x) -w(\lam^{-\frac 12} x)\|_{H^{1}}
\\ &= \Big (\lam_{\th}^{\frac n2}  + \lam_{\th}^{\frac n2-1}  \Big)
\|w_{\th}(x) -
w\big ((\frac {\lam_{\th}}{\lam})^{\frac 12}x\big )\|_{H^{1}}
\\ &\le \Big (\lam_{\th}^{\frac n2}  + \lam_{\th}^{\frac n2-1}\Big )  
\|w_{\th}(x) - w(x)\|_{H^{1}} + 
\|w(x) - w\big ((\frac {\lam_{\th}}{\lam})^{\frac 12}x\big )\|_{H^{1}}
\end{align*}
By \eqref{eq:lep} the constants $\lam_{\th}$ are uniformly bounded 
in terms of the absolute constant $\lam$, which depends only on $w$.
Moreover, $\lam_{\th}\to \lam$ as $\th\to 0$. We appeal again to the 
$H^{1}$ modulus of continuity of the minimizer $w$ and define the function
\footnote{One can show that $\alpha\ne 0$ the Lagrange multiplier 
$\lam\ne 0$. This follows from the following argument. By interpolation
for $p\le \frac {n+2}{n-2}$
$$
\int w^{p+1} \le \|\nabla w\|_{L^{2}}^{n\frac {p-1}2} 
\|w\|_{L^{2}}^{p+1-n\frac {p-1}2}
$$
Thus for $n > 2$ the power $p+1-n\frac {p-1}2 < 2$ and using 
Cauchy-Schwarz, constraint $W[w]=1$ and the assumption that $\alpha\ne 
0$, we can show that $\|\nabla w\|_{L^{2}}\ge c$ for some positive 
constant $c$. Repeating argument determining the Lagrange multiplier
we verify that $\lam\ne 0$}
\be
\label{eq:bmodul}
b_{w}(\th,\de):=\sup_{|\mu- 1|\le C(\de^2 + \th^{\frac {2(p-1)}{p+q}} + a^2_{w}(\th) )}
\|w(x) - w(\mu^{-\frac 12} x)\|_{H^{1}}
\ee
for some positive constant $C$ dependent only on the minimizer $w$.
The function $b_{w}(\th,\de)\to 0$ as $\th,\de\to 0$. 
Therefore, since we have already proved in \eqref{eq:close}
that $w_{\th}$ is close to $w$
in $H^{1}$, we obtain 
\be
\label{eq:final}
\|\phi_{\th} - \phi \|_{H^{1}} \les \de + 
\th^{\frac {p-1}{p+q}} + a_{w}(\th) + b_{w}(\th,\de)
\ee
Since by the claim $\de=\de(\th)\to 0$ as $\th\to 0$ and the functions
$a_{w}(\th)$, $b_{w}(\th,\de)$ also have this property
we obtain the
desired conclusion.

It remains to prove the claim \eqref{eq:uw}, \eqref{eq:nonmin}.
Let $u$ be as in the claim, i.e., $u\in B_{M}$ and
$W_{\th}[u]=1$, and 
\be
\label{eq:faru}
\|u-w(\frac x\mu)\|_{H^{1}}\ge \de.
\ee
for some $\de$ to be chosen below.
Similar to \eqref{eq:muep} we can find a constant $\nu=\nu(u)$
such that 
\begin{align}
&W_{0}[u(\frac x\nu)]=1,\qquad 
J_{0}[u(\frac x\nu)]= J_{\th}[u] + O(\th^{\frac{p-1}{p+q}}),
\label{eq:unu}\\
&|\nu-1|\le \th^{\frac{p-1}{p+q}}\label{eq:nuep}
\end{align}
Using \eqref{eq:muep}, \eqref{eq:faru}, \eqref{eq:nuep}, 
and definition \eqref{eq:modulus} we infer that 
\be
\label{eq:clouw}
\|u(\frac x\nu) - w\|_{H^{1}} \ge 
\| u(\frac x\nu) - w(\frac x{\mu\nu})\|_{H^{1}} -
\| w(\frac x{\nu\mu}) - w(x)\|_{H^{1}}\ge \nu^{\frac n2} \de
- a_{w}(\th)\ge \de - \th^{\frac {p-1}{p+q}}\de - a_{w}(\th)
\ee
We now use Lemma \ref{le:convex} for the variational problem $J=J_{0}$.
This gives a function $\sigma(\ga)$, with the property that 
$\sigma(\ga)\to 0$ as $\ga\to 0$, such that 
for any radial non-increasing positive 
$v$ with the property that 
$\|v-w\|_{H^{1}}\ge \sigma(\ga)$ and $W_{0}[v]=1$ we have 
$J_{0}[v]\ge J_{0}[w] +\ga$. 
We set
$$
\ga = 5\th^{\frac {p-1}{p+q}}, \qquad 
\de = \sigma(\ga) + \th^{\frac {p-1}{p+q}} + a_{w}(\th)
$$ 
It follows from Definition \ref{de:dega} of
$\sigma(\ga)$ and \eqref{eq:clouw} that with these 
choices, function $u(\frac x\nu)$ verifies the inequality
$$
J_{0}[u(\frac x\nu)]\ge J_{0}[w] + 5\th^{\frac {p-1}{p+q}}
$$
Finally, using \eqref{eq:Jeep} and \eqref{eq:unu} we obtain
$$
J_{\ep}[u]\ge J_{\th}[w(\frac x\mu)] + 3\th^{\frac {p-1}{p+q}}
$$
It remains to note that the constant $\de$ in  has been 
chosen 
$$
\de = \sigma (\th^{\frac{p-1}{p+q}}) + \th^{\frac{p-1}{p+q}} +
a_{w}(\th)
$$
and by Lemma \ref{le:convex} and \eqref{eq:modulus} goes to zero as 
$\th\to 0$, as claimed.
\end{proof}
Recall definition of the operator $L_{+}^{\th}$ associated 
with the ground state $\phi_{\th}$.
\be
\label{eq:Laep}
L_{+}^\th = -\frac 12 \Laplace - \beta_{\th}(\phi_{\th}^{2}) - 
2\beta'_{\th}(\phi_{\th}^{2}) \phi_{\th}^{2} + \frac{\alpha^{2}}2
\ee
Denote
\be
\label{eq:Vep}
V_{\th} = \beta_{\th}(\phi_{\th}^{2}) + 
2\beta'_{\th}(\phi_{\th}^{2}) \phi_{\th}^{2}
\ee
Using the definition of $\beta_{\th}$ we compute $V_{\th}$ explicitly
\be
\label{eq:Vexpl}
V_{\th} =  p\phi_{\th}^{p-1} \frac {f(\phi_\th^2)}{\th
+f(\phi_{\th}^{2})} + 
2 \phi_{\th}^{p+1} \frac {\th f'(\phi_{\th}^{2})}{\th +f(\phi_{\th}^{2})}
\ee
Uniform bounds
on ground states $\phi_{\th}$ guaranteed by the Proposition 
\ref{prop:Contin}
imply the following result.
\begin{lemma}
\label{le:Potent}
Let function $f$ obey the assumptions that 
$$
C_1 s^{q+1}\le |f(s^2)|\le C_2 s^{q+1}, \qquad s^2 |f'(s^2)|\le C_1 f(s^2)
$$
for some positive constants $C_1, C_2$ and $q>-1$.
Then for any $p\in (1,1+\frac 4n]$ there exists a $r=r(p)$ in 
the interval $r\in [\frac n2, \infty)$ such that 
\be
\label{eq:potclose}
\|V_{\th} - V_{0}\|_{L^{r}} \to 0, \quad \th \to 0
\ee
\end{lemma}
\begin{proof}
We have a pointwise bound 
$$
|V_{\th} - p \phi_{\th}^{p-1} |\les 
\min \big \{ \phi_{\th}^{p-1}, \th \phi_{\th}^{p-q-2}\big \}\le 
\th^{0+} \phi_{\ep}^{p-1-}
$$
In addition, since $\phi_{\th}\to \phi$ in $H^{1}$ we have 
that $\phi_{\th}^{p-1}\to \phi^{p-1}$ in the space 
$L^{\frac 2{p-1}}\cap L^{\frac {2n}{(n-2)(p-1)}}$. 
Since
$$
|V_{\th} - V_{0}| \le |V_{\th} - p \phi_{\th}^{p-1} | 
+ p |\phi_{\th}^{p-1} - \phi^{p-1}|\les \th^{0+} \phi_{\ep}^{p-1-}+
|\phi_{\th}^{p-1} - \phi^{p-1}|
$$
we obtain the desired conclusion for any $r$ in the 
interval $r\in (\frac 2{p-1}, \frac {2n}{(n-2)(p-1)}]$.
The existence of the Lebesgue exponent $r$ in the desired 
interval now follows from the restriction $p\in (1,1+\frac 4n]$ on the range 
of the exponent $p$.
\end{proof}
\begin{cor}
\label{cor:Ldiffer}
The operators 
\begin{align*}
&(L_{+}^{\th}- L_{+}) (-\Laplace +1)^{-1}:\,\, L^{2} \to L^{2},\\
&(-\Laplace + 1)^{-1} (L_{+}^{\th}- L_{+}) :\,\, L^{2}\to L^{2}
\end{align*}
with the norm converging to $0$ as $\th\to 0$.
\end{cor}
\begin{proof}
The difference $L_{+}^{\th}- L_{+}= V_{\th} - V_{0}$.
The result now follows from Lemma \ref{le:Potent}, 
Sobolev embeddings, and H\"older inequality.
\end{proof}

\subsection {Stability of the convexity condition}
The goal of this section is to prove that the ground states
$\phi_\th$ of the problem 
$$
\frac 12 \Delta \phi +\beta_\th (\phi^2) = \frac {\alpha^2}2 \phi
$$
with the nonlinearity 
\begin{align}
&\beta_\th (\phi^2) = |\phi|^{p-1} \phi \frac {f(\phi^2)}{\th+f(\phi^2)},
\label{eq:Nonl-beta}\\
& C_1 s^{q+1}\le |f(s^2)|\le C_2 s^{q+1}, \qquad 
s^2 |f'(s^2)|\le C_2 |f(s^2)|,\qquad p\in (1,\frac 4n),\,\, q>-1\nn
\end{align}
verify the monotonicity condition \eqref{eq:stabil}, $\la {L_{+}^\th}^{-1} \phi_\th , \phi_\th\ra >0$.
\begin{theorem}
\label{thm:stability}
Let $\phi_{\th}$ be ground states constructed in Proposition
\ref{prop:Contin} with the nonlinearities $\beta_\th$ satisfying \eqref{eq:Nonl-beta}.
Then then for a given nonlinear eigenvalue $\alpha\ne 0$ and 
all sufficiently small $\th$ the ground states 
$\phi_{\th}$ verify the monotonicity condition
\begin{equation}\label{eq:nstab}
\la {L_{+}^\th}^{-1} \phi_\th , \phi_\th\ra >0
\end{equation}
\end{theorem}
\begin{proof}
Condition \eqref{eq:nstab} is meaningful provided that 
$\phi_{\th}$ is orthogonal to the kernel of $L_{+}^{\th}$.
We start by examining the spectrum of the operator $L_{+}$. 
The operator $L^{+}$ has a unique negative eigenvalue,
the zero eigenvalue has multiplicity $n$ and the corresponding
eigenspace is spanned by the function $\frac\pa{\pa x_{i}}\phi$, \cite{W2}.
The rest of the spectrum is contained in the set $[\frac{\alpha^{2}}2,\infty)$.
Therefore, in the case $\alpha\ne 0$ the spectrum $\Sigma (L_{+})$ 
of $L_{+}$ has 
an isolated discrete component (in fact two components).
We can construct an eigenspace projector $P_{0}$ of an isolated
component of the discrete spectrum 
\be
\label{eq:proj}
P_{0}=\frac 1{2\pi i} \int_{\ga} (L_{+}-z)^{-1}\,dz
\ee
with an arbitrary curve $\ga$ encircling the desired spectral set
and such that $\ga\cap \Sigma (L_{+})=0$.
Consider now the resolvent of $L_{+}^{\th}$ at $z$ such that 
dist$(z,\Sigma (L_{+}))\ge C$ for some sufficiently small 
constant $C$, which only depends on $L_{+}$.
We have
\be
\label{eq:resolv}
(L_{+}^{\th}-z)^{-1} = (L_{+} -z)^{-1} -  
(L_{+}^{\th}-z)^{-1} (L_{+}^{\th} - L_{+}) (L_{+} -z)^{-1}
\ee
It is not difficult to show that for such $z$ 
$$
\|(L_{+} -z)^{-1} f\|_{H^{2}}\les \|f\|_{L^{2}}
$$
Therefore, using Corollary \ref{cor:Ldiffer} we can conclude from 
\eqref{eq:resolv} that for all sufficiently small $th\ge 0$
$$
\|(L_{+}^{\th}-z)^{-1}\|\le 2 \|(L_{+} -z)^{-1}\|
$$
and thus $z\not\in \Sigma (L_{+}^{\th})$. Moreover,
\be
\label{eq:Resclose}
\|(L_{+}^{\th}-z)^{-1} -(L_{+} -z)^{-1}\|\le c(\th)
\ee
for any $z:\,\,\,$ dist$(z,\Sigma (L_{+}))\ge C$. By Corollary 
\ref{cor:Ldiffer} the constant
$c(\th)\to 0$ as $\th\to 0$.
Therefore, for the same path $\ga$ as in \eqref{eq:proj} we can define
\be
\label{eq:projep}
P_{\th}=\frac 1{2\pi i} \int_{\ga} (L^{\th}_{+}-z)^{-1}\,dz
\ee
Moreover, for all sufficiently small $\th\ge 0$ the rank of
$P_{\th}$ remains constant.
Thus, for any sufficiently small $\th$ the operator $L_{+}^{\th}$
has a unique simple negative eigenvalue and a zero eigenspace of 
dimension $n$. Since we know that the functions 
$\frac\pa{\pa x_{i}}\phi_{\th}$ are contained in that subspace,
they, in fact, span it. Therefore, $\phi_{\th}$ is orthogonal to the
kernel of $L_{+}^{\th}$ and the expression \eqref{eq:nstab} is well
defined. 

For any sufficiently small $\th\ge 0$ we set $Q_{\th}$ 
to be a projection on the orthogonal complement
of the null eigenspace of $L_{+}^{\th}$.  Let  
$\lam \not\in\cup_{\ep} \Sigma (L_{+}^{\th})$. Define the operators
\be
\label{eq:Kep}
K_\th (\lam): = Q_{0}  (L_{+}^{\th}-\lam)^{-1} Q_{\th} - 
(L_{+}-\lam)^{-1} Q_{0}
\ee
It follows from \eqref{eq:Resclose} and the properties of 
the spectrum of $L_{+}$ that for all small $\th\ge 0$ and {\it {all}}
$\lam$ such that $|\lam|\le C$
\be
\label{eq:Lproj}
\|(L_{+}^{\th}-\lam)^{-1} Q_{\th}\|\le \frac 1{{\text{dist}} \,
\big (\lam, \Sigma (L_{+}^{\th})\setminus \{0\}\big )}\le C'
\ee
for some universal constant $C'$, determined by the 
operator $L_{+}$.
Also note that 
\be
\label{eq:Qclose}
\|Q_{\th} - Q_{0}\|\le c(\th)
\ee
This is a consequence of \eqref{eq:Resclose} and the definition 
$$
Q_{\th}= I - \frac 1{2\pi i} \int_{\ga} (L^{\th}_{+}-z')^{-1}\,dz'
$$
with a short path $\ga$ around the origin.
Using the resolvent identity
$$
(L_{+}^{\th}-z)^{-1} = (L_{+} -z)^{-1} +  
(L_{+}-z)^{-1} (L_{+}^{\th} - L_{+}) (L_{+}^{\th} -z)^{-1}
$$
we obtain that for any $\lam\not\in\cup_{\th} \Sigma (L_{+}^{\th})$
\begin{align*}
K_\th (\lam ) =& Q_{0} (L_{+}-\lam )^{-1} Q_{\th} - 
(L_{+}-\lam )^{-1} Q_{0} + Q_{0}  (L_{+}-\lam )^{-1} 
(L_{+}^{\th} - L_{+}) (L_{+}^{\th} -\lam )^{-1} Q_{\th}\\
=& (L_{+}-\lam )^{-1} Q_{0} (Q_{\th} - Q_{0}) + 
 (L_{+}-\lam )^{-1} Q_{0} 
(L_{+}^{\th} - L_{+}) (L_{+}^{\th} -\lam )^{-1} Q_{\th}
\end{align*}
Using Corollary \ref{cor:Ldiffer}, \eqref{eq:Lproj}, and 
\eqref{eq:Qclose} we infer that for any 
$\lam \le c$ and 
$\lam\not\in\cup_{\th} \Sigma (L_{+}^{\th})$
\be
\label{eq:Knorm}
\|K_{\th}(\lam)\| \le c(\th) 
(1+ \|(L_{+} -\lam )^{-1} Q_{0} (-\Laplace +1)\| )\le c'(\th)
\ee
uniformly in $\lam$. The last inequality follows since the 
operator norm of 
$(L_{+} -\lam )^{-1} Q_{0} (-\Laplace +1)$ is bounded by a universal
constant dependent on $L_{+}$ only. This can be seen as follows. 
Since $V_{0}$ is a smooth potential and $(L_{+} -\lam )^{-1} Q_{0}$
is bounded on $L^{2}$ we can replace the operator 
$(-\Laplace + 1)$ by $(L_{+}-\lam)$ and the result follows immediately.

We now test the operator $K_{\th}(\lam)$ on the ground state
$\phi_{\th}$. Using that $Q_\th \phi_\th = \phi_\th$ we obtain
$$
K_{\th}(\lam) \phi_{\th} = Q_{0} (L_{+}^{\th}-\lam)^{-1}
\phi_{\th} - (L_{+}-\lam)^{-1} \phi + (L_{+}-\lam)^{-1} Q_{0}
(\phi-\phi_{\th})
$$
Coupling the above identity with $\phi$.
\begin{align*}
\la \phi_{\th}, (L_{+}^{\th}-\lam)^{-1}\phi_{\th} \ra -
\la \phi , (L_{+}-\lam)^{-1}\phi \ra  & = 
\la\phi , K_{\th}(\lam) \phi_{\ep}\ra + 
\la (\phi_{\th}-\phi_{0}), (L_{+}^{\th}-\lam)^{-1}\phi_{\th} \ra
\\ &+ \la \phi, (L_{+}-\lam)^{-1} Q_{0}
(\phi-\phi_{\th})\ra = O(c(\th))
\end{align*}
where we have used that $Q_{\th}\phi_{\th} = \phi_{\th}$, the bound
\eqref{eq:Lproj}, and the estimate $\|\phi_{\th}-\phi\|_{H^{1}}\le c(\th)$, 
which follows from Proposition \ref{prop:Contin}. The above holds uniformly for 
all $|\lam |\le c$ and $\lam\not\in\cup_{\th} \Sigma (L_{+}^{\th})$.
Passing to the limit $\lam\to 0$, say from the upper half-plane, we obtain
that for all sufficiently small $\th\ge 0$
$$
\la \phi_{\th}, (L_{+}^{\th})^{-1}\phi_{\th} \ra =
\la \phi , L_{+}^{-1}\phi \ra  + O(c(\th)) < 0
$$
The last inequality follows since by the assumption 
$\phi$ is a stable ground state, i.e.,
$\la \phi , L_{+}^{-1}\phi \ra <0$. 
\end{proof}

\bibliographystyle{amsplain}

\noindent
\textsc{Rodnianski: Institute for Advanced Study and 
Department of Mathematics, Princeton University,  
Princeton N.J. 08544, U.S.A.}\\
{\em email: }\textsf{\bf irod@math.princeton.edu}

\medskip\noindent
\textsc{Schlag: Division of Astronomy, Mathematics, and Physics, 
253-37 Caltech, Pasadena, CA 91125, U.S.A.}\\
{\em email: }\textsf{\bf schlag@its.caltech.edu}

\medskip\noindent
\textsc{Soffer: Mathematics Department, Rutgers University, New 
Brunswick, N.J. 08903, U.S.A.}\\
{\em email: }\textsf{\bf soffer@math.rutgers.edu}

\end{document}